\documentclass[12pt,a4paper]{book}

\usepackage[intlimits,tbtags]{amsmath}
\usepackage{amsfonts, amssymb, amsthm} 
\usepackage{a4,array,cite,fancyhdr}

\makeatletter

\renewcommand{\chapter}{\clearpage\thispagestyle{plain}%
                    \global\@topnum\z@
                    \@afterindentfalse
                    \secdef\@chapter\@schapter}

\def\@chapter[#1]#2{\ifnum \c@secnumdepth >\m@ne
                       \if@mainmatter
                         \refstepcounter{chapter}%
                         \typeout{\@chapapp\space\thechapter.}%
                         \addcontentsline{toc}{chapter}%
                                   {\protect\mathversion{bold}
                                    \numberline{\thechapter}#1}%
                       \else
                         \addcontentsline{toc}{chapter}{#1}%
                       \fi
                    \else
                      \addcontentsline{toc}{chapter}{#1}%
                    \fi
                    \chaptermark{#1}%
                    \addtocontents{lof}{\protect\addvspace{10\p@}}%
                    \addtocontents{lot}{\protect\addvspace{10\p@}}%
                    \if@twocolumn
                      \@topnewpage[\@makechapterhead{#2}]%
                    \else
                      \@makechapterhead{#2}%
                      \@afterheading
                    \fi}

\def\@makechapterhead#1{%
  \vspace*{25\p@}%
  {\parindent \z@ \raggedright \reset@font
    \ifnum \c@secnumdepth >\m@ne
       \if@mainmatter
         \Large\bfseries \@chapapp{} \thechapter
         \par
         \vskip 10\p@
       \fi
       \fi
    \Large \bfseries\mathversion{bold} #1\par
    \nobreak
    \vskip 40\p@
  }}

\def\@makeschapterhead#1{%
  \vspace*{25\p@}%
  {\parindent \z@ \raggedright
    \reset@font
    \Large \bfseries\mathversion{bold}  #1\par
    \nobreak
    \vskip 40\p@
  }}

\renewcommand\tableofcontents{%
    \if@twocolumn
      \@restonecoltrue\onecolumn
    \else
      \@restonecolfalse
    \fi
    \clearpage\thispagestyle{plain}
      {\parindent \z@ \raggedright
        \reset@font
        \Large \bfseries\mathversion{bold}  \contentsname \par
        \@mkboth{\contentsname}{\contentsname}
        \nobreak
        \vskip 25\p@ 
      }
    \@starttoc{toc}%
    \if@restonecol\twocolumn\fi
    }

\renewenvironment{thebibliography}[1]
     {\chapter*{\bibname
        \@mkboth{\bibname}{\bibname}}%
       \addcontentsline{toc}{chapter}{\bibname}%
      \list{\@biblabel{\@arabic\c@enumiv}}%
           {\settowidth\labelwidth{\@biblabel{#1}}%
            \leftmargin\labelwidth
            \advance\leftmargin\labelsep
            \@openbib@code
            \usecounter{enumiv}%
            \let\p@enumiv\@empty
            \renewcommand\theenumiv{\@arabic\c@enumiv}}%
      \sloppy
      \clubpenalty4000
      \@clubpenalty \clubpenalty
      \widowpenalty4000%
      \sfcode`\.\@m}
     {\def\@noitemerr
       {\@latex@warning{Empty `thebibliography' environment}}%
      \endlist}

\renewcommand{\section}{\@startsection{section}{1}{\z@}%
  {-3.25ex\@plus -1ex \@minus -.2ex}%
  {1.5ex \@plus .2ex}%
  {\reset@font\large\bfseries\mathversion{bold}}}
\renewcommand{\subsection}{\@startsection{subsection}{2}{\z@}%
  {-3.25ex\@plus -1ex \@minus -.2ex}%
  {1.5ex \@plus .2ex}%
  {\reset@font\normalsize\bfseries\mathversion{bold}}}
\renewcommand{\subsubsection}{\@startsection{subsubsection}{3}{\z@}%
  {3.25ex \@plus1ex \@minus.2ex}%
  {-1em}%
  {\reset@font\normalsize\bfseries\mathversion{bold}}}

\renewcommand{\chaptermark}[1]{\markboth{\thechapter.\ #1}{}}

\fancypagestyle{plain}{%
 \fancyhf{}}

\makeatother

\newenvironment{squeezedsubequations}{%
  \begin{subequations}
    \\[\parskip]
    \setlength{\abovedisplayskip}{0ex plus0.5ex}
    \setlength{\belowdisplayskip}{0ex plus0.5ex}}{%
  \\[\parskip]\end{subequations}}

\newcommand{\I}{\mathrm{i}}

\newcommand{\tr}{\triangleright}
\newcommand{\tl}{\triangleleft}
\newcommand{\R}{\mathcal{R}}
\newcommand{\RI}{ \mathcal{R}_{\mathrm{I}} }
\newcommand{\RII}{ \mathcal{R}_{\mathrm{II}} }
\newcommand{\id}{\mathrm{id}}
\newcommand{\op}{\mathrm{op}}
\newcommand{\adL}{{\mathrm{ad_L}}}
\newcommand{\adR}{{\mathrm{ad_R}}}
\newcommand{\Proj}{\mathbb{P}}
\newcommand{\Func}{\mathcal{F}}

\newcommand{\Ket}[1]{\lvert #1\rangle} 
 
\newcommand{\rKet}[1]{\lVert #1\rangle} 
\newcommand{\rBra}[1]{\langle #1\rVert}
\newcommand{\Braket}[3]{\langle #1\rvert #2 \lvert #3\rangle}
\newcommand{\rBraket}[3]{\langle #1\rVert #2 \lVert #3\rangle}
\newcommand{\lrAngle}[1]{\langle #1\rangle}
\newcommand{\CGC}[6]{{\mathrm{C}_q(#1,#2,#3 \,|\, #4,#5,#6)}}
\newcommand{\RC}[6]{{\mathrm{R}_q(#1,#2,#3 \,|\, #4,#5,#6)}}
\newcommand{\BolliO}[6]{{B_q^0(#1,#2,#3 \,|\, #4,#5,#6)}}
\newcommand{\Bollil}[6]{{B_q^1(#1,#2,#3 \,|\, #4,#5,#6)}}
\newcommand{\Bas}[2]{E^{#1}_{#2}}
\newcommand{\Base}[1]{E_{#1}}

\newtheorem{Definition}{Definition}
\newtheorem{Theorem}{Theorem}
\newtheorem{Proposition}{Proposition}

\newcommand{\uq}{{\mathcal{U}_q(\mathrm{u}_1)}}

\newcommand{\su}{{\mathcal{U}(\mathrm{su}_2)}}
\newcommand{\slq}{{\mathcal{U}_q(\mathrm{sl}_2)}}

\newcommand{\suq}{{\mathcal{U}_q(\mathrm{su}_2)}}
\newcommand{\slC}{{\mathcal{U}_q(\mathrm{sl}_2(\mathbb{C})) }}
\newcommand{\SLq}{{SL_q(2)}}

\newcommand{\SUq}{{SU_q(2)}}
\newcommand{\SLC}{{SL_q(2,\mathbb{C})}}
\newcommand{\Mink}{{\mathcal{M}_q}}
\newcommand{\Poin}{{\mathcal{P}_q}}
\newcommand{\Euk}{{\mathcal{E}_q}}
\newcommand{\qsphere}{{\mathcal{S}_{q\infty}^\op}}

\newcommand{\Heli}{Z}

\newcommand{\qi}{\xi}

\begin{document}

\thispagestyle{empty}

\begin{center}
  {\LARGE\bfseries\mathversion{bold} Spin Representations of the\\
    $q$-Poincar{\'e} Algebra}
  \\[1cm]
  {Dissertation der Fakult\"at f\"ur Physik der\\
    Ludwig-Maximilians-Universit\"at M\"unchen}
  \\[0.5cm]
  vorgelegt von {\\ \large\bfseries Christian
    Blohmann\footnote{electronic mail:
            Christian.Blohmann{@}Physik.Uni-Muenchen.de}
    \\ } aus Villingen
  im Schwarzwald
  \\[0.5cm]
  M\"unchen, den 8. M\"arz 2001
  \\[1cm]
  Sektion Physik, Universit\"at M\"unchen\\
  Theresienstr. 37, D-80333 M\"unchen
  \\[2ex]
  Max-Planck-Institut f\"ur Physik\\
  F\"ohringer Ring 6, D-80805 M\"unchen
  \\[1cm]
\end{center}

{\centering {\bfseries Abstract}\\[1ex]}
The spin of particles on a non-commutative geometry is investigated
within the framework of the representation theory of the $q$-deformed
Poincar{\'e} algebra. An overview of the $q$-Lorentz algebra is given,
including its representation theory with explicit formulas for the
$q$-Clebsch-Gordan coefficients. The vectorial form of the $q$-Lorentz
algebra (Wess), the quantum double form (Woronowicz), and the dual of
the $q$-Lorentz group (Majid) are shown to be essentially isomorphic.
The construction of $q$-Minkowski space and the $q$-Poincar{\'e}
algebra is reviewed. The $q$-Euclidean sub-algebra, generated by
rotations and translations, is studied in detail. The results allow
for the construction of the $q$-Pauli-Lubanski vector, which, in turn,
is used to determine the $q$-spin Casimir and the $q$-little algebras
for both the massive and the massless case. Irreducible spin
representations of the $q$-Poincar{\'e} algebra are constructed in an
angular momentum basis, accessible to physical interpretation. It is
shown how representations can be constructed, alternatively, by the
method of induction. Reducible representations by $q$-Lorentz spinor
wave functions are considered. Wave equations on these spaces are
found, demanding that the spaces of solutions reproduce the
irreducible representations. As generic examples the $q$-Dirac
equation and the $q$-Maxwell equations are computed explicitly and
their uniqueness is shown.

\tableofcontents

\chapter*{Introduction}
\markboth{Introduction}{Introduction}
\addcontentsline{toc}{chapter}{Introduction}

\subsubsection{Motivation}
From the beginnings of quantum field theory it has been argued that
the pathological ultraviolet divergences should be remedied by
limiting the precision of position measurements by a fundamental
length \cite{Born:1933,Born:1934,March:1936,Heisenberg:1938}. In view
of how position-momentum uncertainty enters into quantum mechanics, a
natural way to integrate such a position uncertainty in quantum theory
would have been to replace the commutative algebra of space
observables with a non-commutative one \cite{Snyder:1947}. However,
deforming the space alone will in general break the symmetry of
spacetime.

In order to preserve a background symmetry the symmetry group must be
deformed together with the space it acts on. It is clear that Lie
groups cannot be continuously deformed within their proper category:
From the classification of semi-simple Lie groups we know that they
form a countable and hence discrete set. Being manifolds, however,
they can be naturally embedded in the category of algebras by the
Gelfand-Neumark map \cite{Gelfand:1943}, the additional group
structure on the manifold side translating into a Hopf structure on
the algebra side. But although Hopf algebras in general had been
familiar to mathematicians for some time
\cite{Hopf:1941,Abe,Sweedler}, hardly any non-trivial examples of Hopf
algebras were known \cite{Pareigis:1981}. This situation changed with
the discovery of quantum groups \cite{Drinfeld:1986}, that is, with
the discovery of generic methods to continuously deform Lie algebras
\cite{Drinfeld:1985,Jimbo:1985} and matrix groups
\cite{Woronowicz:1987,Faddeev:1990,Takeuchi:1990} within the category
of Hopf algebras.

Quantum groups now provided a consistent mathematical framework to
formulate physical theories on non-commutative spaces. Beginning with
the non-commutative plane \cite{Manin:1988}, $q$-deformations of a
variety of objects have since been constructed: differential calculi
on non-commutative spaces \cite{Wess:1991}, Euclidean space
\cite{Faddeev:1990}, Minkowski-Space \cite{Carow-Watamura:1990}, the
Lorentz group and the Lorentz algebra
\cite{Podles:1990,Carow-Watamura:1991,Schmidke:1991,Ogievetskii:1991a},
the Poincar\'e algebra \cite{Ogievetskii:1992a}, to name a few. The
study of these objects has produced interesting results. For example,
it has been found that free theories on non-commutative spaces can be
viewed as theories on ordinary commutative spaces with complicated
interactions \cite{Drinfeld:1989,Fichtmuller:1996,Madore:2000}.

Another result is, that $q$-deformation will in general discretize the
spectra of spacetime observables
\cite{Schwenk:1992,Wess:1997,Cerchiai:1998}, that is, $q$-deformation
puts physics on a spacetime lattice.  This nourishes the hope that
$q$-deformed field theories might be regularized, one of the original
motivations to consider non-commutative geometries. It is not new that
deformation is a method to regularize field theory --- at least, it is
one way to look at the first step of renormalization: In the loop
expansion of transition probabilities some terms turn out to be
infinite, so we regularize them by a sort of deformation process in
order to classify the divergences. In this sense $q$-deformation can
be viewed as attempt to shift the deformation from the end of the
construction of field theory (perturbative expansion) to the beginning
(symmetry structures).

\subsubsection{Aims} 
Given the $q$-Poincar{\'e} algebra as background symmetry, how can we
construct a quantum theory upon it?  If states continue to be
described by vectors of a Hilbert space, it must be specified how the
$q$-Poincar{\'e} algebra acts on them, that is, we must construct
representations of the $q$-Poincare algebra. If we further want to
describe \emph{elementary} particles the representations must be
\emph{irreducible} \cite{Wigner:1939}. If we want to use reducible
representations such as the Dirac spinor representation, we need
additional constraints to eliminate the redundant degrees of freedom.
These constraints are the wave equations. The interpretation of this
quantum theory forces us to consider multi-particle states. These are
described as properly symmetrized (or anti-symmetrized) tensor product
representations.  Symmetrization or anti-symmetrization means that we
need a ray representation of the permutation group on the tensor
product space which is compatible with (intertwines with) the action
of the $q$-Poincar{\'e} algebra. The physical states are the orbits of
this action of the permutation group, while the direct sum of all such
multi-particle spaces is the Fock space. We summarize:
\begin{enumerate}
\item[(i)] Elementary $q$-particles are irreducible representations of the
  $q$-Poincar{\'e} algebra.\label{pg:Programi}
\item[(ii)] Wave equations are the constraints to eliminate the redundant
  degrees of freedom of a reducible representation.
\item[(iii)] $q$-Fields are symmetrized or anti-symmetrized multi-particle
  states.
\end{enumerate}
In the undeformed case these principles completely determine the free
relativistic quantum theory. Therefore, it is reasonable to use them
as program to construct the deformed theory. 

This program has been pursued in previous work
\cite{Cerchiai:1998,Cerchiai,Ogievetskii:1992b,Pillin:1993,Pillin,Schirrmacher:1992,Pillin:1994b,Song:1992,Meyer:1995,Podles:1996}.
In \cite{Ogievetskii:1992b,Pillin:1993,Pillin,Cerchiai,Cerchiai:1998}
irreducible spin zero representations of the $q$-Poincar{\'e} algebra were
constructed. While in \cite{Cerchiai,Cerchiai:1998} the realization of the
$q$-Poincar{\'e} algebra within $q$-Minkowski phase space was considered,
such that the representations were naturally limited to orbital
angular momentum,\footnote{See Eq.~\eqref{eq:Orbital}.}  it is
possible to extend \cite{Pillin:1993,Pillin} to include spin
representations (Sec.~\ref{sec:SpinReps1}).  Various methods to
construct wave equations have been proposed, based on $q$-Clifford
algebras \cite{Schirrmacher:1992}, $q$-deformed co-spinors
\cite{Pillin:1994b}, or differential calculi on quantum spaces
\cite{Song:1992,Meyer:1995,Podles:1996}, leading to mutually different
results. This is unsatisfactory since the construction of wave
equations according to (ii) should determine the wave equations
uniquely as in the undeformed case \cite{BarutRaczka} and should not
demand any additional mathematical structure besides the
$q$-Poincar\'e algebra and the basic apparatus of quantum mechanics.

The aim of the present work is to investigate the nature of spin
within the representation theory of the $q$-Poincar\'e algebra.

\subsubsection{Results}
Our main results are:
\begin{itemize}
\item The $q$-deformed Pauli-Lubanski vector is computed
  (Sec.~\ref{sec:MainContrib2a}), from which the spin Casimir and the
  little algebras can be determined (Sec.~\ref{sec:MainContrib2b}).
\item Irreducible representations with spin are constructed
  (Chap.~\ref{sec:SpinReps}).
\item A practical method to uniquely compute the wave equations is
  developed (Sec.~\ref{sec:MainContrib4a}). As examples the $q$-Dirac
  equation (Sec.~\ref{sec:MainContrib4b}) and the $q$-Maxwell
  equations (Sec.~\ref{sec:MainContrib4c}) are computed.
\end{itemize}
To give a more detailed overview:
 
In chapter~\ref{sec:LorentzConstruction} we review the construction of
the $q$-Lorentz algebra.  We start with the quantum plane, $xy=qyx$,
derive the algebra of coacting quantum matrices $M_q(2)$, introduce
the $q$-spinor metric, the quantum special linear group $\SLq$ and its
real form $\SUq$. We introduce dotted spinors, join an undotted and a
dotted corepresentation to form the quantum Lorentz group $\SLC$.
Using the duality between $\SLq$ and $\slq$ we compute the quantum
Lorentz algebra $\slC$ by dualizing $\SLC$
\cite{Manin:1988,Carow-Watamura:1990,Carow-Watamura:1991,Majid:1993}.
The presentation puts emphasis on the fact that in the construction of
the $q$-Lorentz algebra as it is understood now, hardly any
arbitrariness is involved.

Chapter~\ref{sec:LorentzStructure} explores the structure of the
$q$-Lorentz algebra. The representation theory of the $q$-Lorentz
algebra is reviewed, explicit formulas for the
$q$-Clebsch-Gordan coefficients are given. After a general
consideration of the different sorts of tensor operators, the
vectorial generators of $\slq$ are determined. Three different forms
of the $q$-Lorentz algebra are related by explicit formulas: the dual
of the $q$-Lorentz group \cite{Majid:1993}, the quantum double form
\cite{Podles:1990}, and the vectorial or $RS$-form
\cite{Schmidke:1991,Lorek:1997a,Rohregger:1999}. The isomorphism
between the quantum double form and the vectorial form that is found
(Sec.~\ref{sec:MainContrib1}) relates the work of the Warsaw and the
Munich group.

In chapter~\ref{sec:PoincareStructure} the results of
chapter~\ref{sec:LorentzStructure} are used to construct the algebra
of $q$-Minkowski space \cite{Carow-Watamura:1990}. Commutation
relations of the generators of different forms of the $q$-Lorentz
algebra with the spacetime generators are given. We study the
structure of the $q$-Euclidean algebra consisting of rotations and
translations in order to find a good zero component of the
$q$-Pauli-Lubanski vector. A technique of boosting is used to
calculate the other components (Sec.~\ref{sec:MainContrib2a}). The
$q$-Pauli-Lubanski vector is used to compute the little algebras for
the massive and the massless case (Sec.~\ref{sec:MainContrib2b}).

Chapter~\ref{sec:SpinReps} contains the construction of massive spin
representations of the $q$-Poincar{\'e} algebra. In the first part we
construct irreducible representations in an angular momentum basis,
which is accessible to physical interpretation. The calculations are
considerably simplified by the $q$-Wigner-Eckart theorem.  In the
second part we briefly show how representations of the $q$-Poincar{\'e}
algebra can be constructed using the method of induced
representations.

In chapter~\ref{sec:WaveEqs} we calculate free wave equations. We
start with the representation theoretic interpretation of free wave
equations. Then we consider the generalities of $q$-Lorentz spinor
representations, conjugate spinors, and the relation between momenta
and derivations. Finally, we put things together and uniquely
determine the $q$-Dirac equation including $q$-gamma matrices and
their $q$-Clifford algebra, the $q$-Weyl equations, and the
$q$-Maxwell equations.

\subsubsection{Outlook}
While our approach to the $q$-Poincar{\'e} algebra was representation
theoretic, the problems we had to overcome were mostly on the
algebraic side: A method to boost vector operators, complete sets of
commuting observables, the spin Casimir, the spin symmetry algebras,
spinor conjugation --- all this had to be found before spin
representations and spinorial wave equations could be computed. Now,
that the algebraic tool set is more complete, we are prepared for the
next steps towards a $q$-deformed relativistic quantum theory.

One promising way to continue this work would be to couple the
$q$-Dirac and the $q$-Maxwell field, for which the mathematical
setting has been provided in chapter~\ref{sec:WaveEqs}.

\subsubsection{Notation}
Throughout this work, the deformation parameter $q$ is assumed real,
$q>1$. We frequently use the abbreviations
\begin{xalignat}{2}
  \lambda &:= q - q^{-1} \,,&
  [j] &:= \frac{q^j-q^{-j}}{q-q^{-1}} \,,
\end{xalignat}
where $j$ is a number. In particular, we have $[2] = q+q^{-1}$. Spinor
indices running through $\{-,+\} = \{-\tfrac{1}{2},+\tfrac{1}{2}\}$
are denoted by lower case Roman letters ($a$, $b$, $c$, $d$),
$3$-vector indices running through $\{-,3,+\} = \{-1,0,+1\}$ by upper
case Roman letters ($A$, $B$, $C$), and 4-vector indices running
through $\{0,-,+,3\}$ by lower case Greek letters ($\mu$, $\nu$,
$\sigma$, $\tau$). Quantum Lie groups are written with a subscript $q$
like $\SLq$, quantum enveloping algebras like $\slq$.

\chapter{Construction of the $q$-Lorentz Algebra}
\label{sec:LorentzConstruction}

In undeformed quantum mechanics we can represent a state by a wave
function $\psi:\mathbb{R}^n \rightarrow \mathbb{C}$.  In this
representation, the observables $x_i$, which describe the measurement
of the position of the particle, act on $\psi$ by multiplications with
the functions $x_i:\mathbb{R}^n \rightarrow \mathbb{C}$, $x_i(\vec{r})
= r_i$.  In this sense, geometry is described by the algebra of
functions over a space, $\Func(\mathbb{R}^n)$, rather than by the
space $\mathbb{R}^n$ itself. Replacing a space by its function
algebra, it is natural to replace an endomorphism $f$ by its pullback
$f^*$,
\begin{equation}
  \mathbb{R}^n \stackrel{f}{\longrightarrow} \mathbb{R}^n
  \quad\Rightarrow\quad
  \Func(\mathbb{R}^n) \stackrel{f^*}{\longleftarrow} \Func(\mathbb{R}^n)
  \,, \quad\text{where}\quad (f^* x_i)(\vec{r}) := x_i (f\vec{r}) \,,
\end{equation}
yielding a recipe to translate spaces and homomorphisms of spaces into
algebras and homomorphisms of algebras. In the language of category
theory $\Func$ is called a cofunctor \cite{MacLane}, the prefix ``co''
reminding us that we have to reverse arrows.

For a consistent mathematical framework we must extend this method of
algebraization to any additional structure on $\mathbb{R}^n$. If there
is for example the action $\phi$ of a group $G$ on the space we get
\begin{equation}
  G \otimes \mathbb{R}^n \stackrel{\phi}{\longrightarrow} \mathbb{R}^n
  \quad\Rightarrow\quad
  \Func(G) \otimes \Func(\mathbb{R}^n)
  \stackrel{\phi^*}{\longleftarrow} \Func(\mathbb{R}^n) \,,
\end{equation}
where $\Func(G)$ is the algebra of functions over the group and the
homomorphism of algebras $\rho := \phi^*$ is called the coaction. The
structure maps of the group, multiplication $\mu$, unit $\eta$, and
inverse, translate into comultiplication $\Delta = \mu^*$, counit
$\varepsilon = \eta^*$, and coinverse or antipode $S$. The group
axioms translate into axioms of this co-structure \cite{Hopf:1941}. An
algebra equipped with this co-structure is called a Hopf algebra
\cite{Abe,Sweedler}.

So far, the structure of spaces and groups acting on them has only
been rephrased in a more algebraic but equivalent language. But unlike
the category of Lie groups, the category of algebras allows for
continuous deformation: We can replace the trivial commutation
relations of the algebra of space functions by non-trivial ones, which
depend on a real parameter $q$.  This $q$-deformation of the space
algebra forces us to $q$-deform any Hopf algebra coacting on it, as
well. Reminiscent of their relation to quantum theory, these deformed
algebras are called quantum spaces and quantum groups. Instead of
quantum groups we can consider their Hopf duals
\cite{Takeuchi:1977,vanDaele:1993}, the quantum algebras, which are
deformations of the enveloping Lie algebras. Since quantum algebras
have a familiar undeformed counterpart, they become directly
accessible to physical interpretation. For example, the generators of
the quantum algebra of rotations are the $q$-deformed angular momentum
operators.

\section{$q$-Spinors and $SU_q(2)$}

\subsection{$q$-Spinors and Their Cotransformations}

The simplest quantum space is the deformation of the algebra
$\Func(\mathbb{C}^2)= \mathbb{C}[x,y]$ of polynomial spinor functions.
We replace the trivial commutation relations $xy=yx$ with $xy=qyx$,
where $q$ is a real parameter $q>1$, and call the resulting algebra
\begin{equation}
 \mathbb{C}_q^2 := \mathbb{C}\langle x, y \rangle /
 \langle xy=qyx \rangle
\end{equation}
the algebra of $q$-spinors or the quantum plane \cite{Manin:1988}.

As in the undeformed case we want the spinor algebra to carry a left and
a right matrix corepresentation. We define the vector of spinor
generators 
\begin{equation}
  \psi_a = (\psi_-, \psi_+ ) := (x, y)
\end{equation}
a matrix of generators of the algebra $M_q(2)$ of 2$\times$2-matrices
\begin{equation}
  M^a{}_b = \begin{pmatrix} a & b \\ c & d \end{pmatrix}
\end{equation}
with respect to the indices $\{-,+\}=\{-\tfrac{1}{2},+\tfrac{1}{2}\}$
and the left and right coaction of this matrix on the spinor
\begin{xalignat}{2}
  \rho_\mathrm{L}(\psi_a) &:= M^a{}_{a'} \otimes \psi_{a'}  \,,&
  \rho_\mathrm{R}(\psi_a) &:= \psi_{a'} \otimes M^{a'}{}_a \,,
\end{xalignat}
where we sum over repeated indices, and where the coproduct of
$M_q(2)$ is defined by $\Delta(M^a{}_{c})=M^a{}_{b}\otimes M^b{}_{c}$.

We want the deformed commutation relations between the generators of
$M_q(2)$ to be consistent with those of the $q$-spinor, $xy=qyx$, that
is, the coaction maps must be algebra homomorphisms. This uniquely
determines the relations
\begin{equation}
\label{eq:MqRel}
\begin{gathered}
  ab=qba\,, \quad ac=qca\,, \quad bd = q db\,, \quad cd = q dc \\
  bc = cb\,, \quad ad - da = (q - q^{-1}) bc \,.
\end{gathered}
\end{equation}
The algebra freely generated by $a, b, c, d$ modulo these
relations~\eqref{eq:MqRel} is called $M_q(2)$ the algebra of
2$\times$2 quantum matrices. Introducing the $R$-Matrix
\begin{equation}
   R^{ab}{}_{cd} =
  \begin{pmatrix}
    q & 0 & 0 & 0 \\
    0 & 1 & 0 & 0 \\
    0 & q-q^{-1} & 1 & 0 \\
    0 & 0 & 0 & q 
  \end{pmatrix}
\end{equation}
with respect to the indices $\{--,-+,+-,++\}$, Eqs.~\eqref{eq:MqRel}
can be written in the compact form
\begin{equation}
\label{eq:FRT-Rel}
  R^{ab}{}_{c'd'} M^{c'}{}_{c}
  M^{d'}{}_{d} = M^{b}{}_{b'} M^{a}{}_{a'}
  R^{a'b'}{}_{cd} \,,
\end{equation}
the famous FRT-relations of matrix quantum groups \cite{Faddeev:1990}. 

\subsection{The $q$-Spinor Metric and $SL_q(2)$}

With the spinor metric
\begin{equation}
\label{eq:SpinorMetric}
  \varepsilon_{ab} = - \varepsilon^{ab} =
  \begin{pmatrix} 0 & q^{-1/2} \\ -q^{1/2} & 0 \end{pmatrix}
  ,\quad\text{with}\quad
  \varepsilon_{ab} \varepsilon^{bc}
  = \delta_a^c
\end{equation}
we can write $xy = qyx$ as $ \psi_a \psi_b
\varepsilon^{ab} = 0$. In analogy to the undeformed case the
spinor metric must thus be invariant under $M_q(2)$-transformations up
to a factor. Indeed, we find
\begin{equation}
\label{eq:InvariantMetric}
  M^{a}{}_{a'}
  M^{b}{}_{b'} \varepsilon^{a'b'} =
  (\mathrm{det}_q M)\, \varepsilon^{ab} \,,
\end{equation} 
where $\mathrm{det}_q M= ad-qbc$ is central in $M_q(2)$.  Constraining the
transformations to leave the scalar product $\psi_a \phi_b
\varepsilon^{ab}$ of two spinors strictly invariant, we
obtain $SL_q(2) := M_q(2)/ \langle \mathrm{det}_q M = 1 \rangle$ the
deformation of the function algebra of the group of special linear
transformations.

Finally, Eq.~\eqref{eq:InvariantMetric} can be contracted with the
metric from the right. From the resulting equation
\begin{equation}
  M^{a}{}_{a'} (M^{b}{}_{b'}
  \varepsilon^{a'b'} \varepsilon_{bc}) =
  \delta^{a}_c
\end{equation}
we can read off the antipode 
\begin{equation}
  S(M^{a}{}_{b}) := \varepsilon^{aa'}
  M^{b'}{}_{a'} \varepsilon_{b'b} =
  \begin{pmatrix} d & - q^{-1} b \\ -q c & a \end{pmatrix} \,,
\end{equation}  
playing the role of the inverse, $(M^{-1})^{a}{}_{b} =
S(M^{a}{}_{b})$. This completes the Hopf algebra structure of
$SL_q(2)$.

\subsection{Upper Spinor Indices, Conjugation, and $SU_q(2)$}

Defining a transposition on $SL_q(2)$ by
\begin{equation}
  T(M^{a}{}_{b}) = (M^T)^{a}{}_{b} :=
  M^{b}{}_{a} \,,
\end{equation}
we can consider now a spinor transforming under the congredient
representation $(M^T)^{-1}$. As in the undeformed case we indicate
this transformation property by an upper index.
\begin{equation}
  \rho_\mathrm{R}(\psi^a) = \psi^b \otimes
    ((M^T)^{-1})^b{}_a =
    \psi^b \otimes S(M^a{}_b)
    = \psi^b \otimes \varepsilon^{aa'}
  M^{b'}{}_{a'} \varepsilon_{b'b}
\end{equation}
Contracting this equation with the spinor metric we find
\begin{equation}
  \rho_\mathrm{R}(\varepsilon_{aa'} \psi^{a'})
  = (\varepsilon_{bb'} \psi^{b'}) \otimes
  M^{b}{}_{a} \,,
\end{equation}
telling us that $\varepsilon_{aa'} \psi^{a'}$
transforms as a spinor with lower index. We conclude that we can
raise and lower indices by
\begin{xalignat}{2}
  \psi^a &= \varepsilon^{aa'} \psi_{a'} \,,&
  \psi_a &= \varepsilon_{aa'} \psi^{a'} \,.
\end{xalignat}

When we rewrite the spinor commutation relations as
\begin{equation}
  0= \psi_a \psi_b \varepsilon^{ab} =
  \varepsilon_{aa'}\varepsilon_{bb'}
  \psi^{a'} \psi^{b'} \varepsilon^{ab} =
  \psi^a \psi^b \varepsilon_{ba} =
  - \psi^a \psi^b \varepsilon^{ba} \,,
\end{equation}
we see that a spinor with upper index satisfies commutation relations
opposite to a spinor with lower index. Thus, we can define a
$*$-structure on the spinor algebra $\mathbb{C}_q^2$ by
$(\psi_a)^* := \psi^a$. This induces a $*$-structure on
$SL_q (2)$ as well, by demanding the stars to be compliant with the
coaction, $\rho_\mathrm{R} \circ * = (*\otimes *)\circ
\rho_\mathrm{R}$.

A stared spinor transforms as a spinor with upper index, that is, by
the congredient representation. We conclude that the induced
$*$-operation on $SL_q(2)$ is given by
\begin{equation}
  (M^{a}{}_{b})^* = S(M^{b}{}_{a}) \,.
\end{equation}
In other words, we have $(M^T)^* = M^{-1}$, which can be viewed as a
quantum group analogue of a unitarity condition. Therefore, $SL_q(2)$
with this $*$-structure is called $SU_q(2)$.

\section{The $q$-Lorentz Group}

\subsection{Dotted Spinors} 

We want to construct a deformation of the Lorentz group
$SL(2,\mathbb{C})$, which is, viewed as real manifold, 6-dimensional,
having 6 independent infinitesimal generators. Now, a spinor and its
complex conjugate and thus the corepresentation matrix and its
conjugate are no longer linearly dependent. This means that we have to
add the conjugates $\bar{M}^{a}{}_{b} := (M^{a}{}_{b})^*$ and
$\bar{\psi}_a := (\psi_a)^*$ as extra generators. Of course, the
conjugate spinor cotransforms under the conjugate matrix. As in the
undeformed case, we will indicate that a quantity transforms like a
conjugate spinor by a dotted index.  Thus, writing $\psi_{\dot{a}}$
implies
\begin{equation}
  \rho_\mathrm{R}(\psi_{\dot{a}})
  = \psi_{\dot{b}} \otimes \bar{M}^{b}{}_{a} \,,
\end{equation}
where we think of the dot as belonging to $\psi$ rather than to the
index itself. Since the $*$-operation is by definition an algebra
anti-homomorphism (and a coalgebra homomorphism), the conjugate
generators satisfy the opposite commutation relations of their
pre-images. However, it is more convenient to combine the conjugate
generators $\bar{M}$ linearly to form another matrix $M_2$ defined
implicitly by $(M_2^T)^{-1} := \bar{M}$, that is,
\begin{equation}
  S(M_2^{b}{}_{a}) := \bar{M}^{a}{}_{b} \,.
\end{equation}
$S\circ T$ is an algebra anti-homomorphism (and a coalgebra
homomorphism), so $M_2$ naturally generates a $SL_q(2)$ Hopf algebra.
We now have two sets of generators generating two copies of $SL_q(2)$.
For a consistent notation we will subscript the first set $M = M_1$ as
well. The $*$-operation can then be written as
\begin{equation}
  (M_1^a{}_b)^{*_{ SL_q(2,\mathbb{C}) }} =
  (M_2^{a}{}_{b})^{*_{SU_q(2)}} \,.
\end{equation}

Finally, we introduce upper dotted indices by demanding them to
transform according to 
\begin{equation}
  \rho_\mathrm{R}(\psi^{\dot{a}})
  = \psi^{\dot{b}} \otimes M_2^{b}{}_{a} \,.
\end{equation}
This leads to formulas for raising and lowering dotted indices
\begin{xalignat}{2}
  \psi^{\dot{a}} &= \psi_{\dot{b}}
    \varepsilon^{ba}\,, &
  \psi_{\dot{a}} &= \psi^{\dot{b}}
    \varepsilon_{ba} \,.
\end{xalignat}

\subsection{Commutation Relations of the $q$-Lorentz Group}
\label{sec:LorentzCommutation}

So far, we know that the $q$-Lorentz group must be generated by two
copies of $SL_q(2)$, generated by two sets of generators
$M_1^{a}{}_{b}$ and $M_2^{a}{}_{b}$, respectively.  The only thing we
do not know yet are the commutation relations between $M_1$ and $M_2$.
A priori, there are several choices of commutation relations, from
which we will select one by an additional requirement: We will demand
$SL_q(2,\mathbb{C})$ to possess a substructure of rotational symmetry,
that is, we are looking for a homomorphism of Hopf-$*$
algebras\footnote{During the transition from groups to quantum groups
  the arrows of mappings have to be reversed.} $\mu :
SL_q(2,\mathbb{C}) \rightarrow SU_q(2)$.

Embedding the generators by $M_1^{a}{}_{b} \hookrightarrow
M^{a}{}_{b} \otimes 1$ and $M_2^{a}{}_{b}
\hookrightarrow 1\otimes M^{a}{}_{b}$ in a tensor product of
two $SL_q(2)$, the multiplication map
\begin{equation}
  \mu : SL_q(2)\otimes SL_q(2) \rightarrow SU_q(2)  
\end{equation}
is the obvious choice. Note, that according to the preceding section
$g\otimes h \in SL_q(2)\otimes SL_q(2)$ is to be equipped with the
$*$-structure $(g\otimes h)^* := h^* \otimes g^*$, such that
$\mu((g\otimes h)^*) = \mu(h^* \otimes g^*) = h^* g^* = (gh)^* =
\mu(g\otimes h)^*$. In other words, $\mu$ is already compliant with
the $*$-structures.

For $\mu$ to be a homomorphism of algebras, the images of the
generators, $\mu(M_1^{a}{}_{b}) =
\mu(M_2^{a}{}_{b})= M^{a}{}_{b}$, must satisfy the
$SL_q(2)$ commutation relations~\eqref{eq:FRT-Rel}. This means
that the generators have to satisfy
\begin{equation}
  R^{ab}{}_{c'd'} M_2^{c'}{}_{c}
  M_1^{d'}{}_{d} = M_1^{b}{}_{b'}
  M_2^{a}{}_{a'} R^{a'b'}{}_{cd} \,,
\end{equation}
which completes the algebraic structure of the $q$-Lorentz
group.\footnote{If we drop the requirement of a subsymmetry of
  rotations, we can construct an alternative $q$-Lorentz group with
  two \emph{commuting} copies of $\SLq$. It turns out to be
  unphysical, however, insofar as the according $q$-Poincar\'e algebra
  has no mass Casimir.}

To summarize, let us give a compact and rigorous definition of the
$q$-Lorentz group \cite{Podles:1990,Carow-Watamura:1991}. First we
need to define the so-called coquasitriangular map $R:SL_q(2)\otimes
SL_q(2) \rightarrow \mathbb{C}$ on the generators by
\begin{equation}
\label{eq:R-Map}
  R(M^{a}{}_{c},M^{b}{}_{d}):= q^{-\frac{1}{2}}R^{ab}{}_{cd} \,,
\end{equation}
which can be shown to extend to all of $SL_q(2)$ by linearity in both
arguments and by demanding
\begin{xalignat}{2}
  R(fg,h) &:= R(f,h_{(1)})R(g,h_{(2)}) \,, &
  R(f,gh) &:= R(f_{(1)},h)R(f_{(2)}, g) \,.
\end{xalignat}
The factor $q^{-\frac{1}{2}}$ has been introduced for convenience. The
map $R$ has a unique convolution inverse, that is, a map
$R^{-1}:SL_q(2)\otimes SL_q(2) \rightarrow \mathbb{C}$ with
\begin{equation}
  R(a_{(1)}, b_{(1)}) R^{-1}(a_{(2)}, b_{(2)}) =
  R^{-1}(a_{(1)},b_{(1)}) R(a_{(2)}, b_{(2)}) =
  \varepsilon(a)\varepsilon(b) \,,
\end{equation}
simply defined by
\begin{equation}
  R^{-1}(M^{a}{}_{c},M^{b}{}_{d}):= q^{\frac{1}{2}}(R^{-1})^{ab}{}_{cd} \,.
\end{equation}
Using $R$, the commutation relations of $SL_q(2)$ can be written as
\begin{equation}
  R(a_{(1)},b_{(1)}) a_{(2)} b_{(2)} =
  b_{(1)} a_{(1)} R(a_{(2)}, b_{(2)}) \,.
\end{equation}

\begin{Definition}
  Let $R$ denote the coquasitriangular map of $SL_q(2)$ and $R^{-1}$
  its convolution inverse. The vector space $SL_q(2)\otimes SL_q(2)$
  with tensor product coalgebra structure, $\Delta(g\otimes h) =
  (g_{(1)}\otimes h_{(1)})\otimes (g_{(2)}\otimes h_{(2)})$,
  $\varepsilon(g\otimes h) = \varepsilon(g) \varepsilon(h)$, with
  multiplication
  \begin{equation}
    (g\otimes h)(g'\otimes h') =  gg'_{(2)} \otimes h_{(2)}h'
    \,R^{-1}(h_{(1)},g'_{(1)}) R(h_{(3)},g'_{(3)})
  \end{equation}
  antipode $S(g\otimes h) = (1 \otimes S(h))(S(g)\otimes 1)$, and
  $*$-structure
  \begin{equation}
    (g \otimes h)^{*_{ SL_q(2,\mathbb{C}) }} =
    h^{*_{SU_q(2)}} \otimes g^{*_{SU_q(2)}}
  \end{equation}
  is the $q$-Lorentz group $SL_q(2,\mathbb{C})$.
\end{Definition}

\section{The $q$-Lorentz Algebra as Dual of the $q$-Lorentz Group}

For a symmetry of a quantum mechanical system the mathematical object
with a direct physical interpretation is the enveloping algebra of
the symmetry group's Lie algebra rather than the group itself. The
Hilbert space representations of its generators are the observables of
the conserved quantities corresponding to the symmetry. Consequently,
rather than in the quantum group itself we are interested in its dual,
the quantum enveloping algebra.

\subsection{$\suq$ as dual of $SU_q(2)$}
\label{sec:dualpairing}

We will call two Hopf-$*$ algebras $A$ and $H$ dual to each other if
there is a dual pairing \cite{Takeuchi:1977} between them:
\begin{Definition}
  Let $A$ and $H$ be Hopf-$*$ algebras. A non-degenerate bilinear map
  \begin{equation}
    \lrAngle{\,\cdot\,,\,\cdot\,}: A\times H \longrightarrow \mathbb{C}
    \,, \quad (a,h) \longmapsto \lrAngle{a,h}
  \end{equation}
  is called a dual pairing of $A$ and $H$ if it satisfies
  \begin{equation}
  \begin{aligned}
    (i)&:\quad \lrAngle{\Delta(a),g\otimes h} = \lrAngle{a, gh}
         \,,\quad \lrAngle{a\otimes b,\Delta(h)} = \lrAngle{ab, h}\\
    (ii)&: \quad \lrAngle{a,1} = \varepsilon(a)
         \,,\quad \lrAngle{1,h} = \varepsilon(h)\\
    (iii)&: \quad \lrAngle{a^*, h} = \overline{ \lrAngle{a, (Sh)^*} }\,.
  \end{aligned}
  \end{equation}
\end{Definition}
Remark that for property $(i)$ we have to extend the dual pairing on
tensor products by
\begin{equation}
  \lrAngle{a\otimes b,g\otimes h} := \lrAngle{a, g}\lrAngle{b, h} \,.
\end{equation}
From the properties of the dual pairing it follows that
\begin{equation}    
  \lrAngle{S(a),h} = \lrAngle{a, S(h)} \,.
\end{equation}
The following algebra is dual to $\SUq$ 
\begin{Definition}
  The algebra generated by $E$, $F$, $K$, and $K^{-1}$ with
  commutation relations
  $KK^{-1} = K^{-1}K = 1$ and
  \begin{xalignat}{3}
    KE &= q^2 EK \,,& KF &= q^{-2} FK \,, &
    EF-FE &= \frac{K-K^{-1}}{q-q^{-1}} \,,
  \end{xalignat}
  Hopf structure
  \begin{xalignat}{3}  
  \Delta(E)&= E\otimes K + 1\otimes E \,,& S(E)&=-EK^{-1} \,,&
  \varepsilon(E)&=0 \notag\\
  \Delta(F)&= F\otimes 1 + K^{-1}\otimes F \,, & S(F)&=-KF \,,&
  \varepsilon(F)&=0\\ \notag
  \Delta(K)&= K\otimes K \,,& S(K)&= K^{-1} \,,& \varepsilon(K)&=1
  \end{xalignat}
  and $*$-structure
  \begin{xalignat}{3}
    E^* &= FK \,, & F^* &= K^{-1}E \,, & K^* &= K
  \end{xalignat}
  is called $\suq$, the $q$-deformation of the enveloping algebra
  $\su$ \cite{Kulish:1983,Sklyanin:1985}.
\end{Definition}
The dual pairing of $\suq$ and $\SUq$ is defined on the generators by
\begin{xalignat}{3}
  \lrAngle{E,M^{a}{}_{b}} &:=
    \begin{pmatrix} 0 & 0 \\ q^{\frac{1}{2}} & 0 \end{pmatrix}, &
  \lrAngle{F,M^{a}{}_{b}} &:=
    \begin{pmatrix} 0 & q^{-\frac{1}{2}} \\ 0 & 0 \end{pmatrix},&
  \lrAngle{K,M^{a}{}_{b}} &:=
    \begin{pmatrix} q^{-1} & 0 \\ 0 & q \end{pmatrix}.
\end{xalignat}
There is a universal $\R$-matrix (Sec.~\ref{sec:AppR}) for $\slq$
defined by the formal power series
\begin{equation}
\label{eq:Rmatrix}
  \R = q^{(H\otimes H)/2} \sum_{n=0}^{\infty} R_n(q) (E^n \otimes F^n) \,
\end{equation}
where $R_n(q) := q^{n(n-1)/2} (q-q^{-1})^n ([n]!)^{-1}$, and $K=q^H$
\cite{Drinfeld:1986}.  It is dual to the coquasitriangular map $R$ of
$\SUq$ in the sense that
\begin{equation}
  \lrAngle{\R,g\otimes h} = R(g,h)
\end{equation} 
for all $g,h \in \SUq$. This duality is the reason why we have
introduced the factor $q^{-\frac{1}{2}}$ in the
definition~\eqref{eq:R-Map} of the coquasitriangular map $R$. We will
sometimes write in a Sweedler like notation $\R = \R_{[1]} \otimes
\R_{[2]}$, where the subscripts stand for an index which is summed
over.

\subsection{Computing the Dual of the $q$-Lorentz Group}

The map of the dual pairing of $\slq$ and $\SLq$ naturally extends to
a pairing of the tensor product spaces $\slq \otimes \slq$ and
$SL_q(2,\mathbb{C})\cong \SLq\otimes\SLq$ by
\begin{equation}
  \lrAngle{a\otimes b, g\otimes h}
  :=\lrAngle{a, g} \lrAngle{b, h} 
\end{equation}
for all $a, b \in \slq$ and $g, h\in SL_q(2)$. By
construction, this pairing is non-degenerate. We now want to define a
Hopf algebra structure on $\slq \otimes \slq$ which turns this into 
a pairing of Hopf algebras. Firstly, the multiplication must satisfy 
\begin{equation}
\begin{split}
  \lrAngle{(a\otimes b)(a'&\otimes b'), g\otimes h}
  \stackrel{!}{=} \lrAngle{(a\otimes b)\otimes (a'\otimes b'),
    \Delta(g\otimes h)}\\
  &= \lrAngle{a\otimes a',\Delta(g)}
    \lrAngle{b\otimes b', \Delta(h)}
  = \lrAngle{a a',g} \lrAngle{b b', h}\\
  &= \lrAngle{a a' \otimes b b', g\otimes h}\,.
\end{split}
\end{equation}
Hence, the multiplication on the vector space $\slq\otimes \slq$ must
be defined by $(a\otimes b)(a'\otimes b') = aa' \otimes bb'$, which
means that as an algebra the dual of the $q$-Lorentz group is just the
tensor algebra of two copies of $\slq$. Secondly, we want to define a
coproduct that is consistent with the pairing.
\begin{equation}
\begin{split}
  \lrAngle{\Delta(a&\otimes b),
    (g\otimes h) \otimes (g'\otimes h')}
  \stackrel{!}{=} \lrAngle{a\otimes b,
    (g\otimes h) (g'\otimes h')}\\
  &=\lrAngle{ a\otimes b, gg'_{(2)} \otimes h_{(2)}h' }
    \,R^{-1}(h_{(1)},g'_{(1)}) R(h_{(3)},g'_{(3)})\\
  &=\lrAngle{\Delta(a), g\otimes g'_{(2)} }
    \lrAngle{\Delta(b),h_{(2)}\otimes h'}
    \lrAngle{\R^{-1},h_{(1)}\otimes g'_{(1)}}
    \lrAngle{\R,h_{(3)}\otimes g'_{(3)}}\\
  &=\lrAngle{a_{(1)} , g} \lrAngle{\R^{-1}_{[1]},h_{(1)}}
    \lrAngle{b_{(1)} , h_{(2)}} \lrAngle{\R_{[1]},h_{(3)}}\\
    &\qquad\times
    \lrAngle{\R^{-1}_{[2]},g'_{(1)}} \lrAngle{a_{(2)} , g'_{(2)}}
    \lrAngle{\R_{[2]},g'_{(3)}} \lrAngle{b_{(2)} , h'}\\
  &=\lrAngle{a_{(1)} , g} \lrAngle{\R^{-1}_{[1]} b_{(1)}\R_{[1]},h}
    \lrAngle{\R^{-1}_{[2]} a_{(2)}\R_{[2]},g'} \lrAngle{b_{(2)} , h'}\\
  &=\lrAngle{(a_{(1)}\otimes
    \R^{-1}_{[1]} b_{(1)}\R_{[1]})\otimes
    (\R^{-1}_{[2]} a_{(2)}\R_{[2]}\otimes b_{(2)}),
     (g\otimes h) \otimes (g'\otimes h')}  
\end{split}
\end{equation}
From the last line we read off the coproduct
\begin{equation}
  \Delta(a\otimes b) =\R^{-1}_{23}
  \Delta^{\!\otimes 2} (a\otimes b)\R_{23} \,,
\end{equation}
where $\R_{23}= 1\otimes \R \otimes 1 $ and $\Delta^{\!\otimes 2}
(a\otimes b) = (a_{(1)}\otimes b_{(1)})\otimes (a_{(2)}\otimes
b_{(2)})$. This tells us, that the coproduct of the $q$-Lorentz
algebra is the tensor coproduct of $\slq\otimes \slq$ with the two
inner tensor factors twisted by the universal $\R$-matrix.

Thirdly, the same reasoning for the antipode
\begin{equation}
\begin{split}
  \lrAngle{S(a&\otimes b), g\otimes h}
  \stackrel{!}{=} \lrAngle{a\otimes b, S(g\otimes h)} =
  \lrAngle{a\otimes b,(1\otimes Sh)(Sg\otimes 1)}\\
  &= \lrAngle{a\otimes b,(Sg)_{(2)}\otimes (Sh)_{(2)}}
    \,R^{-1}\bigl((Sh)_{(1)},(Sg)_{(1)}\bigr)
    R\bigl((Sh)_{(3)},(Sg)_{(3)}\bigr)\\
  &=\lrAngle{a , S(g_{(2)})}
    \lrAngle{b , S(h_{(2)})}
    \,R^{-1}\bigl(h_{(3)},g_{(3)}\bigr)
    R\bigl(h_{(1)},g_{(1)}\bigr)\\
  &=\lrAngle{\R_{[2]},g_{(1)}} \lrAngle{S(a) , g_{(2)}}
    \lrAngle{\R^{-1}_{[2]},g_{(3)}} \lrAngle{\R_{[1]},h_{(1)}}
    \lrAngle{S(b),h_{(2)}} \lrAngle{\R^{-1}_{[1]},h_{(3)}}\\
  &=\lrAngle{\R_{[2]}S(a)\R^{-1}_{[2]},g}
    \lrAngle{\R_{[1]}S(b)\R^{-1}_{[1]},h}\\
  &=\lrAngle{\R_{[2]}S(a)\R^{-1}_{[2]}\otimes
    \R_{[1]}S(b)\R^{-1}_{[1]}, g\otimes h}
\end{split}
\end{equation}
leads to
\begin{equation}
  S(a\otimes b) = \R_{21} (Sa\otimes Sb) \R^{-1}_{21} \,,
\end{equation}
where $\R_{21}=\R_{[2]} \otimes \R_{[1]}$.  The antipode is the tensor
antipode twisted by the transposed universal $R$-matrix.

The counit $\varepsilon(a\otimes b)= \varepsilon(a)\varepsilon(b)$
follows directly from the definition of the pairing. Finally, we need
to calculate the star structure.
\begin{equation}
\begin{split}
  \lrAngle{(a&\otimes b)^*, g\otimes h}
  \stackrel{!}{=}
  \overline{ \lrAngle{a\otimes b, S((g\otimes h)^*)} } =
  \overline{ \lrAngle{S(a\otimes b), h^* \otimes g^*} }\\
  &=\overline{
  \lrAngle{\R_{[2]}S(a)\R^{-1}_{[2]},h^*}
  \lrAngle{\R_{[1]}S(b)\R^{-1}_{[1]},g^*} }\\
  &=\lrAngle{\R_{[1]}a^*\R_{[1]}^{-1},h}
    \lrAngle{\R_{[2]}b^*\R_{[2]}^{-1},g}
\end{split}
\end{equation}
Here we have used that $\R$ is real, $\R^{*\otimes *} = \R_{21}$.
Thus, we find
\begin{equation}
  (a\otimes b)^* = \R_{21} (b^*\otimes a^*) \R_{21}^{-1} \,,  
\end{equation}
which completes the structure of the $q$-Lorentz algebra. 

To summarize, we have the following
\begin{Proposition}
\label{th:LorentzAlgebra}
  The tensor product algebra $\slq\otimes \slq$ with the Hopf-$*$-structure
  \begin{equation}
  \label{eq:LorentzHopf}
  \begin{gathered}
    \Delta(a\otimes b) =\R^{-1}_{23}
    \Delta^{\!\otimes 2} (a\otimes b)\R_{23} \,,\quad
    S(a\otimes b) = \R_{21} (Sa\otimes Sb) \R^{-1}_{21}\\
    \varepsilon(a\otimes b)= \varepsilon(a)\varepsilon(b) \,,\quad
    (a\otimes b)^* =\R_{21}(b^*\otimes a^*)\R_{21}^{-1}
  \end{gathered}
  \end{equation}
  is the Hopf-$*$-dual of the $q$-Lorentz group $SL_q(2,\mathbb{C})$.
  Therefore, we will call it the $q$-Lorentz algebra $\slC$
  \cite{Majid:1993}.
\end{Proposition}
There are two universal $\R$-matrices of the $q$-Lorentz algebra,
which are composed of the $\R$-matrix of $\slq$ according to
\begin{xalignat}{2}
\label{eq:LorentzR}
  \RI &= \R^{-1}_{41}\R^{-1}_{31}\R_{24}\R_{23}\,, &
  \RII &= \R^{-1}_{41}\R_{13}\R_{24}\R_{23} \,.
\end{xalignat}
$\RI$ is anti-real while $\RII$ is real.

\chapter{Structure of the $q$-Lorentz Algebra}
\label{sec:LorentzStructure}

\section{Representation Theory of the $q$-Lorentz Algebra}

\subsection{The Clebsch-Gordan Series of $\slq$}

Let us review some facts about the representation theory of $\slq$ and
$\suq$ \cite{Sklyanin:1983}. For any $j\in\frac{1}{2}\mathbb{N}_0$
there is an irreducible representation on the $(2j+1)$-dimensional
Hilbert space $D^j$ with orthonormal basis $\{\Ket{j,m}\,|\,
m=-j,-j+1,\ldots j\}$ and representation map $\rho^j:\slq \rightarrow
\mathrm{Aut}(D^j)$ given by\footnote{There is a second series of
  irreducible representations with negative eigenvalues of $K$, which
  we will not take into account, since they have no undeformed limit.}
\begin{equation}
\label{eq:suqIrreps}
\begin{aligned}
  \rho^j(E)\Ket{j,m} &= q^{(m+1)} \sqrt{[j+m+1][j-m]}\,\Ket{j,m+1}\\
  \rho^j(F)\Ket{j,m} &= q^{-m} \sqrt{[j+m][j-m+1]}\,\Ket{j,m-1}\\
  \rho^j(K)\Ket{j,m} &= q^{2m}\Ket{j,m} \,.
\end{aligned}
\end{equation}
For the real form $\suq$ these are even $*$-representations. $D^0$ is
called the scalar representation, $D^\frac{1}{2}$ the fundamental or
spinor representation, and $D^1$ the vector representation.

Recall that the coproduct of a Hopf algebra enables us to construct
the tensor product of two representations: Let $D^j$ and $D^{j'}$ be
representations of $\slq$ as defined above, with representation maps
$\rho^j$ and $\rho^{j'}$. Then there is a representation on the tensor
product space $D^j\otimes D^{j'}$ with representation map
$(\rho^j\otimes \rho^{j'})\circ \Delta$. We denote this tensor product
of representations also by $D^j\otimes D^{j'}$.

In general, the tensor product of two irreducible representation is no
longer irreducible. In fact, in complete analogy to the classical case
we have an isomorphism of representations
\begin{equation}
  \label{eq:CGSeries}
  D^j\otimes D^{j'} \cong D^{|j-j'|} \oplus D^{|j-j'|+1}
    \oplus \ldots \oplus D^{j+j'} 
\end{equation}
decomposing the tensor product into the Clebsch-Gordan series.
This isomorphism, viewed as a transformation of basis
\begin{equation}
  \label{eq:CGC-def}
  \Ket{j,m} = \sum_{j_1,j_2,m_1,m_2} \CGC{j_1}{j_2}{j}{m_1}{m_2}{m}\,
    \Ket{j_1,m_1}\otimes \Ket{j_2,m_2}
\end{equation}
defines the $q$-Clebsch-Gordan coefficients, which can be calculated
in a closed form (Sec.~\ref{sec:AppClebsch1}). The two most important
cases are the construction of a scalar and the construction of a
vector out of two vector representations, where the right hand side of
Eq.~\eqref{eq:CGC-def} may be viewed as the scalar and the vector
product of two 3-vectors.

\subsection{Clebsch-Gordan Coefficients of the $q$-Lorentz Algebra}

As an algebra the $q$-Lorentz algebra is the tensor product of two
$\slq$. Hence, every finite irreducible representation is composed of
two irreducible representations $D^{j_1}$ and $D^{j_2}$ of $\slq$,
that is, the vector space $D^{(j_1,j_2)} := D^{j_1}\otimes D^{j_2}$
with the representation map $\rho^{(j_1,j_2)} := \rho^{j_1}\otimes
\rho^{j_2}$. Viewing the decomposition of the $q$-Lorentz algebra into
two $\slq$ as chiral decomposition, we call $D^{j_1}$ the left handed
and $D^{j_2}$ the right handed part of the representation.
$D^{(j_1,j_2)}$ is not a $*$-representation, since the $*$-operation
of the $q$-Lorentz algebra is not the tensor product of the $*$'s of
$\slq$. Therefore, all finite irreducible representations are
non-unitary. This is a sign of the non-compactness of the $q$-Lorentz
algebra on a representation theoretic level.

Next, we consider the tensor product of two representations. Again,
its vector space is just the tensor product $D^{(j_1,j_2)}\otimes
D^{(j'_1,j'_2)}$. The representation map is again $\rho =
(\rho^{(j_1,j_2)}\otimes \rho^{(j'_1,j'_2)})\circ \Delta$, where
$\Delta$ is now the coproduct of the $q$-Lorentz algebra as defined in
Eq.~\eqref{eq:LorentzHopf}. The coproduct is calculated by, firstly,
taking the $\slq$ coproduct of the two $\slq$ tensor factors, then
interchanging the 2. and 3. tensor factor, and, finally, conjugating
with the universal $\R$-matrix in the 2. and 3. position of the 4-fold
tensor product. Algebraically, the last step is a complicated inner
automorphism, since $\R$ exists only as an infinite formal power
series. However, when we apply the representation maps, $\R$ becomes a
finite $(j_2 j'_1)\times (j_2 j'_1)$ matrix $R = (R^{ab}{}_{cd})$
\begin{equation}
  \rho^{(j_1,j_2)}\otimes \rho^{(j'_1,j'_2)}(\R_{23})
  = 1\otimes \bigl( (\rho^{j_2}\otimes \rho^{j'_1})(\R)\bigr) \otimes 1
  =: 1 \otimes R \otimes 1 \,,
\end{equation}
and the inner automorphism becomes a simple basis transformation. 

Putting things together, we see how to reduce the product of two
$q$-Lorentz representations. Up to a change of basis we reduce the
tensor product of the 1. with the 3. and the 2. with the 4. of the
$\slq$-subrepresentations, each by means of the Clebsch-Gordan series
of $\slq$.
\begin{equation}
\label{eq:LorentzCGSeries}
  D^{(j_1,j_2)} \otimes D^{(j'_1,j'_2)} \cong
  \bigoplus_{\substack{|j_1 - j'_1| \leq k_1 \leq j_1 + j'_1 \\
                       |j_2 - j'_2| \leq k_2 \leq j_2 + j'_2 }}
  D^{(k_1,k_2)}
\end{equation}
Written out for the important case of the product of two vector
representations, this formula becomes
\begin{equation}
\label{eq:CGSeries1}
  D^{(\frac{1}{2},\frac{1}{2})} \otimes D^{(\frac{1}{2},\frac{1}{2})}
  \cong
  D^{(0,0)} \oplus D^{(1,0)} \oplus D^{(0,1)} \oplus D^{(1,1)}\,,
\end{equation}
which corresponds to the decomposition of a 4$\times$4 matrix viewed
as a second rank Lorentz tensor into the scalar trace part, a left and
a right chiral 3-vector, and the traceless symmetric part
(Sec.~\ref{sec:AppClebsch3}).

So far, the representation theory is in complete accordance with the
undeformed case. New is the appearance of an $R$-matrix, which matters
as soon as we want to write down the above isomorphisms explicitly.
The matrix representing isomorphism~\ref{eq:LorentzCGSeries} is the
product of two Clebsch-Gordan coefficients contracted with the
$R$-matrix. Musing for a while about the right positions of the
indices, we find
\begin{multline}
\label{eq:CGSeries2}
  \Ket{(k_1,k_2),(n_1,n_2)} =
  \sum
  \CGC{j_1}{j'_1}{k_1}{m_1}{b}{n_1} 
  \CGC{j_2}{j'_2}{k_2}{a}{m'_2}{n_2}\\
  \times (R^{-1})^{m_2 m'_1}{}_{ab}\,  
  \Ket{(j_1,j_2),(m_1,m_2)} \otimes
  \Ket{(j'_1,j'_2),(m'_1,m'_2)} \,,
\end{multline}
where we sum over repeated indices, and where the labeling of the free
indices is the same as in Eq.~\eqref{eq:LorentzCGSeries}. This defines
the Clebsch-Gordan coefficients of the $q$-Lorentz algebra
\begin{multline}
  \left[
  \begin{array}{ccc|ccc}
    j_1 & j'_1 & k_1 & m_1 & m'_1 & n_1\\
    j_2 & j'_2 & k_2 & m_2 & m'_2 & n_2
  \end{array}
  \right]_q := \\
  \sum_{a,b}
  \CGC{j_1}{j'_1}{k_1}{m_1}{b}{n_1} 
  \CGC{j_2}{j'_2}{k_2}{a}{m'_2}{n_2}
  (R^{-1})^{m_2 m'_1}{}_{ab} \,.
\end{multline}

\section{Tensor Operators}

\subsection{Tensor Operators in Hopf Algebras}
\label{sec:TensorHopf}

Recall that there is a left and right action of any Hopf algebra $H$
on itself given by
\begin{xalignat}{2}
  \label{eq:HopfAdjointAction}
  \adL(g)\tr h &:= g_{(1)} h S(g_{(2)}) \,,&
  h \tl \adR(g) &:= S(g_{(1)}) h g_{(2)}
\end{xalignat}
for $g,h\in H$, called the left and right Hopf adjoint action,
respectively. In general, this action will be highly reducible. In
fact, if a linearly independent set $\{A_\mu \in H\}$ of operators
generates an invariant subspace $D$ of the left Hopf adjoint action,
this induces a matrix representation map $\rho$ of $H$ by
\begin{equation}
\label{eq:Tensor1}
  \adL (h) \tr A_\mu = A_{\mu'} \rho(h)^{\mu'}{}_\mu \,,
\end{equation}
turning $D$ into a representation. The set of operators $\{A_\mu\}$
with this property is called a left $D$-tensor operator of $H$,
indicated by a lower index. It will be called irreducible if $D$ is
irreducible. If in addition $H$ is equipped with a $*$-operation, we
can demand that $D$ is a $*$-representation.

There are other useful types of tensor operators. If a set of
operators $\{A^\mu\}$ transforms as
\begin{equation}
\label{eq:Tensor2}
  (\adL h)\tr A^\mu = \rho(Sh)^{\mu}{}_{\mu'} A^{\mu'} \,,
\end{equation}
we will call it a left upper or congredient tensor operator, denoted
by an upper index. Its transformation is congredient in the sense that
\begin{equation}
\begin{split}
  (\adL h)\tr(A_\mu B^\mu)
  &= [(\adL h_{(1)})\tr A_\mu] [(\adL h_{(2)})\tr B^\mu] \\
  &= A_{\mu'} \rho(h_{(1)})^{\mu'}{}_\mu
    \,\rho\bigl(S(h_{(1)})\bigr)^{\mu}{}_{\mu''} B^{\mu''} 
  = \varepsilon(h)\,A_\mu B^\mu  \,,
\end{split}
\end{equation}
that is, $A_\mu B^\mu$ is a scalar operator. If $g^{\mu\nu}$ is a
metric for the representation under consideration and $A_\mu$ and
$B_\nu$ are left tensor operator then $g^{\mu\nu} A_\mu B_\nu$ is a
scalar. This is true for the $q$-spinor metric $\varepsilon^{ab}$, the
metric $g^{AB}$ of vector representations of $\suq$ and the
$q$-Minkowski metric $\eta^{\mu\nu}$, as defined in
Eqs.~\eqref{eq:SpinorMetric}, \eqref{eq:epsgdef}, and
\eqref{eq:MinkowskiMetric}, respectively. We conclude, that the
convention for the position of tensor operator indices is consistent
with raising and lowering indices as usual, $A^\mu =
g^{\mu\mu'}A_{\mu'}$. Moreover, we conclude that
\begin{equation}
\label{eq:Tensor6}
  g^{\mu\mu'}g_{\nu'\nu} \rho(h)^{\nu'}{}_{\mu'}
  = \rho(Sh)^{\mu}{}_{\nu} \,.
\end{equation}
If we deal with a Hopf-$*$-algebra and $\rho$ is a $*$-representation,
we can apply $*$ to Eq.~\eqref{eq:Tensor1} and get
\begin{multline}
  (*\circ S)(h_{(2)}) (A_\mu)^* (h_{(1)})^*
  = \bigl[( (Sh)^* )_{(1)}\bigr]
    (A_\mu)^* S\bigl[( (Sh)^* )_{(2)}\bigr]
  = (\adL (Sh)^*)\tr (A_\mu)^* \\
  = (A_{\mu'})^* \overline{ \rho(h)^{\mu'}{}_\mu }
  = (A_{\mu'})^* \rho(h^*)^{\mu}{}_{\mu'} 
  = (A_{\mu'})^* \rho(S[(Sh)^*])^{\mu}{}_{\mu'} \,,
\end{multline}
from which we deduce 
\begin{equation}
  (\adL (Sh)^*)\tr (A_\mu)^* =
  \rho(S[(Sh)^*])^{\mu}{}_{\mu'}\, (A_{\mu'})^*  \,.
\end{equation}
Comparing this with Eq.~\eqref{eq:Tensor2}, we conclude that
$(A_\mu)^*$ is a congredient left tensor operator.

Let us now consider tensor operators $A^{\tilde{\mu}}$ with respect to
the right Hopf-adjoint action
\begin{equation}
\label{eq:Tensor3}
  A^{\tilde{\mu}} \tl (\adR h)
  = S(h_{(1)}) A^{\tilde{\mu}}  h_{(2)}
  = \rho(h)^{\mu}{}_{\mu'} \, A^{\tilde{\mu}'} \,,
\end{equation}
which we call \emph{right} upper tensor operators, distinguished from
left upper tensor operators by putting a tilde over their indices.

Let $A^\mu$ be a left upper tensor operator and let there be an
extension of the antipode of $H$ on $A^\mu$, for example, $A^\mu$
might be an element of $H$. Then we can apply $S$ to
Eq.~\eqref{eq:Tensor2} and obtain
\begin{multline}  
   S(S(h_{(2)})) S(A^\mu )  S(h_{(1)})
   = S((Sh)_{(1)})) S(A^\mu)  (Sh)_{(2)}) \\
   = S(A^\mu) \tl (\adR Sh) 
   = \rho(Sh)^{\mu}{}_{\mu'} \, S(A^{\mu'}) \,.
\end{multline}
Thus, $S(A^\mu)$ is a right upper tensor operator. 

Note, that within a Lie algebra we would have $S(A^\mu) = -A^\mu$.
Hence, in a Lie algebra a right tensor operator is the same as a left
tensor operator. This is why in the undeformed case we need not
distinguish between indices with and without a tilde.

Finally, we define a right lower tensor operator $A_{\tilde{\mu}}$ to
transform as $S(A_\mu)$, that is,
\begin{equation}
\label{eq:Tensor4}
  A_{\tilde{\mu}} \tl (\adR h)
  = A_{\tilde{\mu}'}\rho(S^{-1}h)^{\mu'}{}_{\mu} \,.
\end{equation}
One can check that we can raise and lower indices as usual,
$A_{\tilde{\mu}} = g_{\mu\nu}A^{\tilde{\nu}}$, and that
$A^{\tilde{\mu}}B_{\tilde{\mu}}$ is a scalar operator. Note that being
a left or a right scalar is the same thing: A scalar is an operator
that commutes with $H$.

\subsection{Tensor Operators of $\suq$}

Most tensor operators of $\suq$ that we will consider are $D^0$-tensor
operators, which will be called $\slq$-scalars, and $D^1$-tensor
operators, called $\slq$-vectors. One big advantage of grouping
several operators to a $\slq$-tensor operator lies in the
$q$-Wigner-Eckart theorem:
\begin{Theorem}
\label{th:WignerEckart}
  Let $A_\mu$ be a left $D^\lambda$-tensor operator of $\suq$ and let
  there be a representation of $\suq$ with irreducible
  subrepresentations $D^j$ and $D^{j'}$ with bases $\{\Ket{j,m}\}$ and
  $\{\Ket{j',m'}\}$. Then there exists a number $\rBraket{j'}{A}{j}$
  such that
  \begin{equation}
    \Braket{j',m'}{A_\mu}{j,m} = \CGC{\lambda}{j}{j'}{\mu}{m}{m'}
    \rBraket{j'}{A}{j}
  \end{equation}
  for all $m$, $m'$. This number is called the reduced matrix element
  of the tensor operator $A_\mu$ \cite{Klimyk:1992}.
\end{Theorem}
If we have degeneracy of the $\Ket{j,m}$ basis, the reduced matrix
elements will depend on additional quantum numbers but not on $m$.
Whenever a $q$-Clebsch-Gordan coefficient
$\CGC{\lambda}{j}{j'}{\mu}{m}{m'}$ vanishes for all $m$, $m'$, the
reduced matrix element is not defined uniquely. In that case we set
$\rBraket{j'}{A}{j} := 0$ for convenience.

Looking at the definition~\eqref{eq:HopfAdjointAction}, we see that
$\adL(g)\tr (hh') = (\adL(g_{(1)})\tr h) (\adL(g_{(2)})\tr h')$.
Hence, the product of a $D$- and a $D'$-tensor operator is a $D\otimes
D'$-tensor operator. Just as for the representations of $\slq$ we have
a Clebsch-Gordan decomposition of the product of tensor operators:
\begin{Proposition}
\label{th:TensorDecomp}
  Let $A_\alpha$ be a $D^a$-tensor operator and $B_\beta$ a
  $D^b$-tensor operator of $\slq$. Then
  \begin{equation}
  \label{eq:Tensor5}
    C_\gamma := \sum_{\alpha, \beta}
      \CGC{a}{b}{c}{\alpha}{\beta}{\gamma} A_\alpha B_\beta
  \end{equation}
  is a $D^c$-tensor operator of $\slq$.
\end{Proposition}

If we now take the matrix elements of a tensor operator $C_\gamma$
constructed in this way, we find with the aid of the $q$-Wigner-Eckart
theorem relations between the reduced matrix elements
\begin{equation}
\label{eq:RacahRed}
  \rBraket{j''}{C}{j} = \sum_{j'} \RC{a}{b}{j}{c}{j'}{j''} 
  \rBraket{j''}{A}{j'} \rBraket{j'}{B}{j} \,.      
\end{equation}
Here $\mathrm{R}_q$ denote the $q$-Racah coefficients defined by the
expression
\begin{multline}
  \RC{a}{b}{j}{c}{j'}{j''} :=
    \CGC{c}{j}{j''}{\gamma}{m}{m''}^{-1}\\ \times
    \sum_{\alpha, \beta, m'} \CGC{a}{b}{c}{\alpha}{\beta}{\gamma}
    \CGC{a}{j}{j''}{\alpha}{m'}{m''}\CGC{b}{j}{j'}{\beta}{m}{m'}\,,
\end{multline}
which can be proven not to depend on $m$, $m''$ as the arguments of
$\mathrm{R}_q$ indicate. Values of the $q$-Racah coefficients are
given in Sec.~\ref{sec:AppClebsch1}.

The two cases of Eq.~\eqref{eq:Tensor5} that we encounter most
frequently are the construction of a scalar and the construction of a
vector operator out of two vector operators. This suggests the
definition
\begin{xalignat}{2}
\label{eq:epsgdef}
  g^{AB} &:= -\sqrt{[3]} \CGC{1}{1}{0}{A}{B}{0}\,, &
  \varepsilon^{AB}{}_{C} &=
    -\sqrt{\frac{[4]}{[2]}} \CGC{1}{1}{1}{A}{B}{C} \,,
\end{xalignat}
where the capital Roman indices run through $\{-1,0,1\} = \{-,3,+\}$.
Values are given in Sec.~\ref{sec:AppClebsch2}.
Proposition~\ref{th:TensorDecomp} tells us that we can define a scalar
and a vector product of two vector operators $X_A$ and $Y_A$ by
\begin{xalignat}{2}
\label{eq:CrossDot}
  \vec{X}\cdot \vec{Y} &:= X_A Y_B\, g^{AB} \,,&
  (\vec{X}\times \vec{Y})_C &:= \I X_A Y_B \,\varepsilon^{AB}{}_{C} \,,
\end{xalignat}
where the imaginary unit is needed for the right undeformed
limit.\footnote{See Sec.~\ref{sec:AppClebsch2}, in particular the
  remark above Eq.~\eqref{eq:epsidentities3}.}  By definition, the
scalar product is a scalar and the vector product is a vector operator
in the sense of Eq.~\eqref{eq:Tensor1}.

\subsection{The Vector Form of $\suq$}

For a set of operators to be interpreted as $q$-angular momentum, it
will have to generate the symmetry of rotations on the one hand, but
on the other hand it will itself have to transform like a vector under
rotations. In other words, this set must be a vector operator
generating $\suq$. In the $EFK$-form of $\slq$ it is not obvious, what
this vector operator could be.

We begin our search for such a vector operator by giving the explicit
conditions for $A_\mu$ to be a irreducible $D^j$-tensor operator of
$\slq$: Inserting Eqs.~\eqref{eq:suqIrreps} in Eq.~\eqref{eq:Tensor1}
we get
\begin{equation}
\label{eq:sutensorop}
\begin{aligned}
  E A_\mu - A_\mu E &= q^{(\mu+1)}
    \sqrt{[j+\mu+1][j-\mu]}\,\,A_{\mu+1} K\\
  F A_\mu - q^{-2\mu} A_\mu F &= q^{-\mu}
    \sqrt{[j+\mu][j-\mu+1]}\,\,A_{\mu-1}\\
  K A_\mu &= q^{2\mu} A_\mu K \,.
\end{aligned}
\end{equation}
To find a vector operator in $\suq$ satisfying these conditions we
first look for a highest weight vector $J_+$ and let $\slq$ act on it by
the left Hopf-adjoint action, giving us the subrepresentation generated by
$J_+$. The condition $\adL(E)\tr J_+ = 0$ for $J_+$ to be a highest
weight vector is equivalent to $[E,J_+]=0$.  Thus, $J_+$ must be in the
centralizer of $E$, a very restrictive condition most obviously
satisfied by $E$ itself. The results of the Hopf-adjoint action of
the ladder operators $E$ and $F$ on $E$ are
\begin{equation}
\begin{gathered}
  \adL(F)\tr E = K^{-1}(KFE-EKF) \,, \quad
  \adL(F^2)\tr E = - [2]KF\\ \,,\quad
  \adL(F^3)\tr E = 0 \,,\quad \adL(EF)\tr E = [2]E \\
  \adL(EF^2)\tr E = [2] K^{-1}(KFE-EKF) \,,
\end{gathered}
\end{equation}
which shows that we can indeed interpret $E$ as a highest weight
vector of a vector representation. Comparing the Hopf-adjoint action
with the vector representation as given in Eqs.~\eqref{eq:suqIrreps},
one finds that
\begin{equation}
\label{eq:DefL}
\begin{aligned}
  J_{-} &:= q[2]^{-\frac{1}{2}}KF \\
  J_3    &:= -q[2]^{-1} K^{-1}(KFE-EKF) = [2]^{-1} (q^{-1}EF-qFE) \\
  J_{+} &:= -[2]^{-\frac{1}{2}}E 
\end{aligned}
\end{equation}
form a vector operator.\footnote{Elsewhere \cite{Lorek:1997a}, the
  vector generators have been defined as $L^A = -q^{-3}J_A$.}

How can we describe the subalgebra of $\suq$ generated by $J_A$?
After some calculations we find that the commutation relations of the
$J$'s do not close. Since the commutation relations~\eqref{eq:Tensor1}
are given by the adjoint action of the set of generators on itself,
this is due to the fact that coproduct and antipode of the $J$'s
cannot be expressed by $J$'s again. We can help ourselves out by
introducing the additional generator
\begin{equation}
\label{eq:DefW}
  W := K - \lambda J_3 = K - \lambda [2]^{-1} (q^{-1}EF-qFE)\,,
\end{equation}
so the commutation relations can be written as
\begin{xalignat}{3}
\label{eq:LL-Rel}
  J_A J_B \varepsilon^{AB}{}_C &= W J_C \,, &
  J_A W &= W J_A \,, &
  W^2 - \lambda^2 J_A J_B g^{AB} &= 1 \,,
\end{xalignat}
where the last equation expresses that $W$ and the $J$'s are not
algebraically independent. The $*$-structure reads on the generators
\begin{xalignat}{3}
\label{eq:LL-Star}
  J_{+}^* &= -q J_{-}\,, & J_{3}^* &= J_{3}\,, & W^* &= W \,,
\end{xalignat}
that is, $(J_A)^* = J^A$. We will call the subalgebra of $\slq$
generated by $J_A$, $W$ with relations~\eqref{eq:LL-Rel} and
$*$-structure~\eqref{eq:LL-Star} the vectorial form of $\suq$. Note
that the vectorial form of $\slq$ is a proper subalgebra of $\slq$
since it does not contain $K^{-1}$. We do need $K^{-1}$ to write down
the Hopf structure: the coproduct
\begin{equation}
\begin{aligned}
  \Delta (J_{\pm}) &= J_{\pm} \otimes K + 1 \otimes J_{\pm} \\
  \Delta (J_3) &= J_3 \otimes K + K^{-1}\otimes J_3
  +\lambda (qK^{-1}J_+ \otimes J_- +q^{-1} K^{-1}J_- \otimes J_+) \\
  \Delta (W) &= W \otimes K - \lambda K^{-1}\otimes J_3 
  -\lambda^2 (q K^{-1}J_+ \otimes J_- + q^{-1} K^{-1}J_- \otimes J_+) \,,
\end{aligned}
\end{equation}
the antipode
\begin{equation}
\begin{aligned}
  S(J_{\pm}) &= - J_{\pm} K^{-1} \\
  S(J_3) &= J_3 - \lambda^{-1}(K - K^{-1})\\
  S(W) &= W \,, 
\end{aligned}
\end{equation}
and the counit $\varepsilon(J_A) = 0$, $\varepsilon(W) = 1$.

\section{The $q$-Lorentz Algebra as Quantum Double}
\label{sec:QuantumDouble}

\subsection{Rotations and the $\SUq^\op$ Algebra of Boosts}

In Sec.~\ref{sec:LorentzCommutation} the commutation relations of the
$q$-Lorentz group have been chosen to preserve an $SU_q(2)$
substructure, physically interpreted as rotations. That is, the
multiplication of the two copies of $SL_q(2)$ is a Hopf-$*$-homomorphism
projecting the $q$-Lorentz group onto $SU_q(2)$. On the quantum
algebra level, the dual of multiplication is comultiplication. Hence, the
mapping
\begin{equation}
  i:\suq \overset{\Delta}{\longrightarrow} \slq \otimes \slq = \slC
\end{equation}
ought to define a $\suq$ Hopf-$*$-subalgebra of the $q$-Lorentz
algebra.

Given the properties of the coproduct, it is obvious that $i$ is an
algebra homomorphism. It is less clear, whether $i$ preserves the Hopf
structure and the $*$-structure of $\suq$. For the coproducts we find
\begin{align}
  (\Delta_\slC \circ i)(h) 
    &= \R^{-1}_{23}
    (h_{(1)}\otimes h_{(3)}\otimes h_{(2)}\otimes h_{(4)}) \R_{23} 
    = h_{(1)}\otimes h_{(2)}\otimes h_{(3)}\otimes h_{(4)} \notag\\
    &= \bigl((i\otimes i) \circ \Delta_\slq \bigr)(h) \,,
\end{align}
which shows that $i$ is a coalgebra map. In the same manner we find that
$i$ preserves the counit (trivial), the antipode
\begin{align}
  (S_\slC \circ i)(h) 
  &=\R_{21}\bigl( S(h_{(1)})\otimes S(h_{(2)}) \bigr)\R^{-1}_{21}
  =\R_{21}\bigl( (Sh)_{(2)}\otimes (Sh)_{(1)} \bigr)\R^{-1}_{21} \notag\\
  &= (Sh)_{(1)}\otimes (Sh)_{(2)} = (i \circ S_\slq)(h) \,,
\end{align}
and the $*$-structure
\begin{equation}
\begin{split}
  (i(h))^* 
  &=\R_{21}\bigl( (h_{(2)})^* \otimes (h_{(1)})^* \bigr) \R^{-1}_{21}
  =\R_{21}\bigl( (h^*)_{(2)}\otimes (h^*)_{(1)} \bigr) \R^{-1}_{21}\\
  &= (h^*)_{(1)}\otimes (h^*)_{(2)} = i(h^*) \,.
\end{split}
\end{equation}
We conclude that $i(\suq)$ is indeed a $\suq$ Hopf-$*$ subalgebra of
the $q$-Lorentz algebra.\footnote{It is the appearance of the
  $\R$-matrices in the Hopf structure of $\slC$, which ensures the
  compliance of the embedding $i$ with the Hopf structures. This is
  why $\mathcal{U}_q(\mathrm{so}_4)$ does not possess a Hopf
  subalgebra of rotations.} Since in the undeformed case the embedding
of the rotations in the Lorentz algebra is given by the coproduct,
too, $i(\suq)$ has the right undeformed limit. This strongly suggests
to interpret $i(\suq)$ as the quantum subsymmetry of physical
rotations.

There is another Hopf-$*$ subalgebra of $\slC$. Let $\lrAngle{ ~,~}$
denote the dual pairing of $\slq$ and $\SLq$ as defined in
Sec.~\ref{sec:dualpairing}. We define a map $j:\SUq \rightarrow \slC$
by
\begin{equation}
  j(h) := \lrAngle{\R^{-1}_{31}\R_{23}, h_3} \,,
\end{equation}
where the subscripts denote the position in the tensor product, $h_3
:= 1\otimes 1 \otimes h$, and where the dual pairing acts only on the
third tensor factor.  Let us show some properties of this map.  we
have
\begin{align}
  j(gh)
  &= \lrAngle{\R^{-1}_{31}\R_{23}, g_3 h_3}
  = \lrAngle{\Delta_3 (\R^{-1}_{31}\R_{23}), g_3 h_4} 
  = \lrAngle{\R^{-1}_{41}\R^{-1}_{31} \R_{24}\R_{23},g_3 h_4} \notag\\
  &= \lrAngle{\R^{-1}_{41} \R_{24}\R^{-1}_{31} \R_{23}, g_3 h_4}
  = j(h)j(g) \,,
\end{align}
telling us that $j$ is an algebra anti-homomorphism. Next we consider
the coproduct
\begin{align}
  (\Delta_\slC\circ j)(h)
  &= \lrAngle{\Delta\otimes \id\ (\R^{-1}_{31}\R_{23}), h_3}
  = \lrAngle{\R^{-1}_{23} R^{-1}_{51}
    \R^{-1}_{53} \R_{25}\R_{45} \R_{23}, h_5} \notag\\
  &= \lrAngle{ \R^{-1}_{51} \R_{25} \R^{-1}_{53}\R_{45}, h_5}
  = \bigl( (j\otimes j)\circ\Delta_\SUq\bigr)(h) \,,
\end{align}
so $j$ is a coalgebra homomorphism, too. The calculation for the
counit is trivial.

So far we can say that $j$ is a bialgebra homomorphism from $\SUq^\op$
to the $q$-Lorentz algebra. $\SUq^\op$ becomes a Hopf algebra,
when we equip it with a antipode and $*$-structure according to
\begin{xalignat}{2}
  S^\op &:= S^{-1}\,, & *^\op := * \circ S^2 \,,
\end{xalignat}
where $S$ is the usual antipode of $\SUq$. Let us check now if $j$
preserves this Hopf structure as well.  We begin with the antipode
\begin{align}
  (j\circ S^\op)(h)
  &= \lrAngle{\R^{-1}_{31}\R_{23}, S^{-1}(h_3)}
  = \lrAngle{\R^{-1}_{23}\R_{31}, h_3}
  = \lrAngle{\R_{21}\R_{31}\R^{-1}_{23}\R^{-1}_{21}, h_3 }\notag\\
  &= \lrAngle{\R_{21} [ (S\otimes S \otimes \id)
    (\R^{-1}_{31}\R_{23}) ] \R^{-1}_{21}, h_3 }
  = (S_\slC \circ j)(h) \,,
\end{align}
which is indeed preserved. Finally, we have the $*$-structure
\begin{align}
  j\bigl(h^{*^\op}\bigr)
  &= \lrAngle{\R^{-1}_{31}\R_{23}, (S^2 h_3)^*}
  =\overline{ \lrAngle{ [\id\otimes\id\otimes
      (S^2\circ *\circ S)] (\R^{-1}_{31}\R_{23}), h_3} } \notag\\
  &=\overline{ \lrAngle{ [\id\otimes\id\otimes
      (* \circ S^{-1})] (\R^{-1}_{31}\R_{23}), h_3} }
  =\overline{ \lrAngle{ (*\otimes * \otimes \id )
      (\R_{13} \R^{-1}_{32} ), h_3} } \notag\\
  &=\overline{ \lrAngle{ (*\otimes * \otimes \id )
      (\R^{-1}_{12} \R^{-1}_{32} \R_{13} \R_{12}), h_3} } \notag\\
  &=\overline{ \lrAngle{\R_{21}[ (*\otimes * \otimes \id) 
      (\R^{-1}_{32} \R_{13}) ]\R^{-1}_{21}, h_3} } 
  =\R_{21}\, \lrAngle{ \R^{-1}_{32}\,\R_{13} ,
    h_3}^{(*\otimes *)} \,\R^{-1}_{21} \notag\\ 
  &= \bigl(j(h)\bigr)^* \,.
\end{align}
We conclude that $j$ is a Hopf-$*$ algebra homomorphism from $\SUq^\op$
to the $q$-Lorentz algebra. Hence, $j(\SUq^\op)$ is indeed a Hopf-$*$
subalgebra of $\slC$. We will call it the subalgebra of the boosts.

\subsection{$L$-Matrices and the Explicit Form of the Boost Algebra}
\label{sec:L-Matrices}

To calculate the explicit form of the algebra of boosts we introduce
the computational tool of $L$-Matrices \cite{Faddeev:1990}. Let
$\rho^j$ be the representation map of the $D^j$-representation of
$\slq$. We define matrices of generators by applying $\rho^j$ to one
tensor factor of the universal $\R$-matrix $\R = \R_{[1]} \otimes
\R_{[2]}$,
\begin{xalignat}{2}
\label{eq:Lmatdef}
  (L^j_+)^a{}_b &:= \R_{[1]} \rho^j(\R_{[2]})^a{}_b \,, &
  (L^j_-)^a{}_b &:= \rho^j(\R^{-1}_{[1]})^a{}_b
    \,\R^{-1}_{[2]} \,.
\end{xalignat}
Here, we need the $L$-matrices for $j=\frac{1}{2}$, where we get
\begin{xalignat}{2}
  L_+^\frac{1}{2} &= \begin{pmatrix}
    K^{-\frac{1}{2}} &
    q^{-\frac{1}{2}}\lambda K^{-\frac{1}{2}}E \\
              0      & K^{\frac{1}{2}}   \end{pmatrix}, &
  L_-^\frac{1}{2} &= \begin{pmatrix}
   K^{\frac{1}{2}} & 0 \\ -q^{\frac{1}{2}}\lambda FK^{\frac{1}{2}}
   & K^{-\frac{1}{2}} \end{pmatrix}
\end{xalignat}
with respect to the basis $\{-,+\}$. The appearance of the square
roots of $K$ comes from the fact that $\R$ only exists as formal power
series.

We can derive some properties of the $L$-matrices from the properties
of $\R$: Applying $\id \otimes \rho^j \otimes \rho^j $ to the quantum
Yang-Baxter equation $\R_{12}\R_{13}\R_{23} = \R_{23}\R_{13}\R_{12}$
we obtain
\begin{equation}
\label{eq:LMat1}
  (L_+^j)^{a}{}_{c'} (L_+^j)^{d}{}_{d'} R^{c'd'}{}_{cd} =
  R^{ab}{}_{a'b'} (L_+^j)^{b'}{}_{d} (L_+^j)^{a'}{}_{c} 
\end{equation}
and in an analogous manner 
\begin{equation}
\begin{aligned}
\label{eq:LMat2}
  (L_-^j)^{a}{}_{c'} (L_-^j)^{b}{}_{d'} R^{c'd'}{}_{cd} &=
  R^{ab}{}_{a'b'} (L_-^j)^{b'}{}_{d} (L_-^j)^{a'}{}_{c} \\
  (L_-^j)^{a}{}_{c'} (L_+^j)^{b}{}_{d'} R^{c'd'}{}_{cd} &=
  R^{ab}{}_{a'b'} (L_+^j)^{b'}{}_{d} (L_-^j)^{a'}{}_{c} \,.
\end{aligned}
\end{equation}
From the coproduct properties $(\Delta \otimes \id )(\R) =
\R_{13}\R_{23}$, $(\id\otimes \Delta)(\R^{-1}) =
\R^{-1}_{12}\R^{-1}_{13}$ and from $(\varepsilon \otimes \id)(\R) = 1
= (\id\otimes \varepsilon)(\R^{-1})$ it follows that
\begin{xalignat}{2}
\label{eq:LMat3}
  \Delta \bigl( (L_\pm^j)^{a}{}_{c} \bigr ) &=
  (L_\pm^j)^{a}{}_{b} \otimes (L_\pm^j)^{b}{}_{c} \,, &
  \varepsilon\bigl( (L_\pm^j)^{a}{}_{b} \bigr) = \delta^a_b \,.
\end{xalignat}
Finally, we apply $\id \otimes \rho^j \otimes \id $ to the form
$\R^{-1}_{13}\R^{-1}_{23}\R_{12}\R_{13} = \R_{12}\R^{-1}_{23}$ of the
Yang-Baxter equation in order to get
\begin{equation}
\label{eq:LMat4}
  \R^{-1}\, (L_+^j)^{b}{}_{c}\otimes (L_-^j)^{a}{}_{b} \,\R
  = (L_+^j)^{a}{}_{b} \otimes (L_-^j)^{b}{}_{c} \,.
\end{equation}

Now, we can compute the explicit form of the boosts. Observing that
the dual pairing of $\SUq$ and $\suq$ (Sec.~\ref{sec:dualpairing}) can
be expressed on the matrix $M^a{}_b$ of generators of $\SUq$ by
$\lrAngle{h,M^a{}_b} = \rho^\frac{1}{2}(h)^a{}_b$, we get for the
boost generators
\begin{equation}
\begin{split}
  B^a{}_c &:= j(M^a{}_c) =
  \lrAngle{\R^{-1}_{31}\R_{23}, 1\otimes 1 \otimes M^a{}_c }
  =(\R^{-1}_{[2]} \otimes \R_{[1']})
  \rho^\frac{1}{2}(\R^{-1}_{[1]}\R_{[2']})^a{}_c \\
  &= \R^{-1}_{[2]}\rho^\frac{1}{2}(\R^{-1}_{[1]})^a{}_b \otimes
  \R_{[1']}\rho^\frac{1}{2}(\R_{[2']})^b{}_c
  = \bigl(L_-^\frac{1}{2}\bigr)^a{}_b \otimes
  \bigl(L_+^\frac{1}{2} \bigr)^b{}_c \,,
\end{split}
\end{equation}
explicitly,   
\begin{equation}
\label{eq:boostdef}
  B^{a}{}_{b}  =
  \begin{pmatrix}
    K^{\frac{1}{2}}\otimes K^{-\frac{1}{2}} &
    q^{-\frac{1}{2}}\lambda K^{\frac{1}{2}}
    \otimes K^{-\frac{1}{2}} E \\
    -q^{\frac{1}{2}}\lambda F K^{\frac{1}{2}}
    \otimes K^{-\frac{1}{2}} & \quad
    K^{-\frac{1}{2}}\otimes K^{\frac{1}{2}}
    - \lambda^2 F K^{\frac{1}{2}}\otimes K^{-\frac{1}{2}} E
  \end{pmatrix}
   =: \begin{pmatrix} a & b \\ c & d \end{pmatrix} \,.
\end{equation}
The commutation relations are
\begin{equation}
\label{eq:SUopRel}
\begin{gathered}
  ba=qab, \quad ca=qac, \quad db = q bd, \quad dc = q cd \\
  bc = cb, \quad da - ad = (q - q^{-1}) bc, \quad da-qbc=1 \,.
\end{gathered}
\end{equation}
The coproduct, $\Delta (B^{a}{}_{c}) = B^{a}{}_{b}\otimes
B^{b}{}_{c}$, is the same as for $\SUq$ just as the counit,
$\varepsilon (B^{a}{}_{b}) = \delta^{a}_{b}$.  For the antipode we had
$S^\op = S^{-1}$ and for the $*$-structure $*^\op := *\circ S^2$. Since
$(M^{a}{}_{b})^* = S(M^{b}{}_{a})$, it follows that
$(M^{a}{}_{b})^{*^\op} = S^{\op}(M^{b}{}_{a})$ in $\SUq^\op$ and,
consequently, the unitarity condition $(B^{a}{}_{b})^* =
S(B^{b}{}_{a})$ holds in $\slC$ as well. Written out this is
\begin{xalignat}{2}
\label{eq:StarSBoost}
  S\begin{pmatrix} a & b \\ c & d \end{pmatrix} &=
  \begin{pmatrix} d & -qb \\ -q^{-1}c & a \end{pmatrix}, &
  \begin{pmatrix} a & b \\ c & d \end{pmatrix}^* &=
  \begin{pmatrix} d & -q^{-1}c \\ -qb & a \end{pmatrix}.
\end{xalignat}
If we want to verify that the $B^{a}{}_{b}$ are the generators of a
$\SUq^\op$ subalgebra using the definition of the $q$-Lorentz algebra
only, we find that this is extremely tedious.

\subsection{Commutation Relations between Boosts and Rotations}

Now, we have to figure out the commutation relations between rotations
and boost, embedded into $\slC$ by the maps $i$ and $j$, respectively.
For $l\in \suq$ and $h\in\SUq^\op$ the embedding is
\begin{equation}
\begin{split}
  j(h)i(l)
  &= \lrAngle{\R^{-1}_{31}\R_{23}, h_3} ( l_{(1)}\otimes l_{(2)})\\ 
  &= \R^{-1}_{[2]} l_{(1)}\otimes \R_{[1']} l_{(2)}\,
  \lrAngle{ \R^{-1}_{[1]}\R_{[2']}l_{(3)}S(l_{(4)}), h} \\
  &= \R^{-1}_{[2]} l_{(1)}\otimes l_{(3)}\R_{[1']} \,
  \lrAngle{ \R^{-1}_{[1]} l_{(2)}\R_{[2']}S(l_{(4)}), h} \\ 
  &=  l_{(2)}\R^{-1}_{[2]}\otimes l_{(3)}\R_{[1']} \,
  \lrAngle{ l_{(1)}\R^{-1}_{[1]}\R_{[2']}S(l_{(4)}), h} \\
  &=  l_{(2)}\R^{-1}_{[2]}\otimes l_{(3)}\R_{[1']} \,
  \lrAngle{ l_{(1)}, h_{(1)} }
  \lrAngle{ \R^{-1}_{[1]}\R_{[2']}, h_{(2)}}
  \lrAngle{S(l_{(4)}), h_{(3)}} \\
  &=\lrAngle{ l_{(1)}, h_{(1)} } \, i(l_{(2)})j(h_{(2)}) \,
  \lrAngle{S(l_{(3)}), h_{(3)}} \,.
\end{split}
\end{equation}
The commutation relations which can be read off this equation are
precisely the ones of the quantum double
\cite{Drinfeld:1985,Drinfeld:1986}. For the generators they write out
\begin{align}
  B^a{}_b E &= E B^a{}_{a'} \rho^{\frac{1}{2}}(K^{-1})^{a'}{}_b
    + K \rho^{\frac{1}{2}}(E)^{a}{}_{a'} B^{a'}{}_{b'}
      \rho^{\frac{1}{2}}(K^{-1})^{b'}{}_b
    - B^{a}{}_{a'} \rho^{\frac{1}{2}}(EK^{-1})^{a'}{}_b \notag\\
  B^a{}_b F &= F \rho^{\frac{1}{2}}(K^{-1})^{a'}{}_{a'} B^{a'}{}_{b} 
    - K^{-1} \rho^{\frac{1}{2}}(K^{-1})^{a}{}_{a'} B^{a'}{}_{b'}
      \rho^{\frac{1}{2}}(KF)^{b'}{}_b
    + \rho^{\frac{1}{2}}(F)^{a}{}_{a'} B^{a'}{}_{b} \notag\\
  B^a{}_b K &= \rho^{\frac{1}{2}}(K)^{a}{}_{a'} B^{a'}{}_{b'}
      \rho^{\frac{1}{2}}(K^{-1})^{b'}{}_b \,.
\end{align}
Explicitly, this gives us 
\begin{equation}
\begin{aligned}
  \begin{pmatrix} a & b \\ c & d \end{pmatrix} E
  &= \begin{pmatrix} q E a - q^{\frac{3}{2}} b & q^{-1}Eb \\
    qEc+q^{\frac{3}{2}}Ka- q^{\frac{3}{2}}d &
    q^{-1}Ed+q^{-\frac{1}{2}}Kb \end{pmatrix} \\
  \begin{pmatrix} a & b \\ c & d \end{pmatrix} F
  &= \begin{pmatrix} q F a + q^{-\frac{1}{2}}c &
    qFb-q^{-\frac{1}{2}}K^{-1}a + q^{-\frac{1}{2}}d \\
    q^{-1}Fc & q^{-1}Fd-q^{-\frac{5}{2}}K^{-1}c \end{pmatrix} \\
  \begin{pmatrix} a & b \\ c & d \end{pmatrix} K
  &= K \begin{pmatrix} a & q^{-2}b \\ q^{2} c & d\end{pmatrix}
  \,,\quad\begin{pmatrix} a & b \\ c & d \end{pmatrix} K^{-1}
  = K^{-1} \begin{pmatrix} a & q^{2}b \\ q^{-2} c & d\end{pmatrix}\,.
\end{aligned}
\end{equation}
We summarize:
\begin{Definition}
The Hopf-$*$ algebra generated by $\SUq^\op$ and $\suq$ with cross
commutation relations
\begin{equation}
\label{eq:RotBoostCommute1}
  hl = \lrAngle{ l_{(1)}, h_{(1)} } \,l_{(2)} h_{(2)} \,
  \lrAngle{S(l_{(3)}), h_{(3)}}
\end{equation}
or, equivalently,
\begin{equation}
\label{eq:RotBoostCommute}
 lh = \lrAngle{ S(l_{(1)}), h_{(1)} } \,h_{(2)} l_{(2)} \,
  \lrAngle{l_{(3)}, h_{(3)}}
\end{equation}
for $h\in\SUq^\op$ and $l\in\suq$, is the quantum double form of the
$q$-Lorentz algebra \cite{Podles:1990}.
\end{Definition}
Finally, if we want to invert the embedding $i\otimes j:
\suq\otimes\SUq^\op \rightarrow \slC$ we find
\begin{xalignat}{2}
  E\otimes 1 &= q K^{-\frac{1}{2}}(Ea-q^{\frac{3}{2}}\lambda^{-1}b)\,, &
  1\otimes E &= q^{\frac{1}{2}}\lambda^{-1}a^{-1}b \notag \\
  F\otimes 1 &= -q^{-\frac{1}{2}}\lambda^{-1}ca^{-1}\,, &
  1\otimes F &= q K^{\frac{1}{2}}(Fa+q^{-\frac{1}{2}}\lambda^{-1}c) \\
  K\otimes 1 &= K^{\frac{1}{2}} a \,,&
  1\otimes K &= K^{\frac{1}{2}} a^{-1} \notag
\end{xalignat}
For these expressions to make sense we had to add the generator
$a^{-1}$ to $\SUq^\op$ and $K^{\pm\frac{1}{2}}$ to $\suq$. From the
viewpoint of representation theory this modification seems to be
insignificant.

\section{The Vectorial Form of the $q$-Lorentz Algebra}

\subsection{Tensor Operators of the $q$-Lorentz Algebra}

The definition of tensor operators in Eq.~\eqref{eq:Tensor1} has
been general. We just have to work it out for the $q$-Lorentz
algebra. We begin by calculating for $g\otimes h \in \slC$
\begin{multline}
  (g\otimes h)_{(1)}\otimes S\bigl((g\otimes h)_{(2)}\bigr)
  =(g_{(1)}\otimes \R^{-1}_{[1]} h_{(1)} \R_{[1']})
  \otimes S(\R^{-1}_{[2]} g_{(2)} \R_{[2']} \otimes h_{(2)}) \\
  =(g_{(1)}\otimes \R^{-1}_{[1]} h_{(1)} \R_{[1']})
  \otimes ( \R_{[2'']} S(\R^{-1}_{[2]} g_{(2)} \R_{[2']}) \R^{-1}_{[2''']}
  \otimes \R_{[1'']} S(h_{(2)})\R^{-1}_{[1''']} )\\
  =(g_{(1)}\otimes \R_{[1]} h_{(1)} \R_{[1']})
  \otimes ( \R_{[2'']} S(\R_{[2']}) S(g_{(2)}) \R_{[2]} \R^{-1}_{[2''']}
  \otimes \R_{[1'']} S(h_{(2)})\R^{-1}_{[1''']} ) \,, 
\end{multline}
where in the last step we have used that $(\id \otimes S)(\R^{-1}) =
\R$. Hence, for $T_{\mu\nu} = \sum_n A^n_{\mu\nu} \otimes
B^n_{\mu\nu}$ (no summation of $\mu$ and $\nu$) to be a
$D^{(i,j)}$-tensor operator of $\slC$
\begin{multline}
  T_{\mu'\nu'} \rho^i(g)^{\mu'}{}_\mu \rho^j(h)^{\nu'}{}_\nu
  =\adL(g\otimes h)\tr (T_{\mu\nu}) \\ 
  =\sum_n g_{(1)}A^n_{\mu\nu} \R_{[2'']}
  S(\R_{[2']}) S(g_{(2)}) \R_{[2]} \R^{-1}_{[2''']} 
  \otimes \R_{[1]} h_{(1)} \R_{[1']} B^n_{\mu\nu}
  \R_{[1'']} S(h_{(2)})\R^{-1}_{[1''']} 
\end{multline}
must hold for all $g\otimes h \in \slC$. 

Some tensor operators of $\slC$ can be derived from tensor operators
of $\suq$: If $A_\mu$ is a $D^j$-tensor operator of $\suq$ then $A_\mu
\otimes 1$ is a $D^{(j,0)}$-tensor operator. We check this by
inserting $T_{\mu\nu} = A_\mu\otimes 1$ in the last equation:
\begin{align}
  \adL &(g\otimes h)\tr (A_\mu\otimes 1) \notag\\
  &= g_{(1)}A_\mu S(\R_{[2']} \R^{-1}_{[2'']})
  S(g_{(2)}) \R_{[2]} \R^{-1}_{[2''']} 
  \otimes \R_{[1]} h_{(1)} (\R_{[1']} \R^{-1}_{[1'']})
  S(h_{(2)})\R^{-1}_{[1''']} \notag\\
  &= g_{(1)}A_\mu S(g_{(2)}) \R_{[2]} \R^{-1}_{[2''']} 
  \otimes \R_{[1]} h_{(1)} S(h_{(2)})\R^{-1}_{[1''']}
  = g_{(1)}A_\mu S(g_{(2)})\otimes \varepsilon(h) \notag\\
  &= (A_{\mu'} \otimes 1)\,\rho^j(g)^{\mu'}{}_\mu \,\rho^0(h) \,. 
\end{align}
In the same manner we verify that $\R_{21}(1\otimes A_\mu)
\R^{-1}_{21}$ is a $D^{(0,j)}$-tensor operator:
\begin{align}
  \adL &(g\otimes h)\tr
  \R_{21}(1\otimes A_\mu) \R^{-1}_{21} \notag\\
  &= g_{(1)} \R_{[2'']}
  S(\R_{[2']}) S(g_{(2)}) \R_{[2]} \R^{-1}_{[2''']} 
  \otimes \R_{[1]} h_{(1)} \R_{[1']} \R_{[1'']} A_\mu
  S(h_{(2)})\R^{-1}_{[1''']}\notag\\
  &= g_{(1)} S(g_{(2)}) \R_{[2]} \R^{-1}_{[2''']} 
  \otimes \R_{[1]} h_{(1)} A_\mu
  S(h_{(2)})\R^{-1}_{[1''']}\notag\\
  &= \R_{21}(1\otimes h_{(1)} A_{\mu}S(h_{(2)}) )
  \R^{-1}_{21} \varepsilon(g)\notag\\
  &= \R_{21}(1\otimes A_{\mu'} )
  \R^{-1}_{21} \,\varepsilon(g)\rho^j(h)^{\mu'}{}_\mu \,.
\end{align}

\subsection{The Vectorial Generators}

Now, it is obvious how we can define vectorial generators of the
$q$-Lorentz algebra. Let $J_A$ be the vector generator of $\suq$ as
defined in Eqs.~\eqref{eq:DefL}. We define\footnote{The operators $R$
  and $S$ defined here correspond to the operators $q^2[2] R$ and
  $-q^2[2]S$ of \cite{Rohregger:1999}.}
\begin{xalignat}{2}
\label{eq:RSdef}
  S_A &:= J_A\otimes 1 \,,&
  R_A &:= \R_{21}(1\otimes J_A) \R^{-1}_{21} \,.
\end{xalignat}
From the last section it is obvious that $S_A$ is a $D^{(1,0)}$-tensor
and $R_A$ is a $D^{(0,1)}$-tensor operator, that is, a left and right
chiral vector operator, respectively. Moreover, both $R_A$ and $S_A$
are vector operators with respect to rotations since $D^{(1,0)}$ and
$D^{(0,1)}$ induce a $D^1$ vector representation of the $\suq$
subalgebra.

We can raise the indices with the 3-metric of $\suq$ introduced in
Eq.~\eqref{eq:epsgdef}, $S^A = g^{AB}S_B$, giving us a congredient
vector operator,
\begin{equation}
\begin{split}
  \mathrm{ad_L}(g\otimes h)\tr S^A
  &= \mathrm{ad_L}(g\otimes h)\tr (J_{A'}\otimes 1)g^{AA'} \\ 
  &= (J^{B}\otimes 1) g^{AA'} g_{B'B}\,
    \rho^j(g)^{B'}{}_{A'}\,\varepsilon(h)\\
  &= S^{B} \rho^j(Sg)^{A}{}_{B}\, \varepsilon(Sh)\,,
\end{split}
\end{equation}
and the same for $R_A$. By looking at the definition of the
$*$-structure of $\slC$ we immediately see that
\begin{equation}
  (R_A)^* = S^A \,.
\end{equation}

For the commutation relations of the algebra generated by $R_A$ and
$S_A$ to close we yet have to embed the Casimir operator $W$ of the
vectorial form of $\suq$, as defined in Eq.~\eqref{eq:DefW}, in the
$q$-Lorentz algebra, that is\footnote{The operator $V$ defined here
  corresponds to $U'$ in \cite{Rohregger:1999}.}
\begin{xalignat}{2}
  V &:= W\otimes 1 \,,&
  U &:= \R_{21}(1\otimes W) \R^{-1}_{21} = 1\otimes W\,.
\end{xalignat}
By construction the commutation relations of the $R$'s and $U$ among
each other are the same as for the $L$'s and $W$ as given in
Eqs.~\eqref{eq:LL-Rel}. The same holds for the $S$'s and $V$ since
these generators are embedded by an inner automorphism.  To calculate
the commutation relations of $R_A$ with $S_B$ we first note that
commuting $\R_{21}$ with $1\otimes J_A$ shows us that
\begin{equation}
  R_A = \R_{[2]} \otimes J_{A'} \rho^1(\R_{[1]})^{A'}{}_A \,.
\end{equation}
Then we commute this expression with $S_A$
\begin{equation}
\begin{split}
  R_A S_B &= \R_{[2]} J_B \otimes J_{A'} \rho^1(\R_{[1]})^{A'}{}_A\\
  &= J_{B'} \rho^1(\R_{[2]})^{B'}{}_B \R_{[2']} \otimes 
  J_{A'} \rho^1(\R_{[1']}\R_{[1]})^{A'}{}_A \\
  &=S_{B'} R_{A'} \rho^1(\R_{[1]})^{A'}{}_A\rho^1(\R_{[2]})^{B'}{}_B\,.
\end{split}
\end{equation}
The representation of the universal $\R$-matrix appearing on the last
line is proportional to the $R$-matrix of
$SO_q(3)$, defined in Eq.~\eqref{eq:AppR4}. 
The $RS$-commutation relations can now be written as
\begin{equation}
  R_A S_B = q^2 S_{B'} R_{A'}\, R_{\mathrm{so}_3}^{A'B'}{}_{AB} \,,
\end{equation}
where $R_{\mathrm{so}_3}$ is given explicitly in
Eq.~\eqref{eq:AppR6}. We summarize
\begin{Definition}
\label{th:RS-Algebra}
  The algebra generated by $R_A$, $U$, $S_A$, $V$, where $A$ runs
  through $\{-,+,3\}$, with relations
\begin{subequations}
\label{eq:RRSS-Rel}
\begin{align}
\label{eq:RR-Rel}
  R_A R_B \varepsilon^{AB}{}_{C} &= UR_C \,,&
  R_A U &= U R_A \,,&
  U^2 - \lambda^2 g^{AB} R_A R_B &= 1 \\
\label{eq:SS-Rel}
  S_A S_B \varepsilon^{AB}{}_{C} &= VS_C \,,&
  S_A V &= V S_A \,,&
  V^2 - \lambda^2 g^{AB} S_A S_B &= 1
\end{align}
\begin{gather}
\label{eq:RS-Rel}
  R_C S_D = q^2 S_C R_D - q^{-1}\lambda\, g_{CD}(g^{AB}S_A R_B)
  + \varepsilon_C{}^X{}_D \varepsilon^{AB}{}_X\, S_A R_B \\
  R^A V = V R^A\,, \qquad UV=VU\,, \qquad S^A U = U S^A 
\end{gather}
and $*$-structure
\begin{xalignat}{2}
  R_A^* &=  g^{AB} S_B \,, &  U^* &= V 
\end{xalignat}
\end{subequations}
is called the vectorial or $RS$-form of the $q$-Lorentz algebra
\cite{Rohregger:1999}.
\end{Definition}

\subsection{Relations with the other Generators}
\label{sec:MainContrib1}

Let us first express the vectorial generators $R_A$ and $S_A$ by the
original generators of $\slC$. For $S_A$ and $V$ the case is simple.
We merely have to look up the expressions for $J_A$ and $W$ in
Eqs.~\eqref{eq:DefL} and~\eqref{eq:DefW}. For completeness we write
them down once more
\begin{equation}
\begin{aligned}
  S_{-} &:= q[2]^{-\frac{1}{2}}KF \otimes 1 \\
  S_3    &:= [2]^{-1} (q^{-1}EF-qFE) \otimes 1 \\
  S_{+} &:= -[2]^{-\frac{1}{2}}E \otimes 1\\
  V &:= [K - \lambda [2]^{-1} (q^{-1}EF-qFE)] \otimes 1 \,.
\end{aligned}
\end{equation}
For $R_A$ one might at first sight expect formal power
series, but as we have shown in the preceding section 
\begin{equation}
  R_A = \R_{[2]} \otimes L_{A'} \,\rho^1(\R_{[1]})^{A'}{}_A
  = S\bigl[ (L^1_-)^{A'}{}_A \bigr] \otimes J_{A'} \,.
\end{equation}
We only have to sum over the $L^1_-$-matrix of $\suq$, which has been
computed in Eq.~\eqref{eq:AppR5} where we get 
\begin{equation}
  S\bigl[ (L^1_-)^{A}{}_B \bigr] =
  \begin{pmatrix}
    K^{-1} & 0 & 0 \\
    \lambda [2]^\frac{1}{2} F & 1 & 0 \\
    q^2\lambda^2 K F^2 & q\lambda [2]^\frac{1}{2} KF & K \\
  \end{pmatrix}
\end{equation}
with respect to the $\{-1,0,1\}=\{-,3,+\}$ basis, so the expressions
for the $R$'s become
\begin{equation}
\begin{aligned}
  R_{-} &= q[2]^{-\frac{1}{2}}K^{-1} \otimes KF
  + \lambda[2]^{-\frac{1}{2}} F\otimes (q^{-1}EF-qFE) \\
  &\qquad - q^{2}\lambda^2 [2]^{-\frac{1}{2}} KF^2\otimes E\\
  R_3    &= 1\otimes [2]^{-1}(q^{-1}EF-qFE) 
  - q \lambda KF\otimes E \\
  R_{+} &= -[2]^{-\frac{1}{2}} K\otimes E \\
  U &= 1\otimes [K - \lambda [2]^{-1} (q^{-1}EF-qFE)] \,.
\end{aligned}
\end{equation} 
Next, let us express $R_A$ and $S_A$ by the generators of the quantum
double form of the $q$-Lorentz algebra. For $S_A$ we find 
\begin{equation}
\begin{aligned}
  S_{-} &= -q^{-\frac{1}{2}}\lambda^{-1}[2]^{-\frac{1}{2}}
           K^{\frac{1}{2}}c \\
  S_3 &= q^{-\frac{3}{2}}\lambda^{-1}[2]^{-1}
         K^{-\frac{1}{2}}(qcE-Ec) \\
  &= \lambda^{-1}[2]^{-1}K^{-\frac{1}{2}}
      (q^{-\frac{1}{2}}\lambda Ec +qKa -qd)\\
  S_{+} &= q[2]^{-\frac{1}{2}} \,K^{-\frac{1}{2}}
    (q^{\frac{3}{2}}\lambda^{-1} b- Ea)\\
  V &= [2]^{-1}K^{-\frac{1}{2}}
       (q^{-1} K a - q^{-\frac{1}{2}}\lambda Ec + qd)  \,.
\end{aligned}
\end{equation}
To compute the corresponding expressions for $R_A$ we remember that
$S_{-}^* = -q^{-1}R_{+}$, $S_3^* = R_3$, and $S_{+}^* = -qR_{-}$.
With the $*$-structure of rotations and boosts as given in
Eqs.~\eqref{eq:LL-Star} and~\eqref{eq:StarSBoost} this yields
\begin{equation}
\begin{aligned}
  R_{-} &= [2]^{-\frac{1}{2}} \,K^{-\frac{1}{2}}
    (Fd + q^{-\frac{5}{2}}\lambda^{-1} c)\\
  R_3 &= q^{-\frac{1}{2}}\lambda^{-1}[2]^{-1}
         K^{\frac{1}{2}}(bF - q^3 Fb) \\
  &= \lambda^{-1}[2]^{-1}K^{\frac{1}{2}}
     (-q^{\frac{3}{2}}\lambda Fb - q^{-1}K^{-1} a +q^{-1}d)\\
  R_{+} &= -q^{\frac{1}{2}}\lambda^{-1}[2]^{-\frac{1}{2}}
           K^{\frac{1}{2}}b \\
  U &= [2]^{-1}K^{\frac{1}{2}}
        (q^{\frac{3}{2}}\lambda Fb + q^{-1}K^{-1}a + q d) \,.
\end{aligned}
\end{equation}
We also want to express the generators of boosts and rotations within
the $RS$-algebra. For the vectorial generators of the rotations we
find \cite{Lorek:1997a,Cerchiai:1998}
\begin{equation}
\label{eq:Rot-RS}
\begin{aligned}
  J_C &= VR_C + US_C +q\lambda R_A S_B \,\varepsilon^{AB}{}_C \\
  W   &= UV + q^2\lambda^2 g^{AB} \,R_A S_B \,.
\end{aligned}
\end{equation}
While this yields an expression of $K = W +\lambda J_3$, $K^{-1}$ is
not a member of the $RS$-algebra proper. We must add
$K^{-\frac{1}{2}}$ by hand to the $RS$-algebra to write down expressions for
the boosts
\begin{subequations}
\begin{xalignat}{2}
  a &= K^{-\frac{1}{2}}(V +\lambda S_3) \,, &
  b &= -q^{-\frac{1}{2}}\lambda[2]^{\frac{1}{2}} K^{-\frac{1}{2}} R_+ \\
  c &= -q^{\frac{1}{2}}\lambda[2]^{\frac{1}{2}} K^{-\frac{1}{2}} S_- \,,&
  d &= K^{-\frac{1}{2}}(U +\lambda R_3) \,.
\end{xalignat}
\end{subequations}

\chapter{Algebraic Structure of the $q$-Poincar{\'e} Algebra}
\label{sec:PoincareStructure}

\section{The $q$-Poincar{\'e} Algebra}

\subsection{Construction of the $q$-Minkowski-Space Algebra}
\label{sec:MinkowskiConstruct}

As in the undeformed case, we want to construct the coordinate
functions of Minkowski space to form a matrix $X_{a\dot{b}}$ with a
lower undotted and dotted index. For the cotransformations to be
compliant with the $*$-structure the $*$ has to act on $X_{a\dot{b}}$
as on a product $\phi_a\psi_{\dot{b}}$, that is, $(X_{a\dot{b}})^* :=
X_{b\dot{a}}$. For our purposes it is more convenient to work with the
index structure $X_a{}^{\dot{b}}$,
\begin{xalignat}{2}
\label{eq:MinkowskiStar}
  X_a{}^{\dot{b}} &:= \begin{pmatrix} A & B \\ C & D \end{pmatrix} &
  (X_a{}^{\dot{b}})^*
  &=\begin{pmatrix} -qD & B \\ C & -q^{-1}A \end{pmatrix}\,,
\end{xalignat}
with respect to the $\{-, +\}$ basis. With this index structure the
cotransformation is\footnote{Recall, that we think of the dot as
  belonging to $X$ rather than to the index itself.}
\begin{equation}
  \rho_\mathrm{R}(X_a{}^{\dot{b}}) = X_{a'}{}^{\dot{b}'} \otimes
  (M^{a'}{}_{a}\otimes M^{b'}{}_{b}) \,.
\end{equation}
Upon dualizing, this right coaction of the $q$-Lorentz group becomes a
left $D^{(\frac{1}{2},\frac{1}{2})}$ action of the $q$-Lorentz algebra.

We want to construct the space algebra out of the algebra
$\mathbb{C}\lrAngle{ X_a{}^{\dot{b}} }$ freely generated by the
generators $X_a{}^{\dot{b}}$ divided by some relations. The generators
have the dimension of a length, so we need homogeneous
relations,\footnote{For inhomogeneous relations we would need to
  introduce an additional dimensional parameter.} which for the
correct undeformed limit have to be of second order. We demand the
resulting quotient algebra to be a $q$-Lorentz module algebra.

This last requirement means that the quadratic terms that will be set
zero must be the basis of a $q$-Lorentz submodule. For only if we
divide the free module $\mathbb{C}\lrAngle{ X_a{}^{\dot{b}} }$ by an
ideal generated by a submodule, the quotient will be a module again.
The vector space generated by $X_a{}^{\dot{b}}X_c{}^{\dot{d}}$
naturally forms a $D^{(\frac{1}{2},\frac{1}{2})} \otimes
D^{(\frac{1}{2},\frac{1}{2})}$ representation of the $q$-Lorentz
algebra. By the Clebsch-Gordan-Series~\eqref{eq:CGSeries1} this
representation has the same four subrepresentations as in the
undeformed case. To obtain the correct undeformed limit where the
space-functions commute, it is the submodules $D^{(1,0)}$ and
$D^{(0,1)}$ that have to be set zero. The bases of those two
submodules as computed in Eqs.~\eqref{eq:AppMinkbasis1}
and~\eqref{eq:AppMinkbasis2} yield the relations
\begin{equation}
\begin{aligned}
  0 &= qBA - q^{-1}AB \\
  0 &= DA - AD + \lambda BB + BC - q^{-2}CB  \\
  0 &= DC - CD + \lambda DB \\
  0 &= CA - AC + \lambda BA \\
  0 &= DA - AD + \lambda BB + CB - q^{-2}BC  \\
  0 &= qDB - q^{-1}BD \,,
\end{aligned}
\end{equation}
which can be written more compactly as
\begin{equation}
\label{eq:Minkowski1}
\begin{gathered}
  AB = q^2 BA \,,\qquad BD = q^2 DB \,,\qquad BC=CB\\
  AC-CA= \lambda BA \,,\qquad CD-DC = \lambda DB \\
  AD - DA = \lambda B( B + q^{-1}C)\,.
\end{gathered}
\end{equation}
Now we can give the definition of the $q$-Minkowski-Space Algebra.
\begin{Definition}
  The $*$-algebra generated by $\{A,B,C,D \}$ with $*$-structure as in
  Eq.~\eqref{eq:MinkowskiStar} and commutation
  relations~\eqref{eq:Minkowski1} is called the $q$-Minkowski-Space
  algebra $\Mink$ \cite{Carow-Watamura:1990}.
\end{Definition}
The basis vector of the $D^{(0,0)}$ submodule yields a $q$-Lorentz
scalar, that corresponds to the invariant quadratic length, $X^2$, of
Minkowski space. Up to normalization we get from
Eq.~\eqref{eq:AppMinkbasis3} 
\begin{equation}
\label{eq:Minkowski2}
  X^2 := [2]^{-1}(qDA+q^{-1}AD-q^{-1}BC-q^{-1}CB-q^{-1}\lambda BB) \,,
\end{equation}
which can be simplified with the commutation
relations~\eqref{eq:Minkowski1} to
\begin{equation}
  X^2 = DA - q^{-2}BC \,.
\end{equation}
It turns out that this expression commutes with all generators of
$\Mink$. Hence, it can be viewed as the length Casimir of
$q$-Minkowski space or, within a momentum representation, as mass
Casimir of the $q$-Poincar{\'e} algebra.

\subsection{4-Vectors and the $q$-Pauli Matrices}
\label{sec:FourvectorsPauli}

We have constructed the $q$-Lorentz algebra to possess a $\suq$
Hopf-$*$ subalgebra, viewed as the algebra of rotations. Hence, we are
able to write $X_{a}{}^{\dot{b}}$ in a manifest 4-vector form, that
is, split up its 4 degrees of freedom with respect to rotations into a
scalar and a 3-vector.

The $D^{(\frac{1}{2},\frac{1}{2})}$ representation induces a
representation on the subalgebra of rotations. To compute the
representation map $\rho$ of the latter we have to embed $\suq$ with
$i=\Delta$ and then apply the representation map
$\rho^{(\frac{1}{2},\frac{1}{2})}$ yielding
\begin{equation}
  \rho = (\rho^\frac{1}{2}\otimes\rho^\frac{1}{2})\circ \Delta \,.
\end{equation}
In other words, this induced representation is simply the tensor
representation $D^\frac{1}{2}\otimes D^\frac{1}{2}$ which reduces
according to the Clebsch-Gordan series
\begin{equation}
  D^\frac{1}{2}\otimes D^\frac{1}{2} \cong D^0\oplus D^3
\end{equation}
to the direct sum of a scalar and a vector representation.
Explicitly, this reduction of $X_a{}^{\dot{b}}$ into a 4-vector is
expressed by the $q$-Clebsch-Gordan coefficients,
\begin{equation}
\label{eq:Paulitransform}
  X_0 = q^{-1}[2]^{-\frac{1}{2}}
  \CGC{\tfrac{1}{2}}{\tfrac{1}{2}}{0}{a}{b}{0}X_a{}^{\dot{b}} \,,\quad
  X_C = [2]^{-\frac{1}{2}}
  \CGC{\tfrac{1}{2}}{\tfrac{1}{2}}{1}{a}{b}{C} X_a{}^{\dot{b}} \,,
\end{equation}
where $C$ runs through $(-1,0,1)=(-,3,+)$ and we sum over repeated
indices. The factor $[2]^{-\frac{1}{2}}$ has been introduced to ensure
the right undeformed limit, the factor of $q^{-1}$ in the definition
of $X_0$ is traditional \cite{Lorek:1997a}. Written out, we get
\begin{equation}
\label{eq:FourVectorDef}
\begin{aligned}
  X_0 &= q^{-1}[2]^{-1}(q^{\frac{1}{2}}C - q^{-\frac{1}{2}}B) \\
  X_- &= [2]^{-\frac{1}{2}} A\\
  X_+ &= [2]^{-\frac{1}{2}} D\\
  X_3 &= [2]^{-1}(q^{-\frac{1}{2}}C + q^{\frac{1}{2}}B) \,.
\end{aligned}
\end{equation}
The back transformation is
\begin{subequations}
\label{eq:AX-Trans}
\begin{xalignat}{2}
  A &= [2]^{\frac{1}{2}} X_{-}\,, & 
  B &= q^{\frac{1}{2}}(X_{3} - X_{0}) \\
  C &= q^{-\frac{1}{2}}X_{3} + q^{\frac{3}{2}}X_{0} \,, &
  D &= [2]^{\frac{1}{2}} X_{+} \,.
\end{xalignat}
\end{subequations}
Expressed in terms of the $4$-vector generators, the commutation
relations~\eqref{eq:Minkowski1} become
\begin{equation}
\label{eq:XX-Rel1}
\begin{gathered}
  X_- X_0 = X_0 X_- \,,\quad X_+ X_0 = X_0 X_+ \,,\quad
  X_3 X_0 = X_0 X_3 \\
  q^{-1} X_- X_3 - q X_3 X_- = -\lambda X_- X_0\,,\quad
  q^{-1} X_3 X_+ - q X_+ X_3 = -\lambda X_+ X_0 \\ 
  X_- X_+ - X_+ X_- - \lambda X_3 X_3 = -\lambda X_3 X_0
\end{gathered}
\end{equation}
Using the $q$-deformed $\varepsilon$-tensor~\eqref{eq:epsgdef} this
can be written more compactly as
\begin{xalignat}{2}
\label{eq:XX-Rel2}
  X_0 X_A &= X_A X_0 \,,&
  X_A X_B\,\varepsilon^{AB}{}_{C} &= -\lambda X_0 X_C \,.
\end{xalignat}
For the $*$-structure we get
\begin{xalignat}{2}
  X_0^* &= X_0 \,,& (X_A)^* &= X^A \,,
\end{xalignat}
for the scalar product~\eqref{eq:Minkowski2}
\begin{equation}
  X^2 = X_0^2  + q^{-1} X_- X_+ +  q X_+ X_- - X_3^2
  = X_0^2 - X_A X_B g^{AB}\,.
\end{equation}
From this, we can read off the 4-metric, $X^2 = X_\mu X_\nu
\eta^{\mu\nu}$, with
\begin{xalignat}{2}
\label{eq:MinkowskiMetric}
  \eta^{00}&= 1 \,,& \eta^{AB} &= - g^{AB} 
\end{xalignat}
and zero otherwise. We also could have computed the metric directly
from the formulas of the Clebsch-Gordan coefficients.

If we write the back transformation~\eqref{eq:AX-Trans} as
\begin{equation}
\label{eq:Paulibacktransform}
  X_a{}^{\dot{b}} = \sum_\mu X_\mu (\sigma_\mu)_a{}^{\dot{b}} \,,
\end{equation}
this defines the $q$-Pauli matrices
\begin{equation}
  (\sigma_0)_a{}^{\dot{b}}  = q[2]^{\frac{1}{2}}
         \CGC{\tfrac{1}{2}}{\tfrac{1}{2}}{0}{a}{b}{0}\,,\quad
  (\sigma_C)_a{}^{\dot{b}}  = [2]^{\frac{1}{2}}
         \CGC{\tfrac{1}{2}}{\tfrac{1}{2}}{1}{a}{b}{C} \,.
\end{equation}
For the usual index structure we have to lower the dotted index.
\begin{equation}
  (\sigma_\mu)_{a\dot{b}} =
  (\sigma_\mu)_a{}^{\dot{b}'}\varepsilon_{b'b}
\end{equation}
The $q$-Pauli matrices with lower undotted and dotted indices are
\begin{xalignat}{4}
\label{eq:qPauli}
  \sigma_0 &=
    \begin{pmatrix} q & 0 \\ 0 & q \end{pmatrix}, &
  \sigma_- &= {[2]}^{\frac{1}{2}}
    \begin{pmatrix} 0 & q^{-\frac{1}{2}} \\ 0 & 0 \end{pmatrix}, &
  \sigma_+ &= [2]^{\frac{1}{2}}
    \begin{pmatrix} 0 & 0 \\ -q^{\frac{1}{2}} & 0 \end{pmatrix},&
  \sigma_3 &=
    \begin{pmatrix} -q & 0 \\ 0 & q^{-1} \end{pmatrix} 
\end{xalignat}
with respect to the basis $\{-,+\}$.  If we compare the $q$-Pauli
matrices with the spin-$\frac{1}{2}$ representation of $J_A$ we find
\begin{equation}
  \rho^\frac{1}{2}(J_A) = [2]^{-1} \sigma\!_A \,.
\end{equation}
This tells us that if we raise (and lower) the vector index of
$\sigma\!_A$ as usual by $\sigma^A := g^{AA'} \sigma\!_{A'}$ we get
$(\sigma\!_A)^\dagger = \sigma^A$, that is,
\begin{equation}
  \overline{(\sigma\!_A)_{b\dot{a}}} = (\sigma^A)^{a\dot{b}} \,.
\end{equation}
From Eq.~\eqref{eq:LL-Rel} we deduce
\begin{equation}
  \sigma\!_A\, \sigma_B \,\varepsilon^{AB}{}_C =
  [4][2]^{-1}\, \sigma_C \,.
\end{equation}
Further relations which are not representations of relations within
the algebra of rotations can be found by explicit calculations
\begin{xalignat}{2}
  \sigma^A \sigma^B &= \varepsilon^{BA}{}_C\,\sigma^C + g^{BA}\,, &
  \sigma\!_A \sigma_B &=\sigma_C\, \varepsilon_A{}^C{}_B + g_{AB} \,.
\end{xalignat}

The basis transformation from the matrix generators $X_a{}^{\dot{b}}$
to the 4-vector generators $X_\mu$ defines a matrix representation
$\Lambda$ of the $q$-Lorentz algebra by 
\begin{equation}
\label{eq:FourVectorRep}
  (g\otimes h)\tr X_\mu = X_{\mu'}\, \Lambda(g\otimes h)^{\mu'}{}_\mu \,.
\end{equation}
Using the formulas for the basis transformation
Eqs.~\eqref{eq:Paulitransform} and \eqref{eq:Paulibacktransform} we
get
\begin{equation}
\begin{split}
  (g\otimes h) \tr X_0
  &= (g\otimes h)\tr \bigl( q^{-2}[2]^{-1}
     (\sigma_0)_a{}^{\dot{b}} X_a{}^{\dot{b}} \bigr) \\
  &= q^{-2}[2]^{-1}(\sigma_0)_a{}^{\dot{b}} \rho^{\frac{1}{2}}(g)^{a'}{}_{a}
     \,\rho^{\frac{1}{2}}(h)^{b'}{}_{b} X_{a'}{}^{\dot{b}'} \\
  &= q^{-2}[2]^{-1}(\sigma_0)_a{}^{\dot{b}} \rho^{\frac{1}{2}}(g)^{a'}{}_{a}
  \,\rho^{\frac{1}{2}}(h)^{b'}{}_{b} (\sigma_{\mu})_{a'}{}^{\dot{b}'} X_\mu\\
  &= X_\mu \,\Lambda(g\otimes h)^\mu{}_0
\end{split}
\end{equation}
and
\begin{equation}
\begin{split}
  (g\otimes h) \tr X_A
  &= (g\otimes h)\tr\bigl([2]^{-1}
     (\sigma\!_A)_a{}^{\dot{b}} X_a{}^{\dot{b}} \bigr)\\
  &= [2]^{-1}(\sigma\!_A)_a{}^{\dot{b}} \rho^{\frac{1}{2}}(g)^{a'}{}_{a}
     \,\rho^{\frac{1}{2}}(h)^{b'}{}_{b} X_{a'}{}^{\dot{b}'} \\
  &= [2]^{-1}(\sigma\!_A)_a{}^{\dot{b}} \rho^{\frac{1}{2}}(g)^{a'}{}_{a}
  \,\rho^{\frac{1}{2}}(h)^{b'}{}_{b} (\sigma_{\mu})_{a'}{}^{\dot{b}'} X_\mu \\
  &= X_\mu \,\Lambda(g\otimes h)^\mu{}_A \,,
\end{split}
\end{equation}
for any $(g\otimes h)\in\slC$. From this we can read off explicit
formulas for $\Lambda$ in terms of the
$D^{\frac{1}{2}}$-representation of $\suq$ and the $q$-Pauli matrices
\begin{equation}
\begin{aligned}
  \Lambda(g\otimes h)^\mu{}_0 &=
  q^{-2}[2]^{-1}(\sigma_0)_a{}^{\dot{b}} \rho^{\frac{1}{2}}(g)^{a'}{}_{a}
  \,\rho^{\frac{1}{2}}(h)^{b'}{}_{b} (\sigma_{\mu})_{a'}{}^{\dot{b}'}\\
  \Lambda(g\otimes h)^\mu{}_A &=
  [2]^{-1}(\sigma\!_A)_a{}^{\dot{b}} \rho^{\frac{1}{2}}(g)^{a'}{}_{a}
  \,\rho^{\frac{1}{2}}(h)^{b'}{}_{b} (\sigma_{\mu})_{a'}{}^{\dot{b}'} \,.
\end{aligned}
\end{equation}
The matrices representing the generators of rotations and boosts have
been calculated explicitly in Eqs.~\eqref{eq:AppFourVector1}
and~\eqref{eq:AppFourVector2}.

\subsection{Commutation Relations of the $q$-Poincar{\'e} Algebra}

In order to construct the $q$-Poincar{\'e} algebra we have to view $\Mink$
as the algebra of translations, so we write $P_\mu$ instead of
$X_\mu$. By construction $\Mink$ is a left $\slC$-module $*$-algebra.
Denoting the action of $h\in\slC$ on $p\in\Mink$ by $h\tr p$ this
means
\begin{xalignat}{2}
  h\tr pp' &= (h_{(1)}\tr p)(h_{(2)}\tr p')\,, &
  (h \tr p)^* &= (Sh)^* \tr p^* \,.
\end{xalignat}
As in the undeformed case, $\slC$ and $\Mink$ can then be joined
together in a semidirect product:
\begin{Definition}
\label{th:Poincare}
  The $*$-algebra of the Hopf semidirect product $\Mink \rtimes \slC$,
  that is, the vector space $\Mink \otimes \slC$ with multiplication
  \begin{equation}
  \label{eq:PoincCommute}
    (p\otimes h)(p'\otimes h') :=
    p(h_{(1)}\tr p')\otimes h_{(2)}h'
  \end{equation}
  and $*$-structure $(p\otimes h)^* = (1\otimes h^*)(p^*\otimes 1)$, is
  $\Poin$, the $q$-Poincar{\'e} algebra.
\end{Definition}
By construction we have
\begin{equation}
  (\adL h)\tr P_\mu = h_{(1)} P_\mu S(h_{(2)}) = h \tr P_\mu
  = P_{\mu'} \Lambda(h)^{\mu'}{}_\mu
\end{equation}
for all $h\in\slC$, that is, $P_\mu$ a 4-vector operator.

We want to calculate the commutation relations between $q$-Lorentz
generators and momenta explicitly. By construction of the 4-vectors
the zero component $P_0$ commutes with all rotations. According to
Eq.~\eqref{eq:sutensorop} we get for the 3-vector part
\begin{equation}
\begin{aligned}
  E P_A &=  P_A E + q^{(A+1)} \sqrt{[A+2][1-A]}\,\,P_{A+1} K\\
  F P_A &=  q^{-2A} P_A F +  q^{-A} \sqrt{[1+A][2-A]}\,\,P_{A-1}\\
  K P_A &= q^{2A} P_A K \,,
\end{aligned}
\end{equation}
where $A$ runs through $\{-1,0,1\}=\{-,3,+\}$. In terms of the
vectorial generators this becomes
\begin{equation}
\begin{aligned}
  J_A P_B &=  P_A J_B
  - \varepsilon_A{}^C{}_B \varepsilon^{DE}{}_C \, P_D J_E
  + \varepsilon_A{}^C{}_B \,P_C W \\
  W P_A &=  (\lambda^2 +1 ) P_A W
  - \lambda^2 \varepsilon^{BC}{}_A \, P_B J_C \,.
\end{aligned}
\end{equation}
For the commutation relations between momenta and boosts we use
Eq.~\eqref{eq:AppBoostAct} to write in an obvious matrix notation
\begin{equation}
\begin{aligned}
  \begin{pmatrix} a & b \\ c & d \end{pmatrix}  P_0 &=
    \begin{pmatrix}
      [2]^{-1}\left(\frac{[4]}{[2]}P_0 + q^{-1}\lambda P_3 \right)  &
      q^{-\frac{1}{2}}\lambda [2]^{-\frac{1}{2}} P_+ \\
      -q^{\frac{1}{2}}\lambda [2]^{-\frac{1}{2}} P_- &
      [2]^{-1}\left( \frac{[4]}{[2]} P_0 - q \lambda P_3 \right)
    \end{pmatrix}
    \begin{pmatrix}a&b\\c&d\end{pmatrix} \\
  \begin{pmatrix} a & b \\ c & d \end{pmatrix} P_- &=
    \begin{pmatrix} P_- \quad &
      q^{-\frac{1}{2}}\lambda [2]^{-\frac{1}{2}}(P_3 - P_0) \\
      0 & P_-
    \end{pmatrix}
    \begin{pmatrix}a&b\\c&d\end{pmatrix} \\
  \begin{pmatrix} a & b \\ c & d \end{pmatrix} P_+ &=
    \begin{pmatrix} P_+ & 0 \\
      -q^{\frac{1}{2}}\lambda [2]^{-\frac{1}{2}}(P_3 - P_0) & \quad P_+
    \end{pmatrix}
    \begin{pmatrix}a&b\\c&d\end{pmatrix} \\
  \begin{pmatrix} a & b \\ c & d \end{pmatrix} P_3 &=
    \begin{pmatrix} [2]^{-1}(2P_3 + q\lambda P_0)&
      q^{-\frac{1}{2}}\lambda [2]^{-\frac{1}{2}} P_+ \\
      -q^{\frac{1}{2}}\lambda [2]^{-\frac{1}{2}} P_- & \quad
      [2]^{-1}(2P_3 - q^{-1} \lambda P_0)
    \end{pmatrix}
    \begin{pmatrix}a&b\\c&d\end{pmatrix}.
\end{aligned}
\end{equation}
The commutation relations between momenta and the vectorial
$RS$-generators as defined in Eq.~\eqref{eq:RSdef} are more
complicated but involve only 3-vectors and scalars with respect to
rotations.
\begin{subequations}
\label{eq:RP-Rel}
\begin{align}
\label{eq:RP-Rel-1}
  R_C P_0 &= [4][2]^{-2} P_0 R_C +
    \lambda[2]^{-1}\varepsilon^{AB}{}_C\, P_A R_B - q [2]^{-1} P_C\, U \\
  S_C P_0 &= [4][2]^{-2} P_0 S_C +
    \lambda[2]^{-1}\varepsilon^{AB}{}_C\,P_A S_B + q^{-1}[2]^{-1} P_C\, U \\  
    &\notag\\
  \!\!\!\! R_C P_D &= q P_C R_D
    - \lambda [2]^{-1} \varepsilon_C{}^B{}_D \, P_0 R_B 
    - q^{-1}\lambda [2]^{-1} g_{CD}\, (g^{AB} P_A R_B ) \notag \\
\label{eq:RP-Rel-2} &\quad 
    - 2[2]^{-1}\varepsilon_C{}^X{}_D \varepsilon^{AB}{}_X\, P_A R_B
    - q^{-1}[2]^{-1}  g_{CD}\, P_0 U 
    + [2]^{-1} \varepsilon_C{}^A{}_D \,P_A U \\
  S_C P_D &= q P_C S_D
    - \lambda [2]^{-1} \varepsilon_C{}^B{}_D \, P_0 S_B 
    + q\lambda [2]^{-1} g_{CD}\, (g^{AB} P_A S_B ) \notag \\
    &\quad
    - 2[2]^{-1}\varepsilon_C{}^X{}_D \varepsilon^{AB}{}_X\, P_A S_B
    + q[2]^{-1}  g_{CD} \, P_0 V 
    + [2]^{-1} \varepsilon_C{}^A{}_D \,P_A V \\
 &\notag\\
\label{eq:UP-Rel-1}
  U P_0 &=  [4][2]^{-2} P_0 U - q^{-1}\lambda^2[2]^{-1}(g^{AB} P_A R_B ) \\
  V P_0 &=  [4][2]^{-2} P^0 V + q \lambda^2[2]^{-1}(g^{AB} P_A S_B ) \\
  &\notag\\
\label{eq:UP-Rel-2}
  U P_C &=  [4][2]^{-2} P_A U - q \lambda^2[2]^{-1} P_0 R_A
     - \lambda^2[2]^{-1} \,\varepsilon^{AB}{}_{C}\, P_A R_B\\
  V P_C &=  [4][2]^{-2} P_A V + q^{-1} \lambda^2[2]^{-1} P_0 S_A
     - \lambda^2 [2]^{-1}\,\varepsilon^{AB}{}_{C}\, P_A S_B
\end{align}
\end{subequations}
Finally, we want to indicate how one can boost 4-vector operators. Let
$V_0$ be some element of the $q$-Poincar{\'e} algebra. If we assume that
$V_0$ is the zero component of a left 4-vector operator the action of
the boosts on $V_0$ must be the same as on $P_0$, so according to
Eq.~\eqref{eq:AppBoostAct} we must define the other components by
\begin{equation}
\label{eq:Boost1}
\begin{aligned}
 V_- &:= \adL(-q^{-\frac{1}{2}}\lambda^{-1}[2]^{\frac{1}{2}} \,c) \tr V_0 \\ 
 V_+ &:= \adL(q^{\frac{1}{2}}\lambda^{-1}[2]^{\frac{1}{2}} \,b) \tr V_0 \\ 
 V_3 &:= \adL(\lambda^{-1}\,(d-a)) \tr V_0 \,.
\end{aligned}
\end{equation}
We will make use of this method of computing 4-vectors in
Sec.~\ref{sec:gammaboost} in order to compute the $\gamma$-matrices.
In case we know the zero component $V_{\tilde{0}}$ of a \emph{right}
4-vector the other components must be defined by
\begin{equation}
\label{eq:Boost1b}
\begin{aligned}
 V_{\tilde{-}} &:= V_{\tilde{0}} \tl\adR
        (-q^{\frac{1}{2}}\lambda^{-1}[2]^{\frac{1}{2}} \,b)\\
 V_{\tilde{+}} &:= V_{\tilde{0}} \tl\adR
        (q^{-\frac{1}{2}}\lambda^{-1}[2]^{\frac{1}{2}} \,c)\\
 V_{\tilde{3}} &:= V_{\tilde{0}} \tl\adR (\lambda^{-1}\,(d-a))\,.
\end{aligned}
\end{equation}

\section{The $q$-Pauli-Lubanski Vector and the Spin Casimir}
\label{sec:MainContrib2a}

\subsection{The $q$-Euclidean Algebra}
\label{sec:Euk}

Rotations and translations generate a $*$-subalgebra of the
$q$-Poincar{\'e} algebra, the $q$-Euclidean subalgebra $\Euk$. Since
rotations form a $\suq$ Hopf subalgebra of $\slC$ this $q$-Euclidean
subalgebra is a semidirect product
\begin{equation}
  \Euk = \Mink \rtimes \suq \,. 
\end{equation}
By comparing Eq.~\eqref{eq:XX-Rel2} with Eq.~\eqref{eq:LL-Rel} we note
that $\Mink$ and $\suq$ are very similar as algebras. One could
identify the generators by a map $\xi: \Mink \rightarrow \suq$ with
$\xi(P_A) = \alpha J_A$, $\xi(P_0) = \beta W$, for some numbers
$\alpha$, $\beta$. More precisely, $\xi$ is a homomorphism of
algebras as long as $\alpha / \beta = - \lambda$.

We cannot invert $\xi$, though, since there is no relation like 
\begin{equation}
\label{eq:jgleicheins}
  W^2 - \lambda^2 J_A J_B g^{AB} = 1
\end{equation} 
in $\Mink$. However, for the case of constant positive mass,
$P_\mu P^\mu = m^2$, we find
\begin{equation}
  \xi(P_\mu P^\mu) = \beta^2 W^2-\alpha^2 J_A J_B g^{AB} = m^2 \,.
\end{equation}
We conclude that the image of the constant mass relation in $\Mink$
holds in $\suq$ if $\alpha = - m \lambda$ and $\beta = m $. This is
consistent with the requirement $\alpha / \beta = -\lambda$. We
conclude that $\Mink /\langle P_\mu P^\mu = m^2 \rangle$ is isomorphic
to the vectorial form of $\suq$. Setting aside the lack of $K^{-1}$ in
the vectorial $\suq$ we thus have an isomorphism
\begin{equation}
  \Euk / \langle P_\mu P^\mu = m^2 \rangle
  \stackrel{\phi}{\longrightarrow} \suq \rtimes \suq \,,
\end{equation}
where the action of the semidirect product on the right hand side is
the left Hopf adjoint action of $\suq$ on itself. The isomorphism is
given by $\xi\rtimes 1$ on the momenta and $1 \rtimes\id$ on the
rotations,
\begin{subequations}
\begin{xalignat}{2}
  \phi(P_A) &= -m \lambda J_A \rtimes 1 \,,&
  \phi(P_0) &= m  W \rtimes 1 \\
  \phi(J_A) &= 1 \rtimes J_A \,,& 
  \phi(W) &= 1 \rtimes  W \,.
\end{xalignat}
\end{subequations}
Introducing
\begin{equation}
\label{eq:Jnull}
  J_0 := -\lambda^{-1} W \,,
\end{equation}
we can write $\phi$ more compactly as
\begin{xalignat}{2}
\label{eq:phicompact}
  \phi(P_\mu) &= -m \lambda J_\mu \rtimes 1 \,,&
  \phi(J_\mu) &= 1 \rtimes J_\mu \,, 
\end{xalignat}
where $\mu$ runs through $\{0,-,+,3\}$. Note, however, that $J_\mu$ is
no 4-vector operator. The introduction of $J_0$ merely allows for a
more compact notation. For example, Eq.~\eqref{eq:jgleicheins} can be
written as $\lambda^2 J_\mu J^\mu = 1$. Furthermore, it is convenient
to give the pre-image of $K$ a name
\begin{equation}
  \pi := m\phi^{-1}(K) = (P_0 - P_3) \,.
\end{equation}

For the semidirect product of a Hopf algebra $H$ with itself by the
left Hopf adjoint action we have the following isomorphism of algebras
\begin{equation}
\label{eq:Euklid1}
  \psi:\, H\rtimes H \longrightarrow H \otimes H \,,\qquad
  \psi(g\rtimes h) = gh_{(1)} \otimes h_{(2)} \,.
\end{equation}
First, we prove that $\psi$ is a homomorphism
\begin{align}
  \psi[(g\rtimes h)(g'\rtimes h')]
  &= \psi[g(h_{(1)}\tr g') \rtimes h_{(2)}h' ]
  = g h_{(1)}g'S(h_{(2)})h_{(3)} h'_{(1)} \otimes h_{(4)}h'_{(2)}\notag\\
  &= g h_{(1)}g'h'_{(1)} \otimes h_{(2)}h'_{(2)}
  = (g h_{(1)} \otimes h_{(2)})(g'h'_{(1)} \otimes h'_{(2)})\notag\\
  &=\psi(g\rtimes h)\,\psi(g'\rtimes h') \,.
\end{align}
The invertibility can be shown directly, by defining
\begin{equation}
  \psi^{-1}(g\otimes h) := gS(h_{(1)}) \rtimes h_{(2)} \,, 
\end{equation}
and checking that
\begin{equation} 
\begin{aligned}
  (\psi\circ \psi^{-1})(g\otimes h)
  &= \psi[gS(h_{(1)}) \rtimes h_{(2)}] =
  gS(h_{(1)}) h_{(2)} \otimes h_{(3)} = g \otimes h \\
  (\psi^{-1}\circ\psi)(g\rtimes h)
  &= \psi[gh_{(1)} \otimes h_{(2)}] =
  gh_{(1)} S(h_{(2)}) \rtimes h_{(3)} = g \rtimes h \,.
\end{aligned}
\end{equation}
Thus, Eq.~\eqref{eq:Euklid1} tells us, that we have the sequence of
Isomorphisms
\begin{equation}
\label{eq:EukIso}
  \Euk / \lrAngle{ P_\mu P^\mu = m^2 }
  \stackrel{\phi}{\longrightarrow}\suq \rtimes \suq
  \stackrel{\psi}{\longrightarrow}\suq \otimes \suq\,.
\end{equation}
Through these isomorphisms we get a full understanding of the structure
of the $q$-Euclidean algebra.

One particularly interesting fact is that there is a whole $\suq$
subalgebra of $\Euk$ which commutes with the momenta $\Mink$. This
subalgebra is embedded by the map
\begin{equation}
\label{eq:idef}
  i:\,\suq \longrightarrow \Euk \,,\qquad
  i= \phi^{-1} \circ \psi^{-1} \circ (1\otimes \id) \,,
\end{equation}
which computes to 
\begin{equation}
\begin{aligned}
  i(J_\pm) &= J_\pm + \lambda^{-1} P_\pm \pi^{-1} K \\
  i(J_3) &= m\lambda^{-1} \pi^{-1}K - m^{-1}(\lambda^{-1} P_0 W
            + g^{AB} P_A J_B)\\
  i(W) &=  m^{-1}( P_0 W + \lambda g^{AB} P_A J_B ) \\
  i(K) &=  m\, \pi^{-1}K \,.
\end{aligned}
\end{equation}
Observe that the images of $J_A$ do not exist in $\Euk$ proper, since
they all involve the inverse of $\pi = P_0- P_3$, which is not an
element of $\Mink$.

\subsection{The Center of the $q$-Euclidean Algebra}

We wonder where precisely the condition $P_\mu P^\mu = m^2$ has
entered into our considerations. Which of the results do still hold if
the mass shell condition is relaxed?

Towards this end we list the commutation relations between rotations
and translations
\begin{xalignat}{3}
  [J_-, P_-] &= 0 &
  [J_-, P_3] &=  q^{-1}P_- K &
  [J_-, P_+] &=  P_3  K \notag\\
  [J_+, P_-] &= -P_3 K &
  [J_+, P_3] &= -q P_+ K &
  [J_+, P_+] &= 0\\
  [K ,P_-] &= -q^{-1}\lambda P_- K &
  [K ,P_3] &= 0 &
  [K ,P_+] &= q\lambda P_+ K  \notag
\end{xalignat}
Let us check what relations still hold within $i(\suq)$. We compute
for example
\begin{equation}
\begin{split}
  i(K)i(J_+)
  &= m\pi^{-1}K (J_+ + \lambda^{-1} P_+ \pi^{-1}K)\\
  &= (q^2 J_+ + q P_+ \pi^{-1}K) m\pi^{-1}K +
     m\lambda^{-1} P_+ (\pi^{-1}K)^2 \\
  &= (q^2 J_+ + (q + \lambda^{-1}) P_+ \pi^{-1}K ) m\pi^{-1}K\\
  &= q^2 i(J_+)i(K) \,,
\end{split}
\end{equation}
telling us that the relation $KJ_+ = q^2 J_+ K$ is preserved under $i$.
Similarly, we find that the image of $KJ_- = q^{-2}J_- K$ still holds
in $i(\suq)$. Hence, we did not use the mass shell condition for these
two relations. However, for the relation $\lambda [2]( q J_+ J_- -
q^{-1} J_- J_+) = 1 - K^2$ we find 
\begin{equation}
  i[\lambda [2]( q J_+ J_- - q^{-1} J_- J_+)] =
  1 - \frac{P_\mu P^\mu}{m^2}\, i(K)^2 \,,
\end{equation}
such that this relation holds in $\Euk$ precisely if the mass shell
condition $P_\mu P^\mu = m^2$ holds. We conclude that without the mass
shell condition $i$ is no longer a homomorphism of algebras.

Now, we check if $i(\suq)$ still commutes with all translations.
Setting aside the problem that $\pi^{-1}$ does not exist in $\Euk$
proper, we compute for example
\begin{equation}
\begin{split}
  P_+ i(J_-)
  &= P_+ (J_- + \lambda^{-1} P_- \pi^{-1}K )
  = J_- P_+ + (\lambda^{-1}P_+ P_- - P_3 \pi) \pi^{-1}K\\
  &= J_- P_+ + \lambda^{-1}(P_- P_+) \pi^{-1}K 
  = i(J_-) P_+ \,.
\end{split}
\end{equation}
In the same manner we find, that all of $i(\suq)$ commutes with all
translations. This holds in particular for $i(W)$ which is furthermore
a scalar with respect to rotations, since it is made up of the scalars
$P_0$, $W$, and $\vec{P}\cdot \vec{J}$. In conclusion we
have\footnote{Using the Casimir operators of $\Euk$, the orbital
  angular momentum relation of \cite{Cerchiai:1998} can be
  equivalently written as
\begin{equation}
\label{eq:Orbital}
  \lambda(\vec{P}\cdot\vec{J}) = P_0(1-W)
  \qquad\Leftrightarrow\qquad \Heli =P_0 \,.
\end{equation}}
\begin{Proposition}
  The center of the $q$-Euclidean algebra $\Euk$ is generated by
  $P_\mu P^\mu$, $P_0$ and
  \begin{equation}
  \label{eq:Helidef}
     \Heli := m\,i(W) = P_0 W + \lambda g^{AB} P_A J_B
     = -\lambda P_\mu J^\mu \,.
  \end{equation}
\end{Proposition}

\subsection{The Pauli-Lubanski Vector in the $q$-Deformed Setting}

In the undeformed case one considers the Pauli-Lubanski (pseudo)
vector
\begin{equation}
\label{eq:PauliLubanski1}
  W_\mu^{q=1} := -\frac{1}{2}\, \varepsilon_{\mu\nu\sigma\tau}
  V^{\nu\sigma} P^\tau \,,
\end{equation}
where $V^{\nu\sigma}$ is the matrix of Lorentz generators. Its
usefulness is due to the following two properties:
\begin{itemize}
\item[(i)] $W_\mu$ is a $4$-vector operator of the Poincar{\'e} algebra.
\item[(ii)] Each component $W_\mu$ commutes with all translations
  $P^\tau$.
\end{itemize}
If we demand further that $W_\mu$ be linear in the Lorentz generators
and the translations, conditions (i) and (ii) determine the
Pauli-Lubanski vector up to a constant factor. From (i) and (ii) we
deduce that $W_\mu W^\mu$ is a Casimir operator. Physically, this
Casimir operator turns out to correspond to spin.

In the $q$-deformed case we are tempted to define $W_\mu$ analogously
by Eq.~\eqref{eq:PauliLubanski1} with the $q$-deformed versions of the
epsilon tensor, the matrix of Lorentz generators, and the
translations. By construction, this would be a 4-vector operator.
However, it turns out that with this naive approach property (ii) will
not hold.  Therefore, we will try to find a way to construct $W_\mu$
such that (ii) holds, as well.

Let us start with the zero component $W_{0}$. It has to commute with
all translations to satisfy (ii) and with all rotations since the zero
component of a 4-vector is a scalar with respect to rotations. Thus,
it has to commute with all of the $q$-Euclidean algebra $\Euk$. If we
assume that as in the undeformed case $W_0$ is itself a member of
$\Euk$, we conclude that $W_0$ has to be an element of the center of
the $q$-Euclidean algebra, which we computed in the preceding section.
Since the momenta carry dimensions $W_0$ has to be linear in the
momenta. Hence $W_0$ must be a linear combination of $P_0$ and
$\Heli$. The additional requirement that $W_0$ has to have the right
undeformed limit determines
\begin{equation}
\label{eq:PauliLubnull}
  W_0 := \lambda^{-1}(\Heli - P_0) =
  \lambda^{-1} (W-1)P_0 + g^{AB} J_A P_B 
\end{equation}
up to an overall factor that tends to one as $q\rightarrow 1$.

Now that we have a good candidate for the zero component of the
$q$-Pauli-Lubanski vector we have to see if it can be boosted to a
4-vector. First we have to ask what type of vector operator we would
expect it to be. Recall from Sec.~\ref{sec:TensorHopf} that we have to
distinguish between left and right tensor operators. A short
calculation shows that for any translation $p\in \Mink$ and any
Lorentz transformation $h\in\slC$
\begin{align}
  (W_0 \tl\adR h)\, p
  &= S(h_{(1)})W_0 h_{(2)} p = S(h_{(1)})W_0 (h_{(2)}\tr p) h_{(3)} \notag\\
  &= S(h_{(1)}) (h_{(2)}\tr p) W_0 h_{(3)} 
  = (S(h_{(1)})_{(1)} h_{(2)}\tr p) S(h_{(1)})_{(2)} W_0 h_{(3)} \notag\\
  &= (S(h_{(2)}) h_{(3)}\tr p) S(h_{(1)}) W_0 h_{(4)} \notag\\
\label{eq:RightBoostcommute}
  &= p \,(W_0 \tl\adR h)
\end{align}
Hence, a right boosted $W_0$ commutes with all translations. This is
not be the case for $\adL\tr W_0$. Hence, the $q$-Pauli-Lubanski vector
will satisfy property (ii) only if it is a right vector operator
$W_{\tilde{\mu}}$.

\subsection{Boosting the $q$-Pauli-Lubanski Vector}

If $W_0 = W_{\tilde{0}}$ as defined in \eqref{eq:PauliLubnull} really
is a left 4-vector operator, which is not necessarily so, then the
other components are given, uniquely, by Eqs.~\eqref{eq:Boost1b}. We
will now determine $W_{\tilde{\mu}}$ and rigorously show that it is a
right 4-vector operator. Up to a constant factors $W_{\tilde{0}}$ is
the sum of two parts, $\Heli$ and $P_0$, which we will treat
separately.

\subsubsection{Boosting $\Heli$}
The explicit calculations of the right adjoint action of the boosts on
$\Heli$ by Eqs.~\eqref{eq:Boost1b} turn out to be very lengthy. It is
more efficient to start with a more abstract consideration.

We observe that for all boosts $h\in\SUq^\op$ we have
\begin{equation}
\label{eq:Heliboost1}
  \lrAngle{(J_\mu)_{(1)}, h}(J_\mu)_{(2)}
  = J_{\mu'}\Lambda(h)^{\mu'}{}_\mu \,,
\end{equation}
where $\lrAngle{\,\cdot\,,\,\cdot\,}$ is the dual pairing of $\suq$
and $\SUq$. We exemplify this result for $J_+$,
\begin{equation}
\begin{split}
  \lrAngle{(J_+)_{(1)}, B^a{}_b} (J_+)_{(2)} &=
  \lrAngle{J_+ , B^a{}_b} K + \lrAngle{1, B^a{}_b} J_+ \\
  &= \lambda[2]^{-1}(\sigma_+)_{ab}\,(J_3 - J_0 ) +\delta^a_b\,J_+ \\
  &=  \begin{pmatrix} J_+ & 0 \\
      -q^{\frac{1}{2}}\lambda [2]^{-1/2}(J_3 - J_0) & \quad J_+
    \end{pmatrix}\\
  &=J_{\mu'}\Lambda(h)^{\mu'}{}_+ \,.
\end{split}
\end{equation}
Applying the map $\phi^{-1}$ as defined in Eq.~\eqref{eq:phicompact}
to Eq.~\eqref{eq:Heliboost1} we get 
\begin{equation}
\label{eq:PauliLubanski2}
  \adL h \tr \phi^{-1}(l) = \lrAngle{l_{(1)}, h} \phi^{-1}(l_{(2)}) \,.
\end{equation}
for all $l \in \suq$ and $h\in\SUq^\op$. For example, for $l=
\phi(P_+) = -m\lambda J_+$ the left adjoint action of the boost
generators on $P_+$ can be written as
\begin{equation}
  \adL B^a{}_b \tr P_+
  = \delta^a_b\,P_+ - \lambda \lrAngle{J_+, B^a{}_b} \pi 
  = \delta^a_b\,P_+ + \lambda[2]^{-1}\,
    (\sigma_+)_{ab}\,(P_3- P_0) \,,
\end{equation}
which is the same as in Eq.~\eqref{eq:AppBoostAct}. 

Let $i$ be the map that has been defined in Eq.~\eqref{eq:idef}.  We
try to commute $i(l)$, $l\in\suq$, with a boost $h\in\SUq^\op$ using
Eq.~\eqref{eq:PauliLubanski2}:
\begin{equation}
\begin{split}
  i(l)\, h
  &= \phi^{-1}[S(l_{(1)})] l_{(2)} h
  = \phi^{-1}[S(l_{(1)})] h_{(2)} l_{(3)}
     \lrAngle{ S(l_{(2)}), h_{(1)} } \lrAngle{ l_{(4)} , h_{(3)} }\\
  &= h_{(3)} \{\adL S^{-1}(h_{(2)})\tr \phi^{-1}[S(l_{(1)})] \} l_{(3)}
  \lrAngle{ S(l_{(2)}), h_{(1)} } \lrAngle{ l_{(4)} , h_{(4)} }\\
  &= h_{(3)} \lrAngle{S(l_{(2)}),S^{-1}(h_{(2)})}
  \phi^{-1}[S(l_{(1)})] l_{(4)}
     \lrAngle{ S(l_{(3)}), h_{(1)} } \lrAngle{ l_{(5)} , h_{(4)} }\\
  &= h_{(1)} \phi^{-1}[S(l_{(1)})] l_{(2)} \lrAngle{l_{(3)},h_{(2)}}
\end{split}
\end{equation}
This leads to a remarkably simple formula for the right adjoint action
of a boost on $i(l)$
\begin{equation}
  i(l)\tl\adR h = i( l_{(1)} \lrAngle{l_{(2)},h }) \,.
\end{equation}
For $l = S(J_\mu)$ this formula becomes
\begin{equation}
\begin{split}
  i[S(J_\mu)]\tl\adR h
  &= i[ S(J_\mu)_{(1)} \lrAngle{ S(J_\mu)_{(2)},h }]
  = i[ S( (J_\mu)_{(2)} \lrAngle{ (J_\mu)_{(1)},S^{-1}h } ) ]\\
  &= i[ S(J_{\mu'}) \Lambda(S^{-1}h)^{\mu'}{}_\mu ] \,,
\end{split}
\end{equation}
which tells us that $i(S(J_{\mu}))$ transforms under boosts as a right
lower 4-vector operator.

It remains to check whether $i(S(J_{\mu}))$ transforms as right
4-vector under rotations. We observe that $\phi^{-1}$ maps the
3-vector $J_\mu$ to the 3-vector $P_\mu$ and the $\suq$-scalar $J_0$
to the scalar $P_0$. Hence, for $a,b\in\suq$ we have
\begin{equation}
  \adL b \tr \phi^{-1}(a) = \phi^{-1}(\adL b \tr a) \,.
\end{equation}
Now we are prepared to tackle the right action of a rotation on $i(Sa)$
\begin{equation}
\begin{split}
  i(Sa) \tl\adR b
  &= S(b_{(1)}) \phi^{-1}[S((Sa)_{(1)})] (Sa)_{(2)} b_{(2)} \\
  &= \{ S(b_{(1)})_{(1)} \tr \phi^{-1}[S((Sa)_{(1)})] \}
     \{S(b_{(1)})_{(2)} (Sa)_{(2)} b_{(2)}\} \\
  &= \phi^{-1}[(S(b_{(2)})_{(1)} S((Sa)_{(1)}) S(S(b_{(2)})_{(2)})]
     S(b_{(1)}) (Sa)_{(2)} b_{(3)}\\
  &= \phi^{-1}[S( S(b_{(2)}) (Sa)_{(1)} b_{(3)})]
     S(b_{(1)}) (Sa)_{(2)} b_{(4)}\\
  &= \phi^{-1}[S( S(b_{(1)})_{(1)} (Sa)_{(1)} b_{(2)(1)})]
     S(b_{(1)})_{(2)} (Sa)_{(2)} b_{(2)(2)}\\
  &= i(Sa \tl\adR b) = i(S(\adL S^{-1}b \tr a)) \,.
\end{split}
\end{equation}
This shows that since $S(J_\mu)$ transforms as a right lower 4-vector
under rotations, so does $i(S(J_\mu))$. In conclusion we have
\begin{Proposition}
  The set of operators
  \begin{equation}
    \Heli_{\tilde{\mu}} := -m\lambda\, i(S(J_\mu))
  \end{equation}
  is a right lower 4-vector operator of the $q$-Lorentz algebra. Since
  furthermore $\Heli_{\tilde{0}} = \Heli$, $\Heli_{\tilde{\mu}}$ is
  the unique right lower 4-vector operator with zero component
  $\Heli$.
\end{Proposition}
All that remains to do is to compute $\Heli_{\tilde{\mu}}$
explicitly:
\begin{equation}
\begin{aligned}
  \Heli_{\tilde{0}} &= m\,i(W) = P_0 W + \lambda g^{AB} P_A J_B\\
  \Heli_{\tilde{\pm}} &= P_\pm + \lambda J_\pm K^{-1}\pi \\
  \Heli_{\tilde{3}} &= m\,i(W-K^{-1}) 
  = P_0 W + \lambda  g^{AB} P_A J_B - K^{-1}\pi \,. 
\end{aligned}
\end{equation}
Observe that these expressions do not contain $\pi^{-1}$, hence, they
are proper members of $\Euk$, that is,
\begin{equation}
\begin{aligned}
  \Heli_{\tilde{0}} &=  W P_0  - q \lambda J_+P_-
     - q^{-1}\lambda J_- P_+ + \lambda J_3 P_3 \\
  \Heli_{\tilde{\pm}} &= P_\pm + \lambda J_\pm K^{-1}(P_0 - P_3) \\
  \Heli_{\tilde{3}} &=  (W - K^{-1})  P_0  - q \lambda J_+P_-
     - q^{-1}\lambda J_- P_+ + (\lambda J_3 + K^{-1}) P_3 \,.
\end{aligned}
\end{equation}
Finally, we recall that the square of $\Heli_{\tilde{\mu}}$ must be a
Casimir operator. After lengthy calculations we find 
\begin{equation}
  \Heli^{\tilde{\mu}}\Heli_{\tilde{\mu}} = P_\mu P^\mu \,.
\end{equation}
We conclude that squaring $\Heli_{\tilde{\mu}}$ alone does not yield a
new Casimir operator.

\subsubsection{Boosting $P_0$}
The next step in our calculation of the $q$-Pauli-Lubanski vector is
to find a right 4-vector operator with $P_0$ as zero component. With
a universal $\R$-matrix of the $q$-Lorentz algebra, we can generically
turn a left 4-vector operator into a right 4-vector operator. Defining
for the left 4-vector operator $P_\mu$
\begin{equation}
  j(P_\mu) := S^2(\R_{[1]})(\R_{[2]}\tr P_\mu)
    = S^2\bigl(\R_{[1]}\Lambda(\R_{[2]})^{\mu'}{}_\mu\bigr) P_{\mu'}
\end{equation}
we check that for any $q$-Lorentz transformation $h\in\slC$ we have
\begin{equation}
\begin{split}
  j(P_\mu)\tl\adR h
  &= S(h_{(1)}) S^2(\R_{[1]})(\R_{[2]}\tr P_\mu)h_{(2)}\\
  &= S(h_{(1)})S^2(\R_{[1]}) h_{(3)}
    \bigl( S^{-1}(h_{(2)})\R_{[2]} \tr P_\mu \bigr) \\
  &= S\bigl(S(\R_{[1]})h_{(1)}\bigr)h_{(3)}
    \bigl(S^{-1}(S(\R_{[2]})h_{(2)})\tr P_\mu \bigr) \\
  &= S\bigl(h_{(2)}S(\R_{[1]})\bigr)h_{(3)}
    \bigl(S^{-1}(h_{(1)}S(\R_{[2]}))\tr P_\mu \bigr)\\
  &= S^2(\R_{[1]}) S(h_{(2)})h_{(3)}
    \bigl(\R_{[2]}S^{-1}(h_{(1)})\tr P_\mu \bigr)\\
  &= S^2(\R_{[1]})(\R_{[2]}S^{-1}h\tr P_\mu)\\
  &= j(P_{\mu'}) \Lambda( S^{-1}h )^{\mu'}{}_\mu \,,
\end{split}
\end{equation}
thus, $j(P_\mu)$ is indeed a right 4-vector operator. Recall from
Sec.~\ref{sec:L-Matrices}, that the object
\begin{equation}
  (L^\Lambda_+)^\mu{}_\nu := \R_{[1]}\Lambda(\R_{[2]})^{\mu}{}_\nu
\end{equation}
that appears in the definition of $j(P_\mu)$ is an $L$-matrix.
Furthermore, we recall from Eq.~\eqref{eq:LorentzR} that there are two
universal $\R$-matrices of the $q$-Lorentz algebra, which are composed
of the $\R$-matrix of $\slq$ according to
\begin{xalignat}{2}
  \RI &= \R^{-1}_{41}\R^{-1}_{31}\R_{24}\R_{23}\,, &
  \RII &= \R^{-1}_{41}\R_{13}\R_{24}\R_{23} \,.
\end{xalignat}
We will now compute the $L$-Matrix for $\RI$. We have for the
$(\frac{1}{2},\frac{1}{2})$-form of the vector representation
\begin{equation}
\begin{split}
 \bigl( L^{(\frac{1}{2},\frac{1}{2})}_{\mathrm{I}+} \bigr)^{ab}{}_{cd}
 &= \bigl(\id\otimes\id\otimes
    \rho^{\frac{1}{2}}\otimes\rho^{\frac{1}{2}} \bigr)
   (\RI)^{ab}{}_{cd} \\
 &= \bigl( L^{\frac{1}{2}}_- \bigr)^{b}{}_{b'}
    \bigl( L^{\frac{1}{2}}_- \bigr)^{a}{}_{a'} \otimes
    \bigl( L^{\frac{1}{2}}_+ \bigr)^{b'}{}_{d}
    \bigl( L^{\frac{1}{2}}_+ \bigr)^{a'}{}_{c}\\
 &= B^{b}{}_{d}B^{a}{}_{c} \,,
\end{split}
\end{equation}
where $B^a{}_b \in \SUq^\op$ is the matrix of boosts. For the 4-vector
form of this $L$-matrix we then find
\begin{equation}
\label{eq:LLambda}
  (L^\Lambda_{\mathrm{I}+})^\mu{}_\nu =
  \begin{pmatrix}
    1 & 0 & 0 & 0 \\
    0 & a^2 & b^2 & q^{\frac{1}{2}}[2]^{\frac{1}{2}} ab \\
    0 & c^2 & d^2 & q^{\frac{1}{2}}[2]^{\frac{1}{2}} cd \\
    0 & q^{\frac{1}{2}}[2]^{\frac{1}{2}} ac &
    q^{\frac{1}{2}}[2]^{\frac{1}{2}} bd & (1 + [2]bc) 
  \end{pmatrix}
\end{equation}
with respect to the basis $\{0,-,+,3\}$. This matrix of generators
becomes more familiar if we write it in block diagonal form
\begin{equation}
\label{eq:tcorrep}
  (L^\Lambda_{\mathrm{I}+})^\mu{}_\nu = 
  \begin{pmatrix}
    1 & 0 \\
    0 & t^A{}_B
  \end{pmatrix} \,,
\end{equation}
so we can see that $t^A{}_B$, $A,B\in \{-1,0,1\}$ is the 3-dimensional
corepresentation matrix of $\SUq^\op$
\cite{Koornwinder:1989,Masuda:1991}. From the block diagonal form we
deduce that
\begin{equation}
  j(P_0) = P_0 \,,
\end{equation}
so we get
\begin{Proposition}
  The set of operators
  \begin{equation}
  \label{eq:Pboosted}
    j(P_\mu) :=
    S^2\bigl[(L^\Lambda_{\mathrm{I}+})^{\mu'}{}_\mu\bigr] P_{\mu'}
  \end{equation}
  is a right lower 4-vector operator
  of the $q$-Lorentz algebra. Since furthermore $j(P_0) = P_0$,
  $j(P_\mu)$ is the unique right lower 4-vector operator with zero
  component $P_0$.
\end{Proposition}
With Eq.~\eqref{eq:Pboosted} we find the explicit expressions
\begin{equation}
\begin{aligned}
  j(P_0) &= P_0  \\
  j(P_-) &= a^2\, P_- + q^{-4}c^2\, P_+ +
            q^{-\frac{3}{2}}[2]^{\frac{1}{2}}ac\, P_3 \\
  j(P_+) &= q^{4}b^2\, P_- + d^2\, P_+ +
            q^{\frac{5}{2}}[2]^{\frac{1}{2}}bd\, P_3 \\
  j(P_3) &= q^{\frac{5}{2}}[2]^{\frac{1}{2}}ab\, P_- +
            q^{-\frac{3}{2}}[2]^{\frac{1}{2}}cd\, P_+ +
            (1+[2]bc)\, P_3 \,.
\end{aligned}
\end{equation}

Finally, we want to calculate the square of $j(P_\mu)$ which must be a
Casimir operator. First, we note that since $P_0$ commutes with all
momenta and $j(P_\mu)$ is the right boosted $P_0$, the reasoning of
Eq.~\eqref{eq:RightBoostcommute} applies, that is, all momenta $P_\mu$
commute with $j(P_\nu)$,
\begin{equation}
 P_\mu \, j(P_\nu) = j(P_\nu) \, P_\mu \,.
\end{equation}
Moreover, we have
\begin{equation}
\begin{split}
  (L^\Lambda_{\mathrm{I}+})^{\mu}{}_{\nu}
  (L^\Lambda_{\mathrm{I}+})^{\sigma}{}_{\tau} \,\eta^{\tau\nu}
  &= R_{[1]}R_{[1']} \Lambda(R_{[2]})^\mu{}_\nu
     \Lambda(R_{[2']})^\sigma{}_\tau \,\eta^{\tau\nu}\\
  &= R_{[1]}R_{[1']} \Lambda(R_{[2]})^\mu{}_\nu
     \Lambda(R_{[2']})^{\sigma''}{}_\tau
     \,\eta^{\tau\nu}\eta_{\sigma'\sigma''}\eta^{\sigma\sigma'}\\
  &= R_{[1]}R_{[1']} \Lambda(R_{[2]})^\mu{}_\nu
     \Lambda(S^{-1}R_{[2']})^{\nu}{}_{\sigma'}
     \,\eta^{\sigma\sigma'}\\
  &= R_{[1]}R^{-1}_{[1']} \Lambda(R_{[2]}R^{-1}_{[2']})^\mu{}_{\sigma'}
     \,\eta^{\sigma\sigma'}\\
  &= \eta^{\sigma\mu} \,,
\end{split}
\end{equation}
where we have used Eq.~\eqref{eq:Tensor6}. With the last two equations
we can compute the square of $j(P_\mu)$ quite easily
\begin{equation}
\begin{split}
  j(P^\mu)j(P_\mu)
  &= S^2\bigl((L^\Lambda_{\mathrm{I}+})^{\mu'}{}_\mu\bigr)\,P_{\mu'} \,
     S^2\bigl((L^\Lambda_{\mathrm{I}+})^{\nu'}{}_\nu\bigr)\,P_{\nu'}
     \,\eta^{\nu\mu}\\
  &= \bigl[S^2\bigl((L^\Lambda_{\mathrm{I}+})^{\mu'}{}_\mu\bigr)
     S^2\bigl((L^\Lambda_{\mathrm{I}+})^{\nu'}{}_\nu\bigr)\eta^{\nu\mu}\bigr]
     \,P_{\nu'}P_{\mu'}\\
  &= \eta^{\nu\mu} \,P_{\nu}P_{\mu}\,.  
\end{split}
\end{equation}
Again, the square of one half of the $q$-Pauli-Lubanski vector alone
yields only the mass Casimir.

\subsubsection{The $q$-Pauli-Lubanski Vector}
We come to the following conclusion:
\begin{Proposition}
  The set of operators
  \begin{equation}
    W_{\tilde{\mu}} := \lambda^{-1}[\Heli_{\tilde{\mu}}-j(P_\mu)]
      = -m\,i(S(J_\mu))-\lambda^{-1}j(P_\mu)
  \end{equation}
  has the following properties:
  \begin{itemize}
    \item[(i)]  It is a right lower $4$-vector operator.
    \item[(ii)] Each component $W_\mu$ commutes with all translations
                $P_\tau$.
  \end{itemize}
  Furthermore, it is the unique right lower 4-vector operator with
  zero component $W_0 = \lambda^{-1}(\Heli - P_0)$. We will therefore
  call it the $q$-Pauli-Lubanski vector.
\end{Proposition}
Explicitly, the $q$-Pauli-Lubanski vector is
\begin{equation}
\label{eq:PauLubEnd}
\begin{aligned}
  W_{\tilde{0}} &= \lambda^{-1} (W -1)P_0  - q J_+P_-
     - q^{-1} J_- P_+ + J_3 P_3 \\
  W_{\tilde{-}} &= \lambda^{-1}[\lambda J_- K^{-1} P_0 + (1- a^2)P_-
         -q^{-4}c^2 P_+ - (\lambda J_- K^{-1}
         + q^{-\frac{3}{2}}[2]^{\frac{1}{2}}ac)P_3] \\
  W_{\tilde{+}} &= \lambda^{-1}[\lambda J_+ K^{-1} P_0 - q^{4}b^2 P_-
         + (1 - d^2) P_+ - (\lambda J_+ K^{-1}
         + q^{\frac{5}{2}}[2]^{\frac{1}{2}}bd)P_3] \\
  W_{\tilde{3}} &= \lambda^{-1} (W - K^{-1})  P_0
     - (qJ_+ + q^{\frac{5}{2}}\lambda^{-1}[2]^{\frac{1}{2}} ab) P_-
     - (q^{-1}J_- + q^{-\frac{3}{2}}\lambda^{-1}[2]^{\frac{1}{2}} cd) P_+
     \\ &\quad
     + ( J_3 + \lambda^{-1}K^{-1} - \lambda^{-1}(1+[2]bc)) P_3 \,.
\end{aligned}
\end{equation}

\section{The Little Algebras}
\label{sec:MainContrib2b}

\subsection{Little Algebras in the $q$-Deformed Setting}
\label{sec:Little1}

In classical relativistic mechanics the state of a free particle is
completely determined by its $4$-momentum. In quantum mechanics
particles can have an additional degree of freedom called spin. What
is spin?

Let us assume we have a free relativistic particle described by an
irreducible representation of the Poincar{\'e} algebra. We pick all states
with a given momentum,
\begin{equation}
  \mathcal{H}_p
  := \{ \Ket{\psi}\in \mathcal{H} :
  P_\mu\Ket{\psi} = p_\mu\Ket{\psi} \} \,,
\end{equation}
where $\mathcal{H}$ is the Hilbert space of the particle and
$p=(p_\mu)$ is the $4$-vector of momentum eigenvalues. If the state of
the particle is \emph{not} uniquely determined by the eigenvalues of
the momentum, then the eigenspace $\mathcal{H}_p$ will be degenerate.
In that case we need, besides the momentum eigenvalues, an additional
quantity to label the basis of our Hilbert space uniquely.  This
additional degree of freedom is spin. The spin symmetry is then the
set of Lorentz transformations that leaves the momentum eigenvalues
invariant and, hence, acts on the spin degrees of freedom only,
\begin{equation}
\label{eq:Little1}
  \mathcal{K}'_p := \{ h\in\mathcal{L} :
  P_\mu h \Ket{\psi} = p_\mu h \Ket{\psi}
  \text{ for all } \Ket{\psi}\in \mathcal{H}_p \} \,,
\end{equation}
where $\mathcal{L}$ is the enveloping Lorentz algebra. In mathematical
terms, $\mathcal{K}'_p$ is the stabilizer of $\mathcal{H}_p$. Clearly,
$\mathcal{K}'_p$ is an algebra, called the little algebra.

A priori, there are a lot of different little algebras for each
representation and each vector $p$ of momentum eigenvalues.  In the
undeformed case it turns out that for the physically relevant
representations (real mass) there are (up to isomorphism) only two
little algebras, depending on the mass being either positive or zero
\cite{Wigner:1939}. For positive mass we get the algebra of rotations,
$\mathcal{U}(\mathrm{su}_2)$, for zero mass an algebra that is
isomorphic to the algebra of rotations and translations of the
2-dimensional plane denoted by $\mathcal{U}(\mathrm{iso}_2)$.  The
proof that $\mathcal{K}'_p$ does not depend on the particular
representation but on the mass, does not generalize to the
$q$-deformed case: If we defined for representations of the
$q$-Poincar{\'e} algebra the little algebra as in Eq.~\eqref{eq:Little1},
it could well happen that $\mathcal{K}'_p$ for a spin-$\frac{1}{2}$
particle is not the same as for spin-1. We will therefore define the
$q$-little algebras differently.

In the undeformed case there is an alternative but equivalent
definition of the little algebras. $\mathcal{K}'_p$ is the algebra
generated by the components of the $q$-Pauli-Lubanski vector as defined in
Eq.~\eqref{eq:PauliLubanski1} with the momentum generators replaced by
their eigenvalues. Let us formalize this to see why this definition
works and how it is generalized to the $q$-deformed case.

Let $\mathcal{T}$ be the algebra of translations, $\mathcal{L}$ the
Lorentz algebra, both joined in a semidirect product to form the
Poincar{\'e} algebra $\mathcal{P} = \mathcal{T}\rtimes \mathcal{L}$. Let
$\chi_p$ be the map that maps the momentum generators to the
eigenvalues, $\chi_p(P_\mu) = p_\mu$. Being the restriction of a
representation, $\chi_p$ must extend to a one dimensional
$*$-representation $\chi_p: \mathcal{T} \rightarrow \mathbb{C}$, a
non-trivial condition only in the $q$-deformed case.  Noting that
every element of $\mathcal{P}$ can be written as a sum of products of
Lorentz transformations and translations, $\sum_i l_i t_i$, we extend
$\chi_p$ to a linear map $\tilde{\chi}_p: \mathcal{P} \rightarrow
\mathcal{L}$ by
\begin{equation}
  \tilde{\chi}_p(\sum l_i t_i) := \sum l_i \chi_p(t_i).
\end{equation}
The little algebra can now be alternatively defined as the unital
algebra generated by the images of the $q$-Pauli-Lubanski vector under
$\tilde{\chi}_p$,
\begin{equation}
\label{eq:Little2}
  \mathcal{K}_p
  := \mathbb{C}\lrAngle{\tilde{\chi}_p(W_\mu)}\,.
\end{equation}
Why is this a reasonable definition?  By construction the action of
every element of $\mathcal{P}$ on $\mathcal{H}_p$ is the same as of
its image under $\tilde{\chi}_p$. For any $\Ket{\psi}\in\mathcal{H}_p$
this means
\begin{equation}
  P_\mu \,\tilde{\chi}_p(W_\nu)\Ket{\psi}
  = \tilde{\chi}_p(P_\mu W_\nu)\Ket{\psi}
  = \tilde{\chi}_p(W_\nu P_\mu)\Ket{\psi}
  = p_\mu \,\tilde{\chi}_p(W_\nu)\Ket{\psi} \,,
\end{equation}
which shows that $\mathcal{K}_p \subset \mathcal{K}'_p$. It still
could happen, that $\mathcal{K}_p$ is strictly smaller than
$\mathcal{K}'_p$.  In the undeformed case there are theorems telling
us \cite{Blattner:1969,Dixmier} that this cannot happen, so we really have
$\mathcal{K}_p = \mathcal{K}'_p$. For the $q$-deformed case no such
theorem is known \cite{Schneider:1994}. However, if there were more
generators in the stabilizer of some momentum eigenspace they would
have to vanish for $q\rightarrow 1$. In this sense
Eq.~\eqref{eq:Little2} with the $q$-deformed Pauli-Lubanski vector can
be considered to define the $q$-deformed little algebras.

\subsection{Computation of the $q$-Little Algebras}
\label{sec:Little2}

To begin the explicit calculation of the $q$-deformed little algebras,
we need to figure out if there are eigenstates of $q$-momentum at all.
That is, we want to determine the one-dimensional $*$-representations
of $\Mink$, that is the homomorphisms of $*$-Algebras $\chi:\Mink
\mapsto \mathbb{C}$. Let us again denote the eigenvalues of the
generators by lower case letters $p_\mu := \chi(P_\mu)$. For $\chi$ to
be a $*$-map we must have $p_0$, $p_3$ real and $p_+^* = -q p_-$. To
find the conditions for $\chi$ to be a homomorphism of algebras, we
apply $\chi$ to the relations~\eqref{eq:XX-Rel1} of $\Mink$, yielding
\begin{equation}
  p_A(p_0-p_3) = 0 \,.
\end{equation}
There are two cases. The first is $p_0 \neq p_3$, which immediately
leads to $p_A = 0$, and $p_0 = \pm m$. The second case is $p_0 = p_3$,
leading to $m^2 = - \lvert p_- \rvert^2 - \lvert p_+ \rvert^2$. Hence,
if the mass $m$ is to be real, we must have $p_\pm = 0$.

To summarize, for real mass $m$ we have a massive and a massless type
of momentum eigenstates with eigenvalues given by
\begin{xalignat}{2}
  (p_0,p_-,p_+,p_3) =
  \begin{cases}
  (\pm m ,0,0,0) & m > 0 \\
  (k ,0,0,k) & m = 0,\, k\in\mathbb{R}
  \end{cases}
\end{xalignat}
Now, we need to move the momentum generators in the expressions of the
$q$-Pauli-Lubanski vector to the right and replace them with these
eigenvalues.

\subsubsection{The Massive Case}
In Eqs.~\eqref{eq:PauLubEnd} the momenta have already been moved to
the right, so we can simply replace them with $(P_0,P_-,P_+,P_3)
\rightarrow (m ,0,0,0)$. We get
\begin{equation}
\begin{aligned}
  \tilde{\chi}_p(W_{\tilde{0}}) &= \lambda^{-1} (W -1)m \\
  \tilde{\chi}_p(W_{\tilde{-}}) &= J_- K^{-1} m \\
  \tilde{\chi}_p(W_{\tilde{+}}) &= J_+ K^{-1} m \\
  \tilde{\chi}_p(W_{\tilde{3}}) &= \lambda^{-1} (W - K^{-1}) m \,,
\end{aligned}
\end{equation}
so the set of generators of the little algebra is essentially
$\{W,K^{-1},J_\pm K^{-1}\}$. Since $K^{-1}$ stabilizes the momentum
eigenspace, so does its inverse $K$. Hence, it is safe to add $K$ to
the little algebra which would exist, anyway, as operator within a
representation. We thus get
\begin{equation}
  \mathcal{K}_{m} := \mathcal{K}_{(m ,0,0,0)} = \suq \,,
\end{equation}
completely analogous to the undeformed case.

\subsubsection{The Massless Case}
The massless case is more complicated. Replacing in
Eqs.~\eqref{eq:PauLubEnd} the momentum generators with
$(P_0,P_-,P_+,P_3) \rightarrow (k ,0,0,k)$ we get
\begin{equation}
\begin{aligned}
  \tilde{\chi}_p(W_{\tilde{0}}) &= \lambda^{-1} (K -1) k \\
  \tilde{\chi}_p(W_{\tilde{-}}) &=
    -\lambda^{-1} q^{-\frac{3}{2}}[2]^{\frac{1}{2}} ac \, k\\
  \tilde{\chi}_p(W_{\tilde{+}}) &=
    -\lambda^{-1}  q^{\frac{5}{2}}[2]^{\frac{1}{2}} bd \, k\\
  \tilde{\chi}_p(W_{\tilde{3}}) &=
  \lambda^{-1}\bigl( K - (1 +[2] bc) \bigr) k \,.
\end{aligned}
\end{equation}
The set of generators of the little algebra is essentially $\{K, ac,
bd, bc \}$. The commutation relations of these generators can be
written more conveniently in terms of $K$ and
\begin{xalignat}{3}
\label{eq:Little3}
  N_- &:= q^{\frac{1}{2}}[2]^{\frac{1}{2}} ac \,,&
  N_+ &:= q^{\frac{1}{2}}[2]^{\frac{1}{2}} bd \,,&
  N_3 &:= 1 +[2] bc \,,
\end{xalignat}
or equivalently $N_A = t^3{}_A$, for $t^A{}_B$ as defined by
Eqs.~\eqref{eq:LLambda} and \eqref{eq:tcorrep}. The commutation
relations are
\begin{xalignat}{3}
  N_B N_A\, \varepsilon^{AB}{}_C &= -\lambda N_C \,,&
  N_A N_B\, g^{BA} &= 1 \,,&
  K N_A &= q^{-2A} N_A K \,,
\end{xalignat}
and the conjugation properties 
\begin{xalignat}{2}
  N_A^* &= N_B\, g^{BA} \,, & K^* &= K \,.
\end{xalignat}
In words: The $N_A$ generate the opposite algebra of a unit quantum
sphere, $\qsphere$ \cite{Podles:1987}. $K$, the generator of $\uq$, acts
on $N_A$ as on a right $3$-vector operator. In total we have
\begin{equation}
  \mathcal{K}_{0} := \mathcal{K}_{(k ,0,0,k)} = \uq \ltimes \qsphere \,.
\end{equation}
As opposed to the massive case, $\mathcal{K}_{0}$ is no Hopf algebra.
However, since $N_A = t^3{}_A$ and $\Delta(t^A{}_C) = t^A{}_B\otimes
t^B{}_C$, we have
\begin{equation}
 \Delta(N_{B}) = N_A \otimes t^A{}_B \,,
\end{equation}
hence, $\mathcal{K}_{0}$ is a right coideal.

The only irreducible $*$-representations of $\mathcal{K}_{0}$ are
one-dimensional. They depend on a real parameter $\alpha$ and are
defined on the single basis vector $\Ket{\alpha}$ by
\begin{xalignat}{3}
\label{eq:Nrep}
  K\Ket{\alpha} &= \alpha\Ket{\alpha} \,,&
  N_\pm\Ket{\alpha} &= 0 \,, & N_3 \Ket{\alpha} &= \Ket{\alpha}\,.
\end{xalignat}
Unlike for the undeformed case, no infinite-dimensional irreducible
representation exists.

\chapter{Massive Spin Representations}
\label{sec:SpinReps}

\section{Representations in an Angular Momentum Basis}

\subsection{The Complete Set of Commuting Observables}
\label{sec:SpinReps1}

We want to construct a massive irreducible representations in a basis
that can be given a physical interpretation. Massive irreducible means
that within the representation we have
\begin{equation}
  P_\mu P^\mu = m^2
\end{equation}
for some real positive constant $m$, $P_\mu P^\mu$ being the mass
Casimir operator. We have shown in Sec.~\ref{sec:Little2} that there
are rest states, that is, momentum eigenstates,
$P_\mu\Ket{\psi_0}=p_\mu\Ket{\psi_0}$, with $(p_\mu) = (p_0,p_-,p_+,p_3)
= (m ,0,0,0)$. On these rest states the $q$-Pauli-Lubanski vector acts
as
\begin{xalignat}{2}
  W_{\tilde{0}}\Ket{\psi_0} &= m\lambda^{-1}(W-1)\Ket{\psi_0} \,,&
  W\!_{\tilde{A}}\Ket{\psi_0} &= -m S(J_A) \Ket{\psi_0} \,,
\end{xalignat}
from which it follows that
\begin{equation}
  W^{\tilde{\mu}} W_{\tilde{\mu}} \Ket{\psi_0}
  = 2m^2\lambda^{-2} (1-W)\Ket{\psi_0} \,.
\end{equation}
The spin Casimir $W^{\tilde{\mu}} W_{\tilde{\mu}}$ must be constant,
thus, the angular momentum must be constant within the rest frame.
According to Eq.~\eqref{eq:AppWrep} the possible values are
\begin{equation}
\label{eq:W-fix}
  W\Ket{\psi_0} = [2]^{-1}\bigl(q^{(2s+1)}+q^{-(2s+1)}\bigr)\Ket{\psi_0}\,,
\end{equation}
where $s\in\frac{1}{2}\mathbb{N}_0$ is a half integer.  For the spin
Casimir this means
\begin{equation}
  W^{\tilde{\mu}} W_{\tilde{\mu}} \Ket{\psi_0}
  = - 2 [2]^{-1}m^2 [s+1][s]\Ket{\psi_0}\,. 
\end{equation}
In accordance with the undeformed case we will call $s$ the spin of
the representation. The space of all rest states is stabilized by the
algebra generated by the little algebra for the massive case, $\suq$,
and the momenta, that is, by the $q$-Euclidean algebra $\Euk$. The
observables that are most commonly diagonalized are all elements of
$\Euk$: energy $P_0$, momentum $\vec{P}$, angular momentum $\vec{J}$,
helicity $\vec{J}\cdot\vec{P}$. We opt for an angular momentum basis,
where we diagonalize $J_3$ and $\vec{J}^2=\vec{J}\cdot\vec{J}$. If we
add the Casimir operators of $\Euk$, $P_0$ and $\Heli$ as defined in
Eq.~\eqref{eq:Helidef}, we get a complete set of commuting
observables.\footnote{The authors of \cite{Pillin:1993,Pillin} failed
  to add $Z$ or $\vec{J}\cdot\vec{P}$ to their set of commuting
  observables (cf. \cite{Pillin}, p. 67). This is the reason why they
  only found spin zero representations.}  Instead of $J_3$ and
$\vec{J}^2$ it is more practical to work with $K$ and $W$, whose
possible eigenvalues can be looked up in Sec.~\ref{sec:AppRep1}. From
Sec.~\ref{sec:Euk} we know that $P_0$ and $\Heli$ are Casimir
operators of a $\suq$ algebra, so we know their possible eigenvalues,
as well. Labeling the states of the yet to be constructed
representation by their possible eigenvalues we get including the
Casimirs
\begin{subequations}
\label{eq:Eigenwerte}
\begin{align}
\label{eq:Eigenwert1}
  K\Ket{j,m,n,k} &= q^{2m} \Ket{j,m,n,k} \\
\label{eq:Eigenwert2}
  W\Ket{j,m,n,k} &=
  [2]^{-1}\bigl(q^{(2j+1)}+q^{-(2j+1)}\bigr)\Ket{j,m,n,k} \\
\label{eq:Eigenwert3}
  P_0\Ket{j,m,n,k} &=
  m [2]^{-1}\bigl(q^{(2n+1)}+q^{-(2n+1)}\bigr)\Ket{j,m,n,k}\\
\label{eq:Eigenwert4}
  \Heli \Ket{j,m,n,k} &=
  m[2]^{-1}\bigl(q^{(2k+1)}+q^{-(2k+1)}\bigr)\Ket{j,m,n,k} \\
\label{eq:Eigenwert5}
  P_\mu P^\mu  \Ket{j,m,n,k} &= m^2 \Ket{j,m,n,k} \\
\label{eq:Eigenwert6}
  W^{\tilde{\mu}} W_{\tilde{\mu}} \Ket{j,m,n,k}
  &= - 2 [2]^{-1}m^2 [s+1][s] \Ket{j,m,n,k} \,.
\end{align}
\end{subequations}
The eigenvalues of $W$, $P_0$, and $\Heli$ are all of the same form,
$\qi(j)$, $m\qi(n)$, and $m\qi(k)$, where
\begin{equation}
\label{eq:q1def}
  \qi(j) := [2]^{-1}\bigl(q^{(2j+1)}+q^{-(2j+1)}\bigr) \,.
\end{equation}
For the operators with a more obvious undeformed limit $J_3$,
$\vec{J}^2$, and $\vec{J}\cdot\vec{P}$ we get
\begin{equation}
\begin{aligned}
  J_3\Ket{j,m,n,k}
  &= \bigl(q^{m}[m] -\lambda [2]^{-2}[2j+2][2j]\bigr)\Ket{j,m,n,k} \\
  \vec{J}^{\,2}\Ket{j,m,n,k} &= [2]^{-2}[2j+2][2j] \Ket{j,m,n,k}\\ 
  (\vec{J}\cdot\vec{P})\Ket{j,m,n,k}
  &= \lambda[2]^{-2}\bigl([n+j+k+2][n+j-k]\\
  &\qquad +[n-j+k][n-j-k]\bigr) \Ket{j,m,n,k} \,,
\end{aligned}
\end{equation}
which shows why it is more efficient to work with $K$, $W$, and
$\Heli$ instead.

One further advantage of using an angular momentum basis is, that the
$q$-Wigner-Eckart theorem of Page~\pageref{th:WignerEckart} applies.
The problem of finding the matrix elements of 3-vector or scalar
operators with respect to rotations is reduced to finding the reduced
matrix elements. For 3-vector operators such as $P_A$, $J_A$, $R_A$,
and $S_A$ we get
\begin{equation}
\label{eq:Wigreda}
  \Braket{j',m',n',k'}{P_A}{j,m,n,k}
  = \CGC{1}{j}{j'}{A}{m}{m'} \rBraket{j',n',k'}{\vec{P}}{j,n,k} \,,
\end{equation}
while for scalars with respect to rotations such as $Z$, $W$, $U$, and $V$
we get
\begin{equation}
\label{eq:Wigredb}
  \Braket{j',m',n',k'}{Z}{j,m,n,k}
  = \delta_{mm'}\delta_{jj'} \rBraket{j',n',k'}{Z}{j,n,k} \,.
\end{equation}
The values of the $q$-Clebsch-Gordan coefficients that we will need
are given in Sec.~\ref{sec:AppClebsch1}. Useful relations for the
reduced matrix elements can be derived from Eq.~\eqref{eq:RacahRed},
which has been done explicitly in Eqs.~\eqref{eq:AppRedrel}.

\subsection{Representations of the $q$-Euclidean Algebra}

If we keep $n$ and $k$ constant, we fix the eigenvalues of the
Casimirs operators $P_0$ and $\Heli$ of the $q$-Euclidean algebra
$\Euk$. For constant $n$, $k$ we must thus get an irreducible
representation of $\Euk$. This irreducible representation of $\Euk$ on
the mass shell is by isomorphism~\eqref{eq:EukIso} simply the product
$D^n\otimes D^k$ of two representations of $\suq$. We describe them
briefly in terms of reduced matrix elements.

The reduced matrix element of $J_A$ can be read off
Eq.~\eqref{eq:AppJrep}, 
\begin{equation}
\label{eq:RedL}
  \rBraket{j}{\vec{J}}{j} = -[2]^{-1}\sqrt{[2j+2][2j]} \,.
\end{equation}
Due to the Clebsch-Gordan series~\eqref{eq:CGSeries} $j$ takes on the
values $\{|k-n|,|k-n|+1,\ldots, k+n \}$. Taking the matrix elements of
Eq.~\eqref{eq:Eigenwert4} we find
\begin{equation}
\label{eq:schlumpf1}
  \rBraket{j,n,k}{\vec{P}}{j,n,k} = m\lambda\,
    \frac{[k+n+j+2][j-k+n]-[k+n-j][j+k-n]}{[2]\sqrt{[2j+2][2j]}} \,.
\end{equation}
If we take the diagonal matrix elements of the relation $P_A
P_B\,\varepsilon^{AB}{}_C = -\lambda P_0 P_C$ and of
Eq.~\eqref{eq:Eigenwert6} we get, using Eqs.~\eqref{eq:AppRedrel}, two
equations for the reduced matrix elements from which we can eliminate
the $\rBraket{j}{\vec{P}}{j-1}\rBraket{j-1}{\vec{P}}{j}$ term
\begin{multline}
   [2]\sqrt{[2j+3][2j+1]}\,\rBraket{j}{\vec{P}}{j+1}
   \rBraket{j+1}{\vec{P}}{j} =\\
   [2j]\rBraket{j}{\vec{P}}{j}^2
   +\lambda E \sqrt{[2j+2][2j]}\rBraket{j}{\vec{P}}{j}
   -[2j+2](P^0P^0 - m^2)\,.
\end{multline}
Upon inserting Eq.~\eqref{eq:schlumpf1},
\begin{multline}
  \rBraket{j}{\vec{P}}{j+1}\rBraket{j+1}{\vec{P}}{j} = \\ -m^2 \lambda^2
  \frac{[k+n+j+2][k+n-j][k-n+j+1][n-k+j+1]}{[2][2j+2]\sqrt{[2j+3][2j+1]}}\,,
\end{multline}
and using Eq.~\eqref{eq:RedKonjugation} we finally get
\begin{squeezedsubequations}
\begin{multline}
  \rBraket{j+1,n,k}{\vec{P}}{j,n,k} = \\ m \lambda
  \frac{\sqrt{[k+n+j+2][k+n-j][k-n+j+1][n-k+j+1]}}{\sqrt{[2][2j+3][2j+2]}}
\end{multline}
\begin{multline}
  \rBraket{j-1,n,k}{\vec{P}}{j,n,k} = \\ - m \lambda
  \frac{\sqrt{[k+n+j+1][k+n-j+1][k-n+j][n-k+j]}}{\sqrt{[2][2j][2j-1]}} \,.     
\end{multline}
\end{squeezedsubequations}

\subsection{Possible Transitions of Energy and Helicity}

Next, we will determine the possible transitions of the quantum
numbers $n$ and $k$ under the action of the non-Euclidean generators.
To find restrictions on the possible transitions we consider
Eq.~\eqref{eq:UP-Rel-1} and the contraction of Eq.~\eqref{eq:RP-Rel-1}
with $g^{DC}P_{D}$ from the left
\begin{subequations}
\label{eq:uebiA}
\begin{align}
\label{eq:uebiA1}
  [2]^2 U P^0 &= [4] P^0 U - q^{-1}\lambda^2[2] (\vec{P}\cdot \vec{R})\\
  [2](\vec{P}\cdot \vec{R}) P_0
  &= 2 P^0 (\vec{P}\cdot \vec{R}) - q (\vec{P}\cdot \vec{P}) U \,.
\end{align}
\end{subequations} 
Taking the matrix elements of these equations yields a system of
linear equations
\begin{subequations}
\begin{alignat}{3}
  0&= \, & m \bigl([4]\qi(n')-[2]^2\qi(n)\bigr) & \lrAngle{U} &
    -  q^{-1}\lambda^2 [2] & \lrAngle{\vec{P}\cdot \vec{R}} \\
  0&= & m^2q (1-\qi(n')^2) & \lrAngle{U}& \,\, + \,\, 
     m\bigl( 2\qi(n')- [2]\qi(n)\bigr) & \lrAngle{\vec{P}\cdot \vec{R}}\,, 
\end{alignat}
\end{subequations}
where we have used the abbreviation $\lrAngle{U} :=
\Braket{j,m,n',k'}{U}{j,m,n,k}$ and analogously for
$\lrAngle{\vec{P}\cdot \vec{R}}$. For a nontrivial solution to exist,
the determinant of the coefficient matrix must vanish,
\begin{equation}
\label{eq:uebiC}
\begin{split}
  0 &\stackrel{!}{=}\,\, m^2 [2]^2
     (\qi(n')^2 - [2]\qi(n')\qi(n') +\qi(n)^2) + m^2\lambda^2\\
  & = m^2\lambda^4 \,
  \left[n+\tfrac{1}{2} - n'\right] \left[n-\tfrac{1}{2} - n'\right]
  \left[n+n'+\tfrac{1}{2}\right] \left[n+n'+\tfrac{3}{2}\right] \,,
\end{split}
\end{equation}
which is, since $n\geq 0$, precisely the case for
$n'=n\pm\tfrac{1}{2}$.

To obtain conditions on the transitions of $k$ we contract
Eqs.~\eqref{eq:UP-Rel-2} and \eqref{eq:RP-Rel-1} with $g^{DC}J_D$ from
the left
\begin{subequations}
\begin{align}
\label{eq:uebi1}
  [2]^2 U (\vec{J}\cdot \vec{P}) &= [4] (\vec{L}\cdot\vec{P}) U
    - q \lambda^2 [2] P_0 (\vec{L}\cdot\vec{R})
    + \I \lambda^2 [2]\, \vec{J}\cdot(\vec{P}\times\vec{R}) \\
\label{eq:uebi2}
  [2] (\vec{J}\cdot\vec{R}) P_0 &= -q[2](\vec{J}\cdot\vec{P})U 
  + [4] P_0 (\vec{J}\cdot\vec{R})
  -\I \lambda [2]\,\vec{J}\cdot(\vec{P}\times\vec{R}) \,.
\end{align}
Contracting Eq.~\eqref{eq:RP-Rel-2} with $J^B J^A$ from the right 
and eliminating the $\vec{P}\cdot\vec{R}$ term using
Eq.~\eqref{eq:uebiA1} yields
\begin{multline}
\label{eq:uebi3}
  \lambda (\vec{J}\cdot\vec{R})(\vec{J}\cdot\vec{P}) =
  q^2\{q\lambda [2] \,(\vec{J}\cdot\vec{P}) - \lambda^2 W P_0 \}
  (\vec{J}\cdot\vec{P})
  +  (\vec{J}\cdot\vec{J}) U P_0 \\
  - \{q (\vec{J}\cdot\vec{J}) P_0 - \lambda W (\vec{J}\cdot\vec{P}) \} U
  - 2\I q^{-1}[2]^{-1} \lambda W \, \vec{J}\cdot (\vec{P}\times\vec{R}) \,.
\end{multline}
\end{subequations}
Eliminating the
$\vec{J}\cdot (\vec{P}\times \vec{R})$ term from the last three
equations we obtain
\begin{subequations}
\label{eq:uebiB}
\begin{align}
  \lambda^2\{(\vec{J}\cdot\vec{R})\Heli - q\Heli (\vec{J}\cdot\vec{R})\}
  &= q (P_0 - W\Heli)U - U(P_0 - W\Heli)  \\
  \lambda^2\{(\vec{J}\cdot\vec{R})P_0 - q^{-1}P_0(\vec{J}\cdot\vec{R})\}
  &= q^{-1}(\Heli - W P^0)U- U(\Heli - W P^0) \,.
\end{align}
\end{subequations}
Again we take the matrix elements of these two equations
\begin{alignat}{3}
  0&=\,&\{[\qi(n)-\qi(j)\qi(k)] - q[\qi(n')-\qi(j)\qi(k')]\} &\lrAngle{U}& 
   \,+ \lambda^2 \{\qi(k) - q\qi(k')\}
   &\lrAngle{\vec{J}\cdot\vec{R}} \notag \\
\label{eq:uebiD}
  0&=\,&\{q[\qi(k)-\qi(j)\qi(n)] - [\qi(k')-\qi(j)\qi(n')]\} &\lrAngle{U}&
   \,+ \lambda^2 \{q\qi(n) - \qi(n')\}
   & \lrAngle{\vec{J}\cdot\vec{R}} \,.
\end{alignat}
Provided Eq.~\eqref{eq:uebiC} holds, the determinant condition for a
nontrivial solution is
\begin{align}
  0 &=[2]^2 \{[\qi(k')^2 - [2]\qi(k')\qi(k') +\qi(k)^2]
      - [\qi(n')^2 - [2]\qi(n')\qi(n') +\qi(n)^2]\} \notag\\
    &=[2]^2 [\qi(k')^2 - [2]\qi(k')\qi(k') +\qi(k)^2] + \lambda^2\notag\\
\label{eq:uebiE}
    & = \lambda^4 \,
  \left[k+\tfrac{1}{2} - k'\right] \left[k-\tfrac{1}{2} - k'\right]
  \left[k+k'+\tfrac{1}{2}\right] \left[k+k'+\tfrac{3}{2}\right] \,,
\end{align}
which is fulfilled precisely for $k'=k\pm\tfrac{1}{2}$. We conclude
that the possible transitions of the quantum numbers $n$ and $k$ are
$n\rightarrow n\pm\tfrac{1}{2}$ and $k\rightarrow k\pm\tfrac{1}{2}$.

\subsection{Dependence on Total Angular Momentum}

Eq.~\eqref{eq:uebiD} establishes a correspondence between the reduced
matrix elements of $\vec{J}\cdot\vec{R}$ and $U$. With
Eq.~\eqref{eq:RedL} we get for $j>0$
\begin{equation}
\label{eq:Aeins-Def}
\begin{split}
  \rBraket{j,n',k'}{\vec{R}}{j,n,k} &=
  \left(\frac{\qi(n)-q\qi(n')}{\qi(k)-q\qi(k')} - \qi(j) \right)
  \frac{[2]\rBraket{j,n',k'}{U}{j,n,k}}{\lambda^2\sqrt{[2j+2][2j]}}\\
  &=: A_1(n',k',n,k,j) \rBraket{j,n',k'}{U}{j,n,k} \,.
\end{split}
\end{equation}
The reduced matrix elements of Eq.~\eqref{eq:RP-Rel-1} between 
$\rBra{j+1,n',k'}$ and $\rKet{j,n,k}$, $\rBra{j-1,n',k'}$ and
$\rKet{j,n,k}$ yield
\begin{subequations}
\label{eq:Azweidrei-Def}
\begin{align}
  \rBraket{j+1,n',k'}{\vec{R}}{j,n,k} &=
     A_2(n',k',n,k,j) \rBraket{j,n',k'}{U}{j,n,k}\\
  \rBraket{j-1,n',k'}{\vec{R}}{j,n,k} &=
     A_3(n',k',n,k,j) \rBraket{j,n',k'}{U}{j,n,k} \,,
\end{align}
\end{subequations}
where
\begin{subequations}
\begin{align}
  A_2 &:= \frac{
     (\lambda \sqrt{\frac{[2j]}{[2j+2]}}\, A_1 - q)
     \rBraket{j+1,n',k'}{\vec{P}}{j,n',k'}}%
     {m\bigl([2] \qi(n) - \frac{[4]}{[2]}\qi(n')\bigr)
       + \lambda \sqrt{\frac{[2j+4]}{[2j+2]}}
      \rBraket{j+1,n',k'}{\vec{P}}{j+1,n',k'} } \\
  A_3 &:= \frac{
     -(\lambda \sqrt{\frac{[2j+2]}{[2j]}}\, A_1 + q)
     \rBraket{j-1,n',k'}{\vec{P}}{j,n',k'}}%
     {m\bigl([2]\qi(n) - m\frac{[4]}{[2]}\qi(n')\bigr)
       - \lambda \sqrt{\frac{[2j-4]}{[2j]}}
      \rBraket{j-1,n',k'}{\vec{P}}{j-1,n',k'} } \,.
\end{align}
\end{subequations}
This again can be used to calculate the reduced matrix elements of
Eq.~\eqref{eq:UP-Rel-2} between $\rBra{j+1,n',k'}$ and $\rKet{j,n,k}$ 
\begin{equation}
\label{eq:Avier-Def}
  \rBraket{j+1,n',k'}{U}{j+1,n,k} =
  A_4(n',k',n,k,j) \rBraket{j,n',k'}{U}{j,n,k} \,,
\end{equation}
where
\begin{multline}
  A_4 := \Bigl\{ \Bigl(\tfrac{[4]}{[2]} -
    \lambda^2 \sqrt{\tfrac{[2j]}{[2j+2]}}\, A_1 \Bigr)
    \rBraket{j+1,n',k'}{\vec{P}}{j,n',k'}
    - \lambda^2 A_2 \Bigl(mq\qi(n') \\
    -\sqrt{\tfrac{[2j+4]}{[2j+2]}}\,
      \rBraket{j+1,n',k'}{\vec{P}}{j+1,n',k'} \Bigr)  \Bigr\}
    [2]^{-1}\rBraket{j+1,n,k}{\vec{P}}{j,n,k}^{-1} \,.
\end{multline}
The calculation of the auxiliary functions $A_1$, $A_2$, $A_3$, and
$A_4$ is elementary but lengthy.\footnote{The calculation of the
  auxiliary functions has been done by computer algebra
  \cite{Mathematica}.} The results can be written most compactly
introducing the functions $u(n',k',n,k)$ and $v(n',k',n,k)$ by
\begin{equation}
\label{eq:uvDef}
\begin{aligned}
  u(n+\Delta n, k+\Delta k,n,k) &:= \Delta n\, (2n+1) + \Delta k\, (2k+1)\\
  v(n+\Delta n, k+\Delta k,n,k) &:= \Delta n\, (2n+1) - \Delta k\, (2k+1)\,,
\end{aligned}
\end{equation}
for $\Delta n , \Delta k = \pm\tfrac{1}{2}$, that is,
\begin{equation}
\begin{aligned}
  n'=n-\tfrac{1}{2}\,,\quad k'=k-\tfrac{1}{2} \quad&\Rightarrow\quad
  u = -n-k-1  \,,\quad v = -n+k \\
  n'=n-\tfrac{1}{2}\,,\quad k'=k+\tfrac{1}{2} \quad&\Rightarrow\quad
  u = -n+k  \,,\quad v = -n-k-1 \\
  n'=n+\tfrac{1}{2}\,,\quad k'=k-\tfrac{1}{2} \quad&\Rightarrow\quad
  u = n-k  \,,\quad v = n+k+1 \\
  n'=n+\tfrac{1}{2}\,,\quad k'=k+\tfrac{1}{2} \quad&\Rightarrow\quad
  u = n+k+1  \,,\quad v = n-k \\
\end{aligned}
\end{equation}
Using $u$ we can write $A_4$ as
\begin{equation}
\label{eq:Avier-Werte}
  A_4(n',k',n,k,j)
  =\frac{\sqrt{[j+u+2][j - u +1]}}{\sqrt{[j+u+1][j - u]}}
  = \frac{A_5(n',k',n,k,j+1)}{A_5(n',k',n,k,j)} \,,
\end{equation}
where
\begin{equation}
\label{eq:Afuenf-Werte}
  A_5(n',k',n,k,j) := \sqrt{[j+u+1][j - u]} \,.
\end{equation}
Defining 
\begin{equation}
\label{eq:ZweifachRed-Def}
  \rBraket{n',k'}{U}{n,k} :=
  \frac{\rBraket{j,n',k'}{U}{j,n,k}}{A_5(n',k',n,k,j)}\,,
\end{equation}
Eq.~\eqref{eq:Avier-Def} tells us by induction that
$\rBraket{n',k'}{U}{n,k}$ does not depend on $j$. 
With Eqs.~\eqref{eq:Aeins-Def} and~\eqref{eq:Azweidrei-Def} we
conclude that the $j$-dependence of all reduced matrix elements can be
absorbed in reduction coefficients according to 
\begin{equation}
\begin{aligned}
  \rBraket{j',n',k'}{U}{j,n,k} &=
  \BolliO{j'}{n'}{k'}{j}{n}{k}\rBraket{n',k'}{U}{n,k}\\
  \rBraket{j',n',k'}{\vec{R}}{j,n,k} &=
  \Bollil{j'}{n'}{k'}{j}{n}{k}\rBraket{n',k'}{U}{n,k} \,,
\end{aligned}
\end{equation}
if we define the coefficients as
\begin{subequations}
\begin{align}
  \BolliO{j'}{n'}{k'}{j}{n}{k} &:=
  \begin{cases}
    A_5(n',k',n,k,j)\,, & j'=j>0 \\
    0\,, & \text{else}
  \end{cases}\\
  \Bollil{j'}{n'}{k'}{j}{n}{k} &:=
  \begin{cases}
    A_3(n',k',n,k,j) A_5(n',k',n,k,j)\,, & j'=j-1 \\
    A_1(n',k',n,k,j) A_5(n',k',n,k,j)\,, & j'=j > 0 \\
    A_2(n',k',n,k,j) A_5(n',k',n,k,j)\,, & j'=j+1 \\
    0\,, & \text{else} \,.
  \end{cases}
\end{align} 
\end{subequations}
Explicitly, the formulas for the $B$-coefficients are
\begin{equation}
\label{eq:Bolli-Werte}
\begin{aligned}
  B_q^0(j',n',k'\,|\,j,n,k) &= \delta_{jj'}\sqrt{[j+u+1][j-u]}\\ 
  B_q^1(j-1,n',k'\,|\,j,n,k) &=
    -\frac{q^{-j}\sqrt{[2][j+v][j-v][j-u][j-u-1]}}
    {\lambda \sqrt{[2j][2j-1]}}\\
  B_q^1(j,n',k'\,|\,j,n,k) &=   
    -(q^{(j+1)}[j-v]-q^{-(j+1)}[j+v])
    \frac{\sqrt{[j+u+1][j-u]}}{\lambda\sqrt{[2j+2][2j]}} \\
  B_q^1(j+1,n',k'\,|\,j,n,k) &=  
    -\frac{q^{j+1}\sqrt{[2][j+v+1][j-v+1][j+u+2][j+u+1]}}
    {\lambda \sqrt{[2j+3][2j+2]}} \,,
\end{aligned}
\end{equation}
which can be written more compactly as
\begin{multline}
\label{eq:Bolli-kompaktWerte}
  B^\alpha_q(j',n',k'\,|\,j,n,k) =\\
  (-\lambda)^{-\alpha}\,
  \CGC{\alpha}{j'}{j}{0}{v}{v}\times
  \begin{cases}
    q^{-j}\sqrt{[j'-u+1][j'-u]}  \,,& j'= j-1\\
    \sqrt{[j'+u+1][j'-u]}        \,,& j'= j\\
    -q^{j+1}\sqrt{[j'+u+1][j'+u]} \,,& j'=j+1 \,.
  \end{cases}
\end{multline}

\subsection{Dependence on the other Quantum Numbers}

Using the $B$-coefficients, equations in the reduced matrix elements
of $R$, $U$ can be reduced further to equations in the double reduced
matrix elements $\rBraket{n',k'}{U}{n,k}$ as defined in
Eq.~\eqref{eq:ZweifachRed-Def}. We start by taking the matrix elements of
the $RR$-relations~\eqref{eq:RR-Rel}, $U^2 - \lambda^2\,(\vec{R}\cdot
\vec{R}) = 1$ and $R^A U-U R^A = 0$ between $\rBra{j,n,k}$ and
$\rKet{j,n,k}$. We obtain
\begin{subequations}
\label{eq:nk}
\begin{align}
\label{eq:nk1}
  \sum_{n',k'} A_6(n',k',n,k,j)
  \rBraket{n,k}{U}{n',k'}\rBraket{n',k'}{U}{n,k} &=1 \\
\label{eq:nk2}
  \sum_{n',k'} A_7(n',k',n,k,j)
  \rBraket{n,k}{U}{n',k'}\rBraket{n',k'}{U}{n,k} &=0\,,
\end{align}
\end{subequations}
where the summation indices run through $n'=n\pm\tfrac{1}{2}$,
$k'=k\pm\tfrac{1}{2}$ and
\begin{squeezedsubequations}
\begin{multline}
  A_6(n',k',n,k,j) := 
  \BolliO{j}{n}{k}{j}{n'}{k'}\BolliO{j}{n'}{k'}{j}{n}{k}\\
  - \lambda^2 \sum_{j'=j-1}^{j+1}
  (-1)^{j'-j}\sqrt{\tfrac{[2j'+1]}{[2j+1]}}
  \,\Bollil{j}{n}{k}{j'}{n'}{k'}\Bollil{j'}{n'}{k'}{j}{n}{k} 
\end{multline}
\begin{multline}
  A_7(n',k',n,k,j) := -\frac{\lambda}{[2]}\sqrt{[2j+2][2j]} 
  \Bigl\{\Bollil{j}{n}{k}{j}{n'}{k'}\BolliO{j}{n'}{k'}{j}{n}{k}\\
  -\BolliO{j}{n}{k}{j}{n'}{k'}\Bollil{j}{n'}{k'}{j}{n}{k} \Bigr\}\,.
\end{multline}
\end{squeezedsubequations}
The values of these coefficients are
\begin{equation}
\begin{aligned}
  A_6(n',k',n,k,j) &= 4\Delta k \, \Delta n \, [2][2k'+1][2n'+1] \\
  A_7(n',k',n,k,j) &= [2v][j+u+1][j-u]
    = \lambda^{-2}[2] [2v]\bigl(\qi(j)-\qi(u)\bigr)\,. 
\end{aligned}
\end{equation}
Eq.~\eqref{eq:nk2} must hold for all values of $j$, which turns out to
lead to two independent equations. Thus, Eqs.~\eqref{eq:nk} form a
system of three independent equations in four unknowns of the type
$\rBraket{n,k}{U}{n',k'}\rBraket{n',k'}{U}{n,k}$.  Eliminating two
unknowns in each equation we can interpret them as recursion relations
\begin{subequations}
\label{eq:nkRekurs2}
\begin{align}
\label{eq:nkRekurs2a}
  \rho(\mu,\nu) &= \rho(\mu,\nu-1) + [2\nu+2]\\
\label{eq:nkRekurs2b}
  \omega(\mu,\nu) &= \omega(\mu+1,\nu) + [2\mu]\\
\label{eq:nkRekurs2c}
  \omega(\mu+1,\nu) &= -\rho(\mu,\nu) + [\nu+\mu+2][\nu-\mu+1]
\end{align}
\end{subequations}
where we use the abbreviations $\mu := k-n$, $\nu := k+n$ and 
\begin{squeezedsubequations}
\label{eq:rho-omega-Def}
\begin{multline}
    \rho(\mu,\nu) := [2]^2 [2k+2][2k+1][2n+2][2n+1]\\
    \times \rBraket{n,k}{U}{n+\tfrac{1}{2},k+\tfrac{1}{2}}
    \rBraket{n+\tfrac{1}{2},k+\tfrac{1}{2}}{U}{n,k}
\end{multline}
\begin{multline} 
  \omega(\mu,\nu) := [2]^2 [2k+1][2k][2n+2][2n+1]\\
    \times \rBraket{n,k}{U}{n+\tfrac{1}{2},k-\tfrac{1}{2}}
    \rBraket{n+\tfrac{1}{2},k-\tfrac{1}{2}}{U}{n,k}\,.
\end{multline}
\end{squeezedsubequations}
In order to determine the initial conditions, we recall
Eq.~\eqref{eq:W-fix} which tells us that $n=0$ implies $k=s$. Hence,
matrix elements involving states with $n=0$ and $k\neq s$ have to
vanish, in particular
\begin{equation}
  \rho(s,s-1) = 0 \,.
\end{equation}
The solution of recursion relation~\eqref{eq:nkRekurs2a} with this
initial value is 
\begin{equation}
  \rho(s,\nu) = \sum_{\nu'=s}^\nu [2\nu'+2] = [\nu+s+2][\nu-s+1]
\end{equation}
where we used $\sum_{i'=a}^b [2i'+c] = [a+b+c][b-a+1]$. Inserting this
result in Eq.~\eqref{eq:nkRekurs2c} yields $\omega(s+1,\nu) = 0$.  The
solution of Eq.~\eqref{eq:nkRekurs2b} with this initial value is
\begin{equation} 
  \omega(\mu,\nu) = \sum_{\mu'=\mu}^s [2\mu'] = [\mu+s][s-\mu+1] \,.
\end{equation}
Inserting this again in Eq.~\eqref{eq:nkRekurs2c} results in
\begin{equation}
\label{eq:rho-omega-Werte}
\begin{aligned} 
  \omega(\mu,\nu) &= \omega(\mu) = [\mu+s][s-\mu+1] \\
  \rho(\mu,\nu)   &= \rho(\nu)   = [\nu+s+2][\nu-s+1]
\end{aligned}
\end{equation}
for $-s\leq\mu\leq s+1$ and $s-1\leq\nu$. At the border of this
half-closed strip in $\mu\nu$-space $\rho$ and $\omega$ vanish, so
there are no transitions to the outside. For an irreducible
representation we must not have two disconnected regions, hence,
$\rho$ and $\omega$ must vanish outside this strip. The allowed
quantum numbers form a strip in $nk$-space given by
\begin{xalignat}{2}
\label{eq:nkrange}
  |\mu| = |k - n| &\leq s \,, & \nu = n + k \geq s \,. 
\end{xalignat}

To derive from Eq.~\eqref{eq:rho-omega-Werte} formulas for the matrix
elements we need to take the $RS$-relations~\eqref{eq:RS-Rel} into
account. We begin with the matrix elements of $UV=VU$ between
$\rBra{j,n+\tfrac{1}{2},k+\tfrac{1}{2}}$ and
$\rKet{j,n-\tfrac{1}{2},k-\tfrac{1}{2}}$ using the conjugation $U^*=V$
to obtain
\begin{multline}
\label{eq:Udiagmat}
  \rBraket{n+\tfrac{1}{2},k+\tfrac{1}{2}}{U}{n,k}
  \overline{ \rBraket{n-\tfrac{1}{2},k-\tfrac{1}{2}}{U}{n,k} } =\\
  \overline{ \rBraket{n,k}{U}{n+\tfrac{1}{2},k+\tfrac{1}{2}} }
  \rBraket{n,k}{U}{n-\tfrac{1}{2},k-\tfrac{1}{2}} \,,
\end{multline}
which can be written as
\begin{equation}
\frac{ \overline{ \rBraket{\mu,\nu-1}{U}{\mu,\nu} } }{
       \rBraket{\mu,\nu}{U}{\mu,\nu-1} } =
\frac{ \overline{ \rBraket{\mu,\nu}{U}{\mu,\nu+1} } }{
       \rBraket{\mu,\nu+1}{U}{\mu,\nu} }\,.
\end{equation}
with $\mu := k-n$, $\nu := k+n$ as above. Reading this as recursion
relation, it follows that
\begin{subequations}
\label{eq:Konju0}
\begin{equation}
  \overline{\rBraket{\mu,\nu-1}{U}{\mu,\nu}} = \alpha_\mu
  \rBraket{\mu,\nu}{U}{\mu,\nu-1} \,,
\end{equation}
where the yet to be determined number $\alpha_\mu$ may depend on $\mu$
but not on $\nu$. Taking the matrix elements of $UU'=U'U$ between
$\rBra{j,n+\tfrac{1}{2},k+\tfrac{1}{2}}$ and
$\rKet{j,n-\tfrac{1}{2},k-\tfrac{1}{2}}$, it follows analogously that
\begin{equation}
  \overline{\rBraket{\mu,\nu}{U}{\mu-1,\nu}} = \beta_\nu
  \rBraket{\mu-1,\nu}{U}{\mu,\nu} \,,
\end{equation}
\end{subequations}
with $\beta_\nu$ independent of $\mu$.

Next, we take the diagonal matrix elements of $W = UV +
q^2\lambda^2(\vec{R} \cdot \vec{S})$ as in Eq.~\eqref{eq:Rot-RS}
using the conjugation relations~\eqref{eq:RedKonjugation} to obtain 
\begin{equation}
\label{eq:Aacht}
  \sum_{n',k'} A_8(n',k',n,k,j) |\rBraket{n,k}{U}{n',k'}|^2 = \qi(j) \,,
\end{equation}
where
\begin{equation}
\label{eq:Konju1a}
  A_8(n',k',n,k,j) :=
  |\BolliO{j}{n}{k}{j}{n'}{k'}|^2 + q^2 \lambda^2 \sum_{j'=j-1}^{j+1}
  |\Bollil{j}{n}{k}{j'}{n'}{k'}|^2  \,.
\end{equation}
Eq.~\eqref{eq:Aacht} must hold for all possible values of $j$, thus
yielding two independent equations from which we can derive
\begin{multline}
   [2]^{-2} [\mu+\nu+1]^{-1} =
   q^{-2\mu} [\nu-\mu] |\rBraket{\mu,\nu}{U}{\mu,\nu-1}|^2 \\
   + q^{-2(\nu+1)} [\nu-\mu+2] |\rBraket{\mu,\nu}{U}{\mu-1,\nu}|^2 \,.
\end{multline}
Relations~\eqref{eq:nkrange} tell us that the first term on the right
hand side vanishes for $\nu=s$ while the second vanishes for $\mu=-s$,
that is,
\begin{equation}
\begin{aligned}
  |\rBraket{s,\nu}{U}{s,\nu-1}|^2
  &= \frac{q^{-2s}}{[2]^2 [\nu-s+1][\nu+s]} \\
  |\rBraket{\mu,s}{U}{\mu-1,s}|^2
  &= \frac{q^{2(s+1)}}{[2]^2 [\mu+s+1][s-\mu+2]} \,.
\end{aligned}
\end{equation}
If we compare this with $\rho(-s,\nu-1)$ and $\omega(\mu,s)$ as
computed in Eqs.~\eqref{eq:rho-omega-Werte}, we find
\begin{xalignat}{2}
  \alpha_\mu &= q^{2s}\,, & \beta_\nu = q^{2(s+1)} \,.
\end{xalignat}
With this result Eqs.~\eqref{eq:rho-omega-Werte} can be written as
formulas for the squares of matrix elements. For example,
\begin{align}
  [\mu&+s][s-\mu+1] = \omega(\mu,\nu) \notag\\
  &=  [2]^2 [2k+1][2k][2n+2][2n+1]
    \rBraket{n,k}{U}{n+\tfrac{1}{2},k-\tfrac{1}{2}}
    \rBraket{n+\tfrac{1}{2},k-\tfrac{1}{2}}{U}{n,k}\notag\\
  &=  q^{2(s+1)}[2]^2 [2k+1][2k][2n+2][2n+1]
    \,\,|\rBraket{n+\tfrac{1}{2},k-\tfrac{1}{2}}{U}{n,k}|^2 \,.
\end{align}
This is an equation for the absolute value of the double reduced matrix
elements. In fact, none of the commutation relations of the
$q$-Poincar{\'e} algebra gives us a condition on the phase of the reduced matrix
elements, that is, the phase can be chosen arbitrarily. We choose it,
such that
\begin{equation}
  \rBraket{n+\tfrac{1}{2},k-\tfrac{1}{2}}{U}{n,k} = \frac{q^{-2(s+1)}
    \sqrt{[s+k-n][s-k+n+1]}}{[2]\sqrt{[2k+1][2k][2n+2][2n+1]}} \,.
\end{equation}
Analogously, we determine the other matrix elements. The end result is
\begin{equation}
\label{eq:dopred-Werte}
  \rBraket{n',k'}{U}{n,k} = \frac{ q^{2(n-n')s + (n'-k'-n+k)}
    \sqrt{[s+u+1][s-u]}}{[2]\sqrt{[k'+k+\frac{3}{2}][k'+k+\frac{1}{2}]
      [n'+n+\frac{3}{2}][n'+n+\frac{1}{2}]}} \,.
\end{equation}

\subsubsection{Summary} We summarize the results for the reduced
matrix elements. As before, the abbreviations $u$ and $v$ as defined
in Eq.~\eqref{eq:uvDef} are being used. The relation between the
reduced and the ordinary matrix elements is given by
Eqs.~\eqref{eq:Wigreda} and \eqref{eq:Wigredb}.
\begin{squeezedsubequations}
\begin{multline}
  \rBraket{j',n',k'}{\vec{J}}{j,n,k}
  = -[2]^{-1}\delta_{jj'}\delta_{nn'}\delta_{kk'}\sqrt{[2j+2][2j]}
  \\[-4\lineskip]
\end{multline}
\begin{multline}
  \rBraket{j-1,n',k'}{\vec{P}}{j,n,k} =
  -m\lambda\delta_{nn'}\delta_{kk'}\\ \times
  \frac{\sqrt{[k+n+j+1][k+n-j+1][k-n+j][n-k+j]}}{\sqrt{[2][2j][2j-1]}}
\end{multline}
\begin{multline}
  \rBraket{j,n',k'}{\vec{P}}{j,n,k} =
  m\lambda\delta_{nn'}\delta_{kk'}\\ \times
    \frac{[k+n+j+2][j-k+n]-[k+n-j][j+k-n]}{[2]\sqrt{[2j+2][2j]}}
\end{multline}
\begin{multline}
  \rBraket{j+1,n,k}{\vec{P}}{j,n,k} =
  m \lambda\delta_{nn'}\delta_{kk'}\\ \times
  \frac{\sqrt{[k+n+j+2][k+n-j][k-n+j+1][n-k+j+1]}}{\sqrt{[2][2j+3][2j+2]}}
\end{multline}
\begin{multline}
  \rBraket{j',n',k'}{U}{j,n,k} =
    \delta_{jj'}\, q^{2(n-n')s + (n'-k'-n+k)} \\
    \times\frac{\sqrt{[j+u+1][j-u][s+u+1][s-u]}}
    {[2]\sqrt{[k'+k+\frac{3}{2}][k'+k+\frac{1}{2}]
      [n'+n+\frac{3}{2}][n'+n+\frac{1}{2}]}} \,.
\end{multline}
\begin{multline}
  \rBraket{j',n',k'}{\vec{R}}{j,n,k} =\\
    \frac{ q^{2(n-n')s + (n'-k'-n+k)}
    \sqrt{[s+u+1][s-u]}}{\lambda[2]
    \sqrt{[k'+k+\frac{3}{2}][k'+k+\frac{1}{2}]
      [n'+n+\frac{3}{2}][n'+n+\frac{1}{2}]}}\\
  \times\CGC{1}{j'}{j}{0}{v}{v}\times
  \begin{cases}
    -q^{-j}\sqrt{[j'-u+1][j'-u]}  \,,& j'= j-1\\
    -\sqrt{[j'+u+1][j'-u]}        \,,& j'= j\\
    q^{j+1}\sqrt{[j'+u+1][j'+u]} \,,& j'=j+1 \,.
  \end{cases}
\end{multline}
\begin{multline}
  \rBraket{j',n',k'}{V}{j,n,k} =
    \delta_{jj'}\, q^{2(n'-n)s + (n-k-n'+k')} \\
    \times\frac{\sqrt{[j+u+1][j-u][s+u+1][s-u]}}
    {[2]\sqrt{[k'+k+\frac{3}{2}][k'+k+\frac{1}{2}]
      [n'+n+\frac{3}{2}][n'+n+\frac{1}{2}]}} \,.
\end{multline}
\begin{multline}
  \rBraket{j',n',k'}{\vec{S}}{j,n,k} =\\
    \frac{ q^{2(n'-n)s + (n-k-n'+k') }
    \sqrt{[s+u+1][s-u]}}{\lambda[2]
    \sqrt{[k'+k+\frac{3}{2}][k'+k+\frac{1}{2}]
      [n'+n+\frac{3}{2}][n'+n+\frac{1}{2}]}}\\
  \times\CGC{1}{j'}{j}{0}{-v}{-v}\times
  \begin{cases}
    -q^{j}\sqrt{[j'-u+1][j'-u]}  \,,& j'= j-1\\
    -\sqrt{[j'+u+1][j'-u]}        \,,& j'= j\\
    q^{-(j+1)}\sqrt{[j'+u+1][j'+u]} \,,& j'=j+1
  \end{cases}
\end{multline}
\end{squeezedsubequations}

\section{Representations by Induction}

We want to describe briefly how representations of the $q$-Poincar{\'e}
algebra can be constructed using the method of induced
representations.

\subsection{The Method of Induced Representations of Algebras}

Let us assume that we do have an irreducible representation of the
undeformed Poincar{\'e} algebra $\mathcal{P}$ on a Hilbert space
$\mathcal{H}$,
\begin{equation}
  \sigma:\mathcal{P}\otimes\mathcal{H}\longrightarrow \mathcal{H} \,.
\end{equation}
Let the situation be as in Sec.~\ref{sec:Little1}, where we denoted by
$\mathcal{H}_p$ a momentum eigenspace and by $\mathcal{K}_p$ its
stabilizer (little algebra). By definition, the restriction of $\sigma$ to
$\mathcal{H}_p$ defines representations on translations $\mathcal{T}$ and the
little algebra $\mathcal{K}_p$ by
\begin{equation}
\label{eq:Induced1}
\begin{aligned}
  \chi_p&:\mathcal{T}\longrightarrow \mathbb{R} \,,\quad
  \text{where}\quad \sigma(t\otimes\Ket{\psi_p}) = \chi_p(t)\Ket{\psi_p} \\
  \rho&:\mathcal{K}_p \otimes \mathcal{H}_p
  \longrightarrow \mathcal{H}_p \,,\quad
  \rho(k\otimes\Ket{\psi_p}) = \sigma(k\otimes\Ket{\psi_p})   
\end{aligned}
\end{equation}
for all $\Ket{\psi_p} \in \mathcal{H}_p$. Together, $\chi_p$ and
$\rho$ define a representation of $\mathcal{T}\rtimes\mathcal{K}_p$
on $\mathcal{H}_p$.  Let us assume for a moment that we did not know
about $\sigma$ but were given only $\chi_p$ and $\rho$. There is a
generic method to extend a representation of an subalgebra to a
representation of the whole algebra.
\begin{Definition}
  Let $\mathcal{A}$ be an algebra, $\mathcal{S}$ a subalgebra and $V$
  a left $\mathcal{S}$-module. Then the tensor product of
  $\mathcal{A}$ and $V$ over $\mathcal{S}$,
  $\mathcal{A}\otimes_\mathcal{S} V$ becomes a left
  $\mathcal{A}$-module by left multiplication. It is called the module
  (or representation) induced by $V$.

  Explicitly, $\mathcal{A}\otimes_\mathcal{S} V$ is the vector space
  $\mathcal{A}\otimes V$ (ordinary tensor product over the complex
  numbers), divided by the relations
  \begin{equation}
    as \otimes v = a \otimes sv \,,\quad\text{for all}\quad
    a\in\mathcal{A}\,,s\in\mathcal{S}\,, v\in V \,,
  \end{equation}
  with the left $\mathcal{A}$-action defined by
  \begin{equation}
    a'(a\otimes v) = a'a \otimes v 
  \end{equation}
  and linear extension.
\end{Definition}
For given $\chi_p$, $\mathcal{H}_p$, $\rho$, and $\mathcal{K}_p$, the
induced representation acts on the tensor product
\begin{equation}
  \mathcal{P}\otimes_{\mathcal{T}\rtimes\mathcal{K}_p}\mathcal{H}_p
  = (\mathcal{T}\rtimes\mathcal{L})
  \otimes_{\mathcal{T}\rtimes\mathcal{K}_p}\mathcal{H}_p
  \cong \mathcal{L}\otimes_{\mathcal{K}_p}\mathcal{H}_p. 
\end{equation}
While this construction may look somewhat abstract, its great
practical value lies in the following
\begin{Theorem}
\label{th:Blattner}
Let $\mathcal{P}= \mathcal{T}\rtimes\mathcal{L}$ be the Poincar{\'e}
algebra, $\chi_p$ a one dimensional representation of $\mathcal{T}$,
$\mathcal{K}_p = \{k\in \mathcal{L}\,|\, \chi_p([k,t]) = 0
\quad\text{for all}\quad t\in\mathcal{T} \}$ the according little
algebra, and $\rho$ an \emph{irreducible} representation of
$\mathcal{K}_p$ on the finite vector space $\mathcal{H}_p$. With the
action defined by $\chi_p$ and $\rho$ the space $\mathcal{H}_p$
becomes a left $\mathcal{T}\rtimes\mathcal{K}_p$-module. Then the
induced representation
$\mathcal{P}\otimes_{\mathcal{T}\rtimes\mathcal{K}_p}\mathcal{H}_p$ is
irreducible. Furthermore, all irreducible representations of
$\mathcal{P}$ are of this form \cite{Blattner:1969,Dixmier}.
\end{Theorem}
This means that all we have to do in order to construct the
irreducible representations of $\mathcal{P}$ is
\begin{enumerate}
\item determine the little algebras,
\item construct the irreducible representations of the little algebras,
\item induce these representations.
\end{enumerate}
Using the Lie group version of this method, Wigner \cite{Wigner:1939}
was the first to construct all irreducible representations of the
Poincar{\'e} group (see also \cite{Mackey}). Theorem~\ref{th:Blattner}
cannot be generalized to Hopf semidirect products but in very special
cases \cite{Blattner:1986,Schneider:1990,Schneider:1994}. The method
of induced representations, however, works for any algebra.

\subsection{Induced Representations of the $q$-Poincar{\'e} Algebra}

We will deal only with the massive case, $p = (p_\mu) = (m,0,0,0) =
\chi_p(P_\mu)$, where we have $\mathcal{K}_p = \suq$, as calculated in
Sec.~\ref{sec:Little2}.  Let $D^j$ be an irreducible $\suq$-module.
Recall (p.~\pageref{th:Poincare}) the definition of the $q$-Poincar{\'e}
algebra $\mathcal{P}_q = \Mink_q \rtimes \slC$. In the quantum double
form (Sec.~\ref{sec:QuantumDouble}) the $q$-Lorentz algebra is $\slC
\cong \SUq^\op\otimes \suq$ as vector space. We conclude that the
induced representation of $D^j$ acts on the vector space
\begin{equation}
\begin{split}
  \mathcal{P}_q \otimes_{\Mink\rtimes\suq} D^j
  &= [\Mink\rtimes\slC] \otimes_{\Mink\rtimes\suq} D^j \\
  &\cong \slC \otimes_{\suq} D^j \\
  &\cong (\SUq^\op\otimes \suq) \otimes_{\suq} D^j \\
  &\cong \SUq^\op \otimes D^j.
\end{split}
\end{equation}
Let ${e_m}$ be a basis of $D^j$.  The action of some boost
$h'\in\SUq^\op$ on $h\otimes e_m \in \SUq^\op \otimes D^j$ is simply
given by left multiplication
\begin{equation}
  h'(h\otimes e_m) = h'h \otimes e_m \,.
\end{equation}
For the action of a rotation $l\in\suq$ we have to commute $lh$ using
Eq.~\eqref{eq:RotBoostCommute} and let $l$ act on $e_m$
\begin{equation}
 l (h\otimes e_m)
 = \lrAngle{ S(l_{(1)}), h_{(1)} } \lrAngle{l_{(3)}, h_{(3)}} \,
  (h_{(2)} \otimes e_{m'}\,\rho^j(l_{(2)})^{m'}{}_m) \,.
\end{equation}
Finally, for the action of $P_\mu \in\Mink$ we must use
Eq.~\eqref{eq:PoincCommute},
\begin{equation}
  P_\mu (h\otimes \psi) = 
  p_{\mu'} \Lambda(S^{-1}h_{(1)})^{\mu'}{}_\mu \,
  (h_{(2)} \otimes e_m) \,,
\end{equation}
where $p_\mu = \chi_p(P_\mu)$ are the momentum eigenvalues. 

We can equip this representation with a scalar product using the Haar
measure of $\SUq$ (\cite{Woronowicz:1987}, see also \cite{Schmuedgen},
pp.~111-117). An orthogonal basis is provided by the Peter-Weyl
theorem (\cite{Schmuedgen}, pp.~106-111).

\chapter{Free Wave Equations}
\label{sec:WaveEqs}

\section{General Wave Equations}
\label{sec:MainContrib4a}

\subsection{Wave Equations by Representation Theory}

On the way from free theories to theories with interaction we need to
leave the mass shell. The space of on-shell states is clearly too
small as to allow for interactions where energy and momentum can be
transfered from one sort of particle onto another. Moreover, we need a
way to describe several particle types and their coupling in one
common formalism.

These issues are resolved by introducing Lorentz spinor wave
functions, that is, tensor products of the algebra of functions on
spacetime with a finite vector space containing the spin degrees of
freedom, the whole space carrying a tensor representation of the
Lorentz symmetry. The additional mathematical structure we need to
describe coupling is provided by the multiplication within the algebra
of space functions. This structure is equally present in the
undeformed as in the deformed case.

Using such Lorentz spinors has some consequences that have to be
dealt with:
\begin{itemize}
\item[(a)] The Lorentz spinor representations cannot be irreducible.
  Otherwise they would have to be on shell and the spinor degrees of
  freedom would have to carry a representation of the little algebra.
\item[(b)] The Lorentz spinor representations cannot be unitary since
  the spin degrees of freedom carry a finite representation of the
  non-compact Lorentz algebra.
\end{itemize}
The solution to these problems are:
\begin{itemize}
\item[(a)] We consider only an irreducible subrepresentation to be the
  space of physical states. This subrepresentation is described as
  kernel of a linear operator $\mathbb{A}$, that is, we demand all
  physical states $\psi$ to satisfy the wave equation $\mathbb{A}\psi
  = 0$.
\item[(b)] We introduce a non-degenerate but indefinite pseudo scalar
  product, such that the spinor representation becomes a
  $*$-representation with respect to the corresponding pseudo adjoint.
  This amounts to introducing a new conjugation $j$ on states and
  operators. 
\end{itemize}
For $\ker \mathbb{A}$ to be a subrepresentation, the operator must
satisfy
\begin{equation}
\label{eq:WaveCondition1}
  \mathbb{A}\psi = 0 \quad\Rightarrow\quad \mathbb{A}h\psi = 0
\end{equation}
for all $q$-Poincar{\'e} transformations $h$. Depending on the particle
type under consideration we might include charge and parity
transformations. $\mathbb{A}$ is not unique since the wave
equations for $\mathbb{A}$ and $\mathbb{A}'$ must be considered
equivalent as long as their solutions are the same, $\ker(\mathbb{A})
= \ker(\mathbb{A}')$. 

Ideally, the operator $\mathbb{A}$ is a projector, $\mathbb{A}=\Proj$,
with $\Proj^2=\Proj$, $\Proj^* = \Proj$.
Condition~\eqref{eq:WaveCondition1} is then equivalent to
\begin{equation}
\label{eq:WaveCondition2}
  [\Proj,h] = 0
\end{equation}
for all $q$-Poincar{\'e} transformations $h$. Whether the wave equation is
written with a projection is a matter of convenience. The Dirac
equation is commonly written with such a projection which is determined
uniquely (up to complement) by condition~\eqref{eq:WaveCondition2}.
For the Maxwell equations a projection can be found but yields a second
order differential equation. For this reason, the Maxwell equations
are commonly described by a more general operator $\mathbb{A}$, which
leads to a first order equation. So far, all considerations pertain
equally to the deformed as to the undeformed case.

\subsection{$q$-Lorentz Spinors}

We define a general, single particle $q$-Lorentz spinor wave function
as element of the tensor product $\mathcal{S}\otimes \Mink$ of a
finite vector space $\mathcal{S}$ holding the spin degrees of freedom
and the space of $q$-Minkowski space functions $\Mink$
(Sec.~\ref{sec:MinkowskiConstruct}).

Let $\{ e_k \}$ be a basis of $\mathcal{S}$ transforming under a
$q$-Lorentz transformation $h\in\slC$ as $h\tr e_j = e_i
\,\rho(h)^i{}_j$, where $\rho : \slC \rightarrow
\mathrm{End}(\mathcal{S})$ is the representation map. Any spinor
$\psi$ can be written as
\begin{equation}
  \psi = e_j \otimes \psi^j \,,
\end{equation}
where $j$ is summed over and the $\psi^j$ are elements of $\Mink$.
The total action of $h\in\slC$ on a spinor is
\begin{equation}
  h\psi = (h_{(1)}\tr e_j) \otimes (h_{(2)}\tr \psi^j)
  = e_i\otimes \rho(h_{(1)})^i{}_j  (h_{(2)}\tr \psi^j) \,.
\end{equation}
This tells us that, if we want to work directly with the
$\Mink$-valued components $\psi^j$, the action of $h$ is
\begin{equation}
\label{eq:Wave1}
  h\psi^i = \rho(h_{(1)})^i{}_j  (h_{(2)}\tr \psi^j) \,.
\end{equation}
Do not confuse the total action $h\psi^i$ with the action of $h$
on each component of $\psi^i$ denoted by $h\tr \psi^i$. The
transformation of $\psi^i$ can easily be generalized to the case where
$\mathcal{S}$ carries a tensor representation of two finite
representations, that is, we have spinors with two or more indices
\begin{equation}
\label{Wave2}
  h\psi^{ij} = \rho(h_{(1)})^i{}_{i'}
                  \rho'(h_{(2)})^j{}_{j'}  (h_{(3)}\tr \psi^{i'j'}) \,,
\end{equation}
where $\rho$ and $\rho'$ are the representation maps of the first and
second index, respectively. 

Furthermore, we can derive spinors by the action of tensor operators:
Let $T^i$ be a upper left $\rho$-tensor operator and $\psi = e_j
\otimes \psi^j$ a $\rho'$-spinor field. Any operator $T^i$ can be
written as $T^i = \sum_k A_k^i \otimes B_k^i \in
\mathrm{End}(\mathcal{S})\otimes \mathrm{End}(\Mink)$ such that the
action of $T^i$ becomes
\begin{equation}
  T^i \psi
  = e_j \otimes \sum_k \rho(A_k^i)^j{}_{j'} B_k^i \tr \psi^{j'}
  = e_j \otimes (T^i \psi^j) =: e_j \otimes  \phi^{ij} \,.
\end{equation}
How does this new array of wave functions $\phi^{ij} = T^i\psi^j$
transform under $q$-Lorentz transformations? Letting act $h$ from the
left, we find
\begin{equation}
  h\phi^{ij} = \rho(h_{(1)})^j{}_{j'}(h_{(2)}\tr \phi^{ij'}) \,,
\end{equation}
that is, $h$ acts only on the index that came from the wave functions
$\psi^j$. However, if we transform $\phi^{ij}$ by transforming
$\psi^j$ inside, we find
\begin{equation}
\begin{split}
  T^i (h\psi^j) &= (T^i h)\psi^j 
  = h_{(2)}  [ \adL S^{-1}(h_{(1)}) \tr T^i ] \psi^{j} \\ 
  &= \rho(h_{(1)})^i{}_{i'}  h_{(2)}T^{i'} \psi^{j} 
  = \rho(h_{(1)})^i{}_{i'}  h_{(2)} \phi^{i'j} \\
  &= \rho(h_{(1)})^i{}_{i'} \rho'(h_{(2)})^j{}_{j'} 
     (h_{(3)} \tr \phi^{i'j'})  \,.
\end{split}
\end{equation}
In other words, if $\psi^j$ is transformed $\phi^{ij} = T^i\psi^j$
will transform as a $\rho\otimes\rho'$-spinor. Note, that for the last
calculation the order in the tensor product $\mathcal{S}\otimes \Mink$
is essential. This reasoning would not have worked out as nicely if we
had constructed the spinor space as $\Mink\otimes\mathcal{S}$.  Chief
examples of this construction are the gauge term $P^\mu\phi$ of the
vector potential $A^\mu$, or the derivatives of the vector potential
$P^\mu A^\nu$ which are used to construct the electromagnetic field
strength tensor $F^{\mu\nu}$.

We have not said yet how the momenta $P^\mu$ act on $q$-Lorentz
spinors. One might be tempted to assume that they act on the wave
function part only, that is, as $1\otimes P^\mu$ on the tensor
product. However, this is not possible, as in general $1\otimes
P^\mu$ is no $4$-vector operator and thus cannot represent 4-momentum.
We can turn $1\otimes P_\mu$ into a vector operator, though, by
twisting $1\otimes P_\mu$ with an $\R$-matrix of the $q$-Lorentz
algebra,
\begin{equation}
  P^\mu := \R^{-1}(1\otimes P^\mu)\R =
           (L^\Lambda_+)^{\mu}{}_{\mu'} \otimes P^{\mu'} \,,
\end{equation}
with the $L$-matrix for the $4$-vector representation as defined in
Eq.~\eqref{eq:Lmatdef}. Of the two universal $\R$-matrices of the
$q$-Lorentz algebra we opt for $\RI$, because only then the twisting
is compatible with the $*$-structure. The momenta act on a
$\rho$-spinor as
\begin{equation}
  P^\mu \psi^i
  = \rho\bigl( (L^\Lambda_{\mathrm{I}+})^{\mu}{}_{\mu'}\bigr)^i{}_j
  \,(P^{\mu'} \tr \psi^j) \,,
\end{equation}
where the $L$-matrix has been calculated in Eq.~\eqref{eq:LLambda}.
The action of $P^\mu$ on each component of $\psi^j$ can be viewed as
derivation within the algebra of $q$-Minkowski space functions
$\Mink$. The $q$-derivation operators are
\begin{equation}
  \partial^\mu := 1\otimes \I P^\mu \,.
\end{equation}
Now we can interpret an operator linear in the momenta as
$q$-differential operator. If $C_\mu = C_\mu\otimes 1$ are operators
that act on the spinor indices only,
\begin{equation}
  \I\,C_\mu  P^\mu
  = C_\mu\,\rho\bigl( (L^\Lambda_{\mathrm{I}+})^{\mu}{}_{\mu'}\bigr)
  \partial^{\mu'} = \tilde{C}_{\mu'} \partial^{\mu'} \,,
\end{equation} 
where
\begin{equation}
\label{eq:Optilde}
  \tilde{C}_{\mu'} :=
  C_\mu\,\rho\bigl( (L^\Lambda_{\mathrm{I}+})^{\mu}{}_{\mu'}\bigr)
\end{equation}
such that $\tilde{C}_{\mu'}$ still acts on the spinor index only,
while $\partial^{\mu'}$ acts componentwise, so the two operators
commute $[\tilde{C}_{\mu},\partial^{\nu}] = 0$. It remains to
calculate the transformation $C_\mu \rightarrow \tilde{C}_\mu$ for
particular representations $\rho$. Finally, we remark that for the
mass Casimir we have $P_\mu P^\mu = \R^{-1}(1\otimes P_\mu P^\mu)\R =
1\otimes P_\mu P^\mu$, hence, $P_\mu P^\mu = - \partial_\mu
\partial^\mu$. This means, that mass irreducibility for a spinor is
the same as mass irreducibility for each component of the spinor.

\subsection{Conjugate Spinors}

One of the effects of using Lorentz spinors is that the underlying
representations can no longer be unitary, since there are no unitary
finite representations of the non-compact Lorentz algebra. However, we
can introduce non-degenerate but indefinite bilinear forms playing the
role of the scalar product. With respect to these pseudo scalar
products the spinors carry $*$-representations, that is, the
$*$-operation on the algebra side is the same as the pseudo adjoint on
the operator side.

The problem of non-unitarity arises from the finiteness of the spin
part, $\mathcal{S}$, within the space of spinor wave functions
$\mathcal{S}\otimes \Mink$, so we can assume that the wave function
part $\Mink$ does carry a $*$-representation. It is then sufficient to
redefine the scalar product on $\mathcal{S}$ only. Consider a
$D^{(j,0)}$-representation of $\slC$ with orthonormal basis $\{e_m\}$
and the canonical scalar product $\lrAngle{e_m | e_n } = \delta_{mn}$.
We want to define a pseudo scalar product by
\begin{equation}
\label{eq:pseudounitary}
  ( e_m | e_n ) := A_{mn} \quad\text{such that}\quad
  ( e_m | (g\otimes h)\tr e_n ) = ((g\otimes h)^*\tr e_m | e_n )
\end{equation}
for any $g\otimes h \in \slC$. For a pseudo scalar product we must
suppose $A_{mn}$ to be a non-degenerate, hermitian, but not
necessarily positive definite matrix. Inserting the definition of the
pseudo scalar product, the pseudo-unitarity
condition~\eqref{eq:pseudounitary} reads
\begin{equation}
\label{eq:pseudounitary2}
\begin{split}
  ( e_m | (g\otimes h) \tr e_n )
  &= ( e_m |\,e_{n'} \rho^j(g)^{n'}{}_{n}\,\varepsilon(h) )
  = A_{mn'} \rho^j(g)^{n'}{}_{n}\,\varepsilon(h) \\
  &{\stackrel{!}{=}} ( (g\otimes h)^* \tr e_m |\, e_n )
  = ( e_{m'}\varepsilon(g^*) \rho^j(h^*)^{m'}{}_{m} |\, e_{n} ) \\
  &= A_{m'n} \,\overline{\varepsilon(g^*)\rho^j(h^*)^{m'}{}_{m} }
  =  A_{m'n} \,\varepsilon(g)\rho^j(h)^{m}{}_{m'} \,,
\end{split}
\end{equation}
where we have used the definition~\eqref{eq:LorentzHopf} of $(g\otimes
h)^*$ observing that $\varepsilon(\R_{[1]})\R_{[2]}=1$. Traditionally,
the scalar product is not described by a matrix $A_{mn}$ but by
introducing a conjugate spinor basis $\{\bar{e}_m\}$ demanding
\begin{equation}
  ( e_m | e_n ) = \lrAngle{ \bar{e}_m | e_n } \quad\Rightarrow\quad
  \bar{e}_m = e_{m'} A_{m'm} \,.
\end{equation}
Using~\eqref{eq:pseudounitary2} the conjugate basis turns out to
transform as
\begin{equation}
\begin{split}
  (g\otimes h) \tr \bar{e}_n 
  &= e_{m'} \rho^j(g)^{m'}{}_{m} \varepsilon(h) \,A_{mn}
  =e_{m'} A_{m'n'} \varepsilon(g)\rho^j(h)^{n'}{}_{n} \\
  &= \bar{e}_{n'} \varepsilon(g) \rho^j(h)^{n'}{}_{n} \,,
\end{split}
\end{equation}
that is, $\bar{e}_m$ ought to transform according to a
$D^{(0,j)}$-representation. $D^{(j,0)}$ and $D^{(0,j)}$ being
inequivalent representations, the conjugate basis $\bar{e}_m$ cannot
be expressed as a linear combination of the original basis vectors
$e_m$. In order to allow for a conjugate spinor basis we must consider
a representation that contains both, $D^{(j,0)}$ and $D^{(0,j)}$, and
thus at least their direct sum $D^{(j,0)}\oplus D^{(0,j)}$ as
subrepresentation.

So far it seems that everything is almost trivially analogous to the
undeformed case. It is not. If we consider irreducible representations
of mixed chirality, $D^{(i,j)}$, we find that the appearance of the
$\R$-matrix in $(g\otimes h)^*$ makes it impossible
to define conjugate spinors. It only works for $D^{(j,0)}$, because
$\rho^0 = \varepsilon$ and $\varepsilon(\R_{[1]})\R_{[2]} = 1$.
Fortunately, we do have conjugate spinors for the most interesting
cases: Dirac spinors ($D^{(\frac{1}{2},0)}\oplus D^{(0,\frac{1}{2})}$)
and the Maxwell tensor ($D^{(1,0)}\oplus D^{(0,1)}$). For these cases
everything is analogous to the undeformed case.

Let us consider a $D^{(j,0)}\oplus D^{(0,j)}$ representation with
basis $\{e^\mathrm{L}_m\}$ for the left chiral subrepresentation
$D^{(j,0)}$ and the basis $\{e^\mathrm{R}_m\}$ for $D^{(0,j)}$. We
define the conjugate basis by $\overline{e^\mathrm{L}_m} :=
e^\mathrm{R}_m$ and $\overline{e^\mathrm{R}_m} = e^\mathrm{L}_m$. Let
us call $\mathcal{P}$ the parity operator that exchanges the left and
right chiral part. Its matrix representation in the basis
$\{e^\mathrm{L}_m , e^\mathrm{R}_m\}$ is
\begin{equation}
  \mathcal{P}_{mn}=
  \begin{pmatrix} 0 & 1 \\ 1 & 0 \end{pmatrix}\,,
\end{equation}
where $1$ is the $(2j+1)$-dimensional unit matrix. This is the matrix
that represents our new pseudo scalar product as a bilinear form.  The
pseudo Hermitian conjugate of some operator $A$ can now be written as
\begin{equation}
  j(A) := \mathcal{P}A^\dagger \mathcal{P} \,,
\end{equation}
which is an involution because $\mathcal{P} = \mathcal{P}^\dagger$ and
an algebra anti-homomorphism because $\mathcal{P} = \mathcal{P}^{-1}$.

We apply this result to the whole space of spinor wave functions
$\mathcal{S}\otimes \Mink$. Let us assume that the scalar product of
two wave functions $f,g\in\Mink$ can be written (at least formally) as
some sort of integral $\lrAngle{f | g } = \int f^* g$.  The pseudo
scalar product of two $D^{(j,0)}\oplus D^{(0,j)}$ spinors $\psi$,
$\phi$ becomes
\begin{equation}
\begin{split}
  (\psi|\phi) &= (e_m\otimes \psi^m|e_n \otimes \phi^n)
  =(e_m|e_n) \lrAngle{ \psi^m|\phi^n} \\
  &= {\textstyle\int} (\psi^m)^* \mathcal{P}_{mn} \phi^n
  = {\textstyle\int} \bar{\psi}^n \phi^n \,,
\end{split}
\end{equation}
with the conjugate spinor wave function defined as
\begin{equation}
\label{eq:Wave3}
  \bar{\psi}^n := (\psi^m)^* \mathcal{P}_{mn} \,.
\end{equation}

To summarize, we have convinced ourselves that in the case of
$D^{(j,0)}\oplus D^{(0,j)}$ representations the conjugation of
spinors, spinor wave functions and operators works exactly as in the
undeformed case.

\section{The $q$-Dirac Equation}
\label{sec:MainContrib4b}

\subsection{The $q$-Dirac Equation in the Rest Frame}

In this section we consider $q$-Dirac spinors $\psi = e_j
\otimes\psi^j$ with the spin part transforming according to a
$D^{(\frac{1}{2},0)}\oplus D^{(0,\frac{1}{2})}$ representation. We
hope that we can write the projection onto an irreducible component of
the space of $q$-Dirac spinors as expression which involves momenta
only to first order terms, corresponding to a first order differential
equation. The general expression for such a $q$-Dirac equation would
be
\begin{equation}
  \Proj \psi := \frac{1}{2m}(m + \gamma_\mu P^\mu)\psi = 0 \,,
\end{equation}
with $\gamma_\mu$ being some operators acting on $\psi^j$. We can
already say that $\gamma_\mu$ must be a left 4-vector operator.  If it
were not, $\gamma_\mu P^\mu$ would not be scalar and, hence, would not
commute with the $q$-Lorentz transformations as required in
Eq.~\eqref{eq:WaveCondition2}.

We consider here a massive $q$-Dirac spinor representation, so there
is a rest frame (Sec.~\ref{sec:Little2}), that is, a set of states
$\psi^j$, which the momenta act upon as $P^0\psi^j = m \psi^j$,
$P^A\psi^j = 0$. We start the search for a projector $\Proj$ that
reduces the $q$-Dirac representation by computing how it has to act on
the rest frame, where we have
\begin{equation}
  \Proj_0 = \tfrac{1}{2}(1 + \gamma_0) \,,
\end{equation}
the zero indicating that this is a projector within the rest frame
only. We assume that we can realize the operator $\gamma_0$ as
4$\times$4-matrix that acts on the spin degrees of freedom only. This
is not unreasonable, for if $\gamma_\mu$ is a set of matrices that
form a $4$-vector operator in the $D^{(\frac{1}{2},0)}\oplus
D^{(0,\frac{1}{2})}$ representation then $\gamma_\mu \otimes 1$ will
also be a $4$-vector operator in the representation of spinor wave
functions. So, let us assume we can write $\Proj_0 = \Proj_0 \otimes
1$ in block form as
\begin{equation}
  \Proj_0 = \begin{pmatrix} A & B \\
    C & D \end{pmatrix} \,,
\end{equation}
where $A$, $B$, $C$, $D$ are 2$\times$2-matrices.

Recall that $\Proj_0$ must satisfy
condition~\eqref{eq:WaveCondition2}. This tells us in particular that
$\Proj_0$ must commute with rotations, the symmetry of the rest frame.
A rotation $l$ is represented by
\begin{equation}
  \rho(l) = \begin{pmatrix}
    \rho^{\frac{1}{2}}(l) & 0 \\ 0 & \rho^{\frac{1}{2}}(l)
  \end{pmatrix} \,.
\end{equation} 
Since the $\rho^{\frac{1}{2}}$-representations of the rotations
generate all 2$\times$2-matrices (the $q$-Pauli matrices are a basis),
$\Proj_0$ will only commute with all rotations if $A$, $B$, $C$, $D$
are numbers, that is, complex multiples of the unit matrix.

Furthermore, $\Proj_0$ has to be a projector, $\Proj_0^2 = \Proj_0$,
$\Proj_0^\dagger = \Proj_0$, and, as in the undeformed case, we
require it to commute with the parity operator,
$[\Proj_0,\mathcal{P}]=0$.  Altogether these conditions fix $\Proj_0$
and hence $\gamma_0$ uniquely to be
\begin{equation}
  \gamma_0 = \begin{pmatrix} 0 & 1 \\ 1 & 0 \end{pmatrix} \,,
\end{equation}
the same as in the undeformed case.

\subsection{The $q$-Gamma Matrices and the $q$-Clifford Algebra}
\label{sec:gammaboost}

If $\gamma_0$ is to be a 4-vector operator, we have to define the
other gamma matrices as in Eq.~\eqref{eq:Boost1} by 
\begin{equation}
\begin{aligned}
 \gamma_- &= \adL(-q^{-\frac{1}{2}}\lambda^{-1}[2]^{\frac{1}{2}} \,c)
             \tr \gamma_0 \\ 
 \gamma_+ &= \adL(q^{\frac{1}{2}}\lambda^{-1}[2]^{\frac{1}{2}} \,b)
             \tr \gamma_0 \\ 
 \gamma_3 &= \adL(\lambda^{-1}\,(d-a)) \tr \gamma_0 \,,
\end{aligned}
\end{equation}
where the adjoint action is understood with respect to the $q$-Dirac
representation. To compute this, explicitly, we have to calculate
the representations of the boosts first.
\begin{subequations}
\begin{xalignat}{2}
  \rho(a) &=  \begin{pmatrix}
    \rho^{\frac{1}{2}}(K^{\frac{1}{2}}) & 0 \\
    0 & \rho^{\frac{1}{2}}(K^{-\frac{1}{2}}) \end{pmatrix}, &
  \rho(b) &=  \begin{pmatrix} 0 & 0 \\
    0 & q^{-\frac{1}{2}}\lambda \rho^{\frac{1}{2}}(K^{-\frac{1}{2}}E)
  \end{pmatrix}\\
  \rho(c) &=  \begin{pmatrix} 
    -q^{\frac{1}{2}}\lambda \rho^{\frac{1}{2}}(FK^{\frac{1}{2}})
    & 0 \\ 0 & 0 \end{pmatrix}, &
  \rho(d) &=  \begin{pmatrix}
    \rho^{\frac{1}{2}}(K^{-\frac{1}{2}}) & 0 \\
    0 & \rho^{\frac{1}{2}}(K^{\frac{1}{2}}) \end{pmatrix}
\end{xalignat}
\end{subequations}
To demonstrate the simplicity of the technique of boosting, let us
demonstrate it with an example.
\begin{equation}
\begin{split}
  \gamma_+
  &= \adL(q^{\frac{1}{2}}\lambda^{-1}[2]^{\frac{1}{2}}\,b)\tr\gamma_0
  = q^{\frac{1}{2}}\lambda^{-1}[2]^{\frac{1}{2}}
    [\rho(b)\gamma_0\rho(a) - q \rho(a)\gamma_0\rho(b) ] \\
  &= [2]^{\frac{1}{2}}\left[
    \begin{pmatrix} 0 & 0 \\
      \rho^{\frac{1}{2}} ( K^{-\frac{1}{2}} E K^{\frac{1}{2}} ) & 0
    \end{pmatrix} 
    -q \begin{pmatrix}
      0 & \rho^{\frac{1}{2}}(E) \\ 0 & 0
    \end{pmatrix} \right] 
  = \begin{pmatrix} 0 & q \,\sigma_+ \\
    -q^{-1}\sigma_+ & 0 \end{pmatrix} 
\end{split}
\end{equation}
Here, $\sigma_+$ is the $q$-Pauli matrix
(Sec.~\ref{sec:FourvectorsPauli}). After doing the other calculations
we get
\begin{xalignat}{2}
  \gamma_0 &=  \begin{pmatrix} 0 & 1 \\ 1 & 0 \end{pmatrix}, &
  \gamma_A &=  \begin{pmatrix} 0 & q\, \sigma\!_A \\
      -q^{-1}\sigma\!_A & 0 \end{pmatrix} \,,
\end{xalignat}
where $A$ runs as usual through $\{-,+,3\}$.

This result can be easily generalized to higher spin massive
particles. All we have to do for a massive $D^{(j,0)}\oplus
D^{(0,j)}$-spinor is to replace $\rho^{\frac{1}{2}}$ with $\rho^j$ in
the above calculations. The result are higher dimensional
$\gamma$-matrices
\begin{xalignat}{2}
\label{eq:GammaGeneral}
  \gamma_0^{(j)} &= \begin{pmatrix} 0 & 1 \\ 1 & 0 \end{pmatrix}, &
  \gamma_A^{(j)} &=  [2]\begin{pmatrix} 0 & q\,\rho^j(J_A)  \\
    -q^{-1}\rho^j(J_A) & 0 \end{pmatrix} \,.
\end{xalignat}

Now we want to write the $q$-Dirac equation as $q$-differential
equation. Towards this end we need to calculate $\tilde{\gamma}_\mu$
by formula~\eqref{eq:Optilde}. Using Eqs~\eqref{eq:AppR2}
and~\eqref{eq:AppR1} we get for the $q$-Pauli matrices
\begin{xalignat}{2}
  \sigma\!_A\,\rho^{(\frac{1}{2},0)}
  \bigl((L^\Lambda_{\mathrm{I}+})^{A}{}_{B}\bigr)
  &= q^2\tilde{\sigma}\!_B \,,& 
  \sigma\!_A\,\rho^{(0,\frac{1}{2})}
  \bigl((L^\Lambda_{\mathrm{I}+})^{A}{}_{B}\bigr)
  &= q^{-2}\tilde{\sigma}_B \,,
\end{xalignat}
where
\begin{xalignat}{3}
  \tilde{\sigma}_- &= {[2]}^{\frac{1}{2}}
    \begin{pmatrix}
      0 & q^{\frac{1}{2}} \\ 0 & 0
    \end{pmatrix}, &
  \tilde{\sigma}_+ &= [2]^{\frac{1}{2}}
    \begin{pmatrix}
      0 & 0 \\ -q^{-\frac{1}{2}} & 0
    \end{pmatrix},&
  \tilde{\sigma}_3 &=
    \begin{pmatrix} -q^{-1} & 0 \\ 0 & q \end{pmatrix} 
\end{xalignat}
with respect to the $\{-,+\}$ basis. We can write this more compactly
as
\begin{equation}
  \tilde{\sigma}\!_A = -[2]\,\rho^\frac{1}{2}(SJ_A) \,.
\end{equation}
In this sense the transformed $q$-Pauli matrices,
$\tilde{\sigma}\!_A$, can be viewed as antipodes of the original ones.
For the transformed $q$-gamma matrices we obtain
\begin{xalignat}{2}
\label{eq:gamma2}
  \tilde{\gamma}_0 &=
    \begin{pmatrix} 0 & 1 \\ 1 & 0 \end{pmatrix}\,, &
  \tilde{\gamma}_A &=
    \begin{pmatrix} 0 & q^{-1}\, \tilde{\sigma}\!_A \\
      -q \tilde{\sigma}\!_A & 0 \end{pmatrix} \,, 
\end{xalignat}
so the $q$-Dirac equation written as $q$-differential equation becomes
\begin{equation}
  (m - \I \tilde{\gamma}_\mu \partial^\mu)\psi = 0 \,.
\end{equation}
What commutation relations do the gamma matrices satisfy? Using
Eqs.~\eqref{eq:gamma2} and~\eqref{eq:sigmatilderel} we find after some
lengthy calculations
\begin{equation}
\label{eq:Clifford}
  \tilde{\gamma}_c  \tilde{\gamma}_d
  = \eta_{dc} + \tilde{\gamma}_a \tilde{\gamma}_b
  \Proj_\mathrm{A}^{ba}{}_{dc} \,,
\end{equation}
where $\Proj_\mathrm{A}$ is the antisymmetric projector defined in
Eq.~\eqref{eq:AppClebsch1}.  This is the $q$-deformation of the
Clifford algebra. Using the relations of the $q$-Clifford algebra it
can be shown that the square of $q$-Dirac operator is indeed the mass
Casimir,
\begin{equation}
  (\tilde{\gamma}_\mu \partial^\mu)^2 = \partial_\mu \partial^\mu
  = -P_\mu P^\mu \,.
\end{equation}
As in the undeformed case we conclude that a solution $\psi$ to the
$q$-Dirac equation satisfies automatically the mass shell condition
$P_\mu P^\mu \psi = m^2 \psi$, and that the operator $\Proj =
\frac{1}{2m}(m+\gamma_\mu P^\mu)$ really is a projector. The
$q$-Clifford relations~\eqref{eq:Clifford} can be written in
equivalent but more familiar forms as
\begin{equation}
\label{eq:Clifford2}
  \tilde{\gamma}_a \tilde{\gamma}_b
  \Proj_\mathrm{S}^{ba}{}_{dc} = \eta_{dc} \,,
  \qquad\text{or}\qquad
  \tilde{\gamma}_c \tilde{\gamma}_d + \tilde{\gamma}_a \tilde{\gamma}_b
  R_{\mathrm{II}}^{ab}{}_{dc} = q[2]\eta_{cd} \,,
\end{equation}
with the symmetrizer~\eqref{eq:AppClebsch1} and the
$R$-matrix~\eqref{eq:AppR3}. 

One could have started directly from these relations trying to find
matrices that satisfy them \cite{Schirrmacher:1992}. This approach has
a number of disadvantages: a) It is computationally much more
cumbersome than boosting $\gamma_0$. b) The result is not unique, that
is, we would get many solutions to the $q$-Clifford algebra not
knowing which representations they belong to. c) Having determined a
solution $\tilde{\gamma}_\mu$, the covariance of the $q$-Dirac
equation remains unclear as $\tilde{\gamma}_\mu$ cannot be a 4-vector
operator.

\subsection{The Zero Mass Limit and the $q$-Weyl Equations}

The zero mass limit of the $q$-Dirac equations, $(m+\gamma_\mu
P^\mu)\psi = 0$, is formally
\begin{equation}
   \gamma_\mu P^\mu \psi = 0 \,,
\end{equation}
where $\gamma_\mu$ is defined as in Eq.~\eqref{eq:GammaGeneral}. The
operator $\mathbb{A} := \gamma_\mu P^\mu$ is no longer a projection.
For $m\rightarrow 0$ the wave equation decouples into two independent
equations for a left handed $D^{(\frac{1}{2},0)}$-spinor
$\psi_\mathrm{L}$ and a right handed $D^{(0,\frac{1}{2})}$-spinor
$\psi_\mathrm{R}$,
\begin{xalignat}{2}
  \sigma\!_A P^A \psi_\mathrm{L} &= q^{-1}P^0 \psi_\mathrm{L} \,,&
  \sigma\!_A P^A \psi_\mathrm{R} &= -q P^0 \psi_\mathrm{R} \,,
\end{xalignat}
the $q$-Weyl equations for massless left and right handed
spin-$\frac{1}{2}$ particles. Written as $q$-differential equation
they become
\begin{xalignat}{2}
  \tilde{\sigma}\!_A \partial^A \psi_\mathrm{L}
  &= -q \partial^0 \psi_\mathrm{L} \,,& 
  \tilde{\sigma}\!_A \partial^A \psi_\mathrm{R}
  &= q^{-1} \partial^0 \psi_\mathrm{R} \,.
\end{xalignat}
The operator $\mathbb{A}$ inherits property~\eqref{eq:WaveCondition1}
from the massive $q$-Dirac projector $\Proj$, so $\mathbb{A}\psi = 0$
is a viable wave equation. Let us see what it looks like in the
momentum eigenspace $\mathcal{H}_p$ for the momentum eigenvalues $p =
(p_0,p_-,p_+,p_3) = (k,0,0,k)$ for some real parameter $k$
(Sec.~\ref{sec:Little2}). On this subspace $\mathbb{A}$ acts as
\begin{equation} 
  \mathbb{A}\rvert_{\mathcal{H}_p}
  = k \begin{pmatrix}
    0 & 1-q\sigma_3 \\
    1+q^{-1}\sigma_3 & 0 \end{pmatrix}
  = k[2]\begin{pmatrix}
    0 & 0 & q & 0 \\
    0 & 0 & 0 & 0 \\
    0 & 0 & 0 & 0 \\
    0 & q^{-1} & 0 & 0 \end{pmatrix} \,.
\end{equation}
The kernel of this operator is 2-dimensional leaving us with two
states corresponding to helicity $\pm\frac{1}{2}$. 

If we generalize these considerations to higher spin Dirac type
spinors, we find that the corresponding operator $\mathbb{A}$ has a
zero kernel, $\ker\mathbb{A} = 0$, which can be easily verified in the
momentum eigenspace $\mathcal{H}_p$. In other words: the wave equation
for massive $D^{(j,0)}\oplus D^{(0,j)}$ spinor wave functions leads
for $m\rightarrow 0$ to a wave equation that has no solutions.  This
applies in particular to $q$-Maxwell spinors. Therefore, we need a
different approach to find the $q$-Maxwell equations.

\section{The $q$-Maxwell Equations}
\label{sec:MainContrib4c}

\subsection{The $q$-Maxwell Equations in the Momentum Eigenspaces}

In this section we consider massless spinors $\psi^j$ with the spinor
index transforming according to a $D^{(1,0)}\oplus D^{(0,1)}$
representation. According to the Clebsch-Gordan
series~\eqref{eq:CGSeries1} this type of spinor is equivalent to
considering an antisymmetric tensor $F^{\mu\nu}$ with two 4-vector
indices. These are the types of spinor wave functions commonly used to
describe the electromagnetic field, a massless field of spin-1.

We start our calculations in the massless momentum eigenspace
$\mathcal{H}_p$ with momentum eigenvalues $\mathrm{p} =
(p_0,p_-,p_+,p_3) = (k,0,0,k)$ for some real parameter $k$. In
Sec.~\ref{sec:Little2} we have shown this eigenspace to be invariant
under the little algebra $\mathcal{K}_0$, whose generators $K$, and
$N_A$ have been defined in Eq.~\eqref{eq:Little3}. Within
$\mathcal{H}_p$ the little algebra acts only on the spinor index. The
$D^{(1,0)}\oplus D^{(0,1)}$ matrix representation of the generators
are given by
\begin{equation}
\label{eq:Maxwell1}
\begin{gathered}
  N_- = -q[2]
    \begin{pmatrix} \rho^1(J_-) & 0 \\ 0 & 0 \end{pmatrix}, \quad
  N_+ = -q^{-1}[2]
    \begin{pmatrix} 0 & 0 \\ 0 & \rho^1(J_+)  \end{pmatrix}, \quad
  N_3 =
    \begin{pmatrix} 1 & 0 \\ 0 & 1 \end{pmatrix} \\
  K =
    \begin{pmatrix} \rho^1(K) & 0 \\ 0 & \rho^1(K) \end{pmatrix} \,,
\end{gathered}
\end{equation}
where $\rho^1$ is the vector representation map of $\suq$. 

We seek a projector $\Proj = \Proj \otimes 1$ that projects
onto an irreducible subrepresentation of the little algebra. We write
it in block form as
\begin{equation}
  \Proj = \begin{pmatrix} A & B \\
    C & D \end{pmatrix} \,,
\end{equation}
where $A$, $B$, $C$, $D$ are 3$\times$3-matrices. We must have
$\Proj_0^\dagger =\Proj_0$, so $A$ and $D$ must be Hermitian matrices
and $C=B^\dagger$. Recall from Eq.~\eqref{eq:Nrep} that within an
irreducible representation of $\mathcal{K}_0$ we have $N_\pm = 0$.
Therefore, we must have
\begin{equation}
  N_\pm \Proj = 0
\end{equation}
within the $D^{(1,0)}\oplus D^{(0,1)}$ spinor representation. This
leads to the conditions
\begin{xalignat}{4}
   \rho^1(J_-)\,A &= 0\,,  & \rho^1(J_+)\,D &= 0 \,, &
   \rho^1(J_-)\,B &= 0\,, & \rho^1(J_+)\,B^\dagger &= 0 \,.
\end{xalignat}
To satisfy these conditions $A$, $B$, and $D$ must be of the form
\begin{xalignat}{3}
  A &= \begin{pmatrix}
    \alpha & 0 & 0 \\
    0 & 0 & 0 \\
    0 & 0 & 0  \end{pmatrix}, &
  B &= \begin{pmatrix}
    0 & 0 & \beta \\
    0 & 0 & 0 \\
    0 & 0 & 0  \end{pmatrix}, &
  D &= \begin{pmatrix}
    0 & 0 & 0 \\
    0 & 0 & 0 \\
    0 & 0 & \delta  \end{pmatrix},
\end{xalignat}
for $\alpha$, $\delta$ real and $\beta$ complex. Furthermore,
$\Proj$ must project on an eigenvector of $K$. From this it
follows that $\beta = 0$ and either $\alpha=1$, $\delta=0$ or
$\alpha=0$, $\delta=1$. To summarize, there are two possible
projectors 
\begin{xalignat}{2}
  \Proj_\mathrm{L} &= \begin{pmatrix}
    1 &   &   &   \\
      & 0 &   &   \\
      &   & \diagdown  &   \\
      &   &   & 0 \end{pmatrix}, &
  \Proj_\mathrm{R} &= \begin{pmatrix}
    0 &   &   &   \\
      & \diagdown &   &   \\
      &   & 0 &   \\
      &   &   & 1 \end{pmatrix} 
\end{xalignat}
projecting each on a irreducible one-dimensional representation of the
little algebra $\mathcal{K}_0$. The image of $\Proj_\mathrm{L}$
is part of the left handed $D^{(1,0)}$ component while
$\Proj_\mathrm{R}$ projects to the right handed $D^{(0,1)}$
component of the spinor. Physically, this corresponds to left and
right handed circular waves. We want to allow for parity
transformations exchanging the left and right handed parts, so we need
both parts
\begin{equation}
  \Proj = \Proj_\mathrm{L} + \Proj_\mathrm{R} \,.
\end{equation}
With the parity transformation included, the two dimensional space
which $\Proj$ projects onto is irreducible.

\subsection{Computing the $q$-Maxwell Equation}

We would like to find the $q$-Maxwell equation in the form of a first
order differential equation
\begin{equation}
\label{eq:Maxwell2}
  \mathbb{A}\psi = C_\mu P^\mu \, \psi = 0 \,,
\end{equation}
hoping that the operators $C_\mu$ can be chosen to act on the spinor
index only, $C_\mu = C_\mu \otimes 1$. This wave equation has to
fulfill condition~\eqref{eq:WaveCondition1}: The $q$-Lorentz transform
of a solution must again be a solution. For this, it would be
sufficient but not necessary, if $\mathbb{A}$ were a scalar operator,
as it has been the case for the $q$-Dirac equation and its zero mass
limit, the $q$-Weyl equations.

Recall from the last section, that as long as we do not include parity
transformations, we must have two independent equations for the right
and the left handed part of the spinor, $\psi_\mathrm{L}$ carrying a
$D^{(1,0)}$ representation and $\psi_\mathrm{R}$ carrying a
$D^{(0,1)}$ representation
\begin{xalignat}{2}
  \mathbb{A}_\mathrm{L}\psi_\mathrm{L} &= 0 \,,&
  \mathbb{A}_\mathrm{R}\psi_\mathrm{R} &= 0 \,.
\end{xalignat}
Let us try to choose $\mathbb{A}_\mathrm{L}= C_\mu^\mathrm{L} P^\mu$
and $\mathbb{A}_\mathrm{R} = C_\mu^\mathrm{R} P^\mu$, so they commute
with rotations. For this to be possible $C_0^\mathrm{L}$,
$C_0^\mathrm{R}$ must be scalars with respect to rotations while
$C_A^\mathrm{L}$, $C_A^\mathrm{R}$ must transform as 3-vectors. The
only scalar operators within the $D^1$-representation of rotations are
multiples of the unit matrix, while every 3-vector operator is
proportional to $\rho^1(J_A)$. Hence, up to an overall constant factor
our wave equations can be written as
\begin{xalignat}{2}
  \bigl(P^0 + \alpha_\mathrm{L}\,\rho^1(J_A)P^A\bigr)\psi_\mathrm{L} &= 0\,,&
  \bigl(P^0 + \alpha_\mathrm{R}\,\rho^1(J_A)P^A\bigr)\psi_\mathrm{R} &= 0 \,,
\end{xalignat}
where $\alpha_\mathrm{L}$, $\alpha_\mathrm{R}$ are constants. To
determine these constants, we consider the wave equations in the
momentum eigenspace, where they take the form
\begin{xalignat}{2}
  \bigl(1 + \alpha_\mathrm{L}\,\rho^1(J_3)\bigr)\psi_\mathrm{L} &= 0 \,,&
  \bigl(1 + \alpha_\mathrm{R}\,\rho^1(J_3)\bigr)\psi_\mathrm{R} &= 0 \,.
\end{xalignat}
The space of solutions of each of these equations must equal the image
of the projectors $\Proj_\mathrm{L}$ and $\Proj_\mathrm{R}$,
respectively. This requirement fixes the constants to
$\alpha_\mathrm{L} = q^{-1}$ and $\alpha_\mathrm{R}= -q$.

Although this determines our candidate for the $q$-Maxwell equation,
condition~\eqref{eq:WaveCondition1} has yet to be checked for the
boosts.  Let $\psi_0\in\mathcal{H}_p$ be an element of the momentum
eigenspace, $P_\mu \psi_0 = p_\mu\psi_0$, with $p_\mu =
(p_0,p_-,p_+,p_3) = (k,0,0,k)$. Using the commutation relations
between boosts and momentum generators we find
\begin{subequations}
\begin{xalignat}{2}
  P_\mu (a \psi_0) &= q^{-1} p_\mu(a\psi_0) \,,&
  P_\mu (b \psi_0) &= q^{-1} p_\mu(b\psi_0) \\
  P_\mu (c \psi_0) &= q p_\mu(c\psi_0)  \,,&
  P_\mu (d \psi_0) &= q p_\mu(d\psi_0)\,.
\end{xalignat}
\end{subequations}
By induction it follows, that for any monomial in the boosts, $h=a^i
b^j c^k d^l$, we have $P_\mu (h \psi_0) = q^{k+l-i-j}\,
p_\mu(h\psi_0)$. Thus, for any such $\psi:=h\psi_0$, the wave
equation~\eqref{eq:Maxwell2} takes the form
\begin{equation}
\label{eq:Maxwell3}
  (C_0 - C_3)\psi = 0 \,.
\end{equation} 
Looking separately at the left and right handed part of $\psi =
\psi_\mathrm{L} + \psi_\mathrm{R}$ this equation writes out
\begin{xalignat}{2}
  \begin{pmatrix} 0&0&0 \\0&q^{-2}&0\\0&0&q^{-1}[2]\end{pmatrix}
  \begin{pmatrix} \psi_\mathrm{L}^- \\
                  \psi_\mathrm{L}^3 \\
                  \psi_\mathrm{L}^+ \end{pmatrix} &= 0\,, &
  \begin{pmatrix} q[2]&0&0 \\0&q^{2}&0\\0&0&0\end{pmatrix}
  \begin{pmatrix} \psi_\mathrm{R}^- \\
                  \psi_\mathrm{R}^3 \\
                  \psi_\mathrm{R}^+ \end{pmatrix} &= 0 \,,
\end{xalignat}
which is equivalent to
\begin{xalignat}{2}
  \psi_\mathrm{L}^3 = \psi_\mathrm{L}^+ &= 0 \,,&
  \psi_\mathrm{R}^- = \psi_\mathrm{R}^3 &= 0 \,.
\end{xalignat}
If we now have a solution of Eq.~\eqref{eq:Maxwell3}, that is, a spinor
$\psi$ whose only non-vanishing components are $\psi_\mathrm{L}^-$
and $\psi_\mathrm{R}^+$, could it happen that by boosting it gets
other non-vanishing components, thus turning a solution into a
non-solution? The answer to this question is no. We exemplify this,
applying formula~\eqref{eq:Wave1} for the action of the boost
generator $c$ on a left handed spinor, 
\begin{equation}
\begin{split}
  c\, \psi_\mathrm{L}^A
  &= \rho^{(1,0)}(c_{(1)})^A{}_{A'}
     \bigl(c_{(2)}\tr \psi_\mathrm{L}^{A'}\bigr) \\
  &= \rho^{(1,0)}(c)^A{}_{A'} \bigl(a\tr \psi_\mathrm{L}^{A'}\bigr)
    +\rho^{(1,0)}(d)^A{}_{A'} \bigl(c\tr \psi_\mathrm{L}^{A'}\bigr) \\
  &= -q^{\frac{1}{2}}\lambda \rho^{1}(FK^{\frac{1}{2}})^A{}_{A'}
      \bigl(a\tr \psi_\mathrm{L}^{A'}\bigr)
    +\rho^{1}(K^{-\frac{1}{2}})^A{}_{A'}
      \bigl(c\tr \psi_\mathrm{L}^{A'}\bigr) \\
  &= -q^{\frac{1}{2}}\lambda[2]^{\frac{1}{2}}
    \begin{pmatrix}0&1&0\\0&0&1\\0&0&0\end{pmatrix}
    \begin{pmatrix} a\tr \psi_\mathrm{L}^-\\
                    a\tr \psi_\mathrm{L}^3\\
                    a\tr \psi_\mathrm{L}^+ \end{pmatrix}
    +\begin{pmatrix}q&0&0\\0&1&0\\0&0&q^{-1}\end{pmatrix}
    \begin{pmatrix} c\tr \psi_\mathrm{L}^-\\
                    c\tr \psi_\mathrm{L}^3\\
                    c\tr \psi_\mathrm{L}^+ \end{pmatrix} \\
  &=\begin{pmatrix}
      -q^{\frac{1}{2}}\lambda[2]^{\frac{1}{2}} a\tr \psi_\mathrm{L}^3
      +q\,c\tr \psi_\mathrm{L}^- \\
      -q^{\frac{1}{2}}\lambda[2]^{\frac{1}{2}} a\tr \psi_\mathrm{L}^+
      +c\tr \psi_\mathrm{L}^3 \\
      q^{-1} c\tr \psi_\mathrm{L}^+
    \end{pmatrix},
\end{split}
\end{equation}
which clearly shows that, if $\psi_\mathrm{L}^3$ and
$\psi_\mathrm{L}^+$ vanish, so do $c\psi_\mathrm{L}^3$ and
$c\psi_\mathrm{L}^+$. Similar calculations can be done for the other
boost generators and right handed spinors.

By induction we conclude, that if $\psi_0\in\mathcal{H}_p$ is a
solution of Eq.~\eqref{eq:Maxwell3} and $h=a^i b^j c^k d^l$ is a
monomial in the boosts, $h=a^i b^j c^k d^l$, the spinor $\psi =
h\psi_0$ will be a solution, as well. The algebra of all boosts,
$\SUq^\op$, is generated as linear space by monomials, thus, $h\psi_0$
is a solution for any boost $h\in\SUq^\op$.  Since furthermore every
$q$-Lorentz transformation can be written as a sum of products of
rotations and boost, $h\psi_0$ is a solution for \emph{any}
$q$-Lorentz transformation $h$. We assume that the space of solutions,
$\ker\mathbb{A}$, is an irreducible representation. This means in
particular that the $q$-Lorentz algebra acts transitively on
$\ker\mathbb{A}$, so any solution can be written as $h\psi_0$.  Hence,
the wave equations
\begin{xalignat}{2}
\label{eq:Maxwell4}
  \rho^1(J_A)P^A\psi_\mathrm{L} &= -qP_0 \psi_\mathrm{L}\,, &
  \rho^1(J_A)P^A\psi_\mathrm{R} &= q^{-1}P_0 \psi_\mathrm{R} 
\end{xalignat}
do indeed satisfy property~\eqref{eq:WaveCondition1}. 

We want to write these equations, $C_\mu P^\mu \psi = 0$ as
$q$-differential equations $\tilde{C}_\mu \partial^\mu \psi = 0$,
where $\tilde{C}_\mu$ is defined in Eq.~\eqref{eq:Optilde}.  After
lengthy calculations using Eqs.~\eqref{eq:AppRep1}, \eqref{eq:AppR2},
and~\eqref{eq:AppR1} we get for the left and right handed part
separately
\begin{equation}
\begin{aligned}
  \rho^1(J_{A'})^B{}_{C'}\,\rho^{(1,0)}
  \bigl((L^\Lambda_+)^{A'}{}_{A}\bigr)^{C'}{}_C
  &= -q^2 \varepsilon_C{}^B{}_A  \\
  \rho^1(J_{A'})^B{}_{C'}\,\rho^{(0,1)}
  \bigl((L^\Lambda_+)^{A'}{}_{A}\bigr)^{C'}{}_C
  &= -q^{-2} \varepsilon_C{}^B{}_A \,,
\end{aligned}
\end{equation}
so the wave equations can be written as
\begin{xalignat}{2}
\label{eq:Maxwell5}
  \vec{\partial}\times \vec{\psi}_\mathrm{L}
  &= \I q^{-1}\partial_0 \vec{\psi}_\mathrm{L}\,, &
  \vec{\partial}\times \vec{\psi}_\mathrm{R}
  &= -\I q\,\partial_0 \vec{\psi}_\mathrm{R}\,,
\end{xalignat}
where $\vec{\psi}_\mathrm{R} = (\psi^A_\mathrm{R})$,
$\vec{\psi}_\mathrm{L} = (\psi^A_\mathrm{L})$, and where the cross
product is defined in Eq.~\eqref{eq:CrossDot}. A spinor
$\vec{\psi}_\mathrm{L}$ which is a solution to this equation must yet
satisfy the mass zero condition. Using the
identities~\eqref{eq:epsidentities2} for the cross product, the
commutation relations of the derivations $\vec{\partial} \times
\vec{\partial} = -\I\lambda\partial_0\vec{\partial}$, and the wave
equation~\eqref{eq:Maxwell5}, we rewrite the mass zero condition as
\begin{equation}
\label{eq:Butschi}
\begin{split}
  0 &= \partial_\mu \partial^\mu \vec{\psi}_\mathrm{L}
  = (\partial_0^2 - \vec{\partial}\cdot\vec{\partial})
    \vec{\psi}_\mathrm{L}\\
  &= \partial_0^2 \vec{\psi}_\mathrm{L}
    -(\vec{\partial}\times\vec{\partial})\times \vec{\psi}_\mathrm{L}
    +\vec{\partial}\times(\vec{\partial}\times \vec{\psi}_\mathrm{L})
    -\vec{\partial}(\vec{\partial}\cdot \vec{\psi}_\mathrm{L}) \\
  &=\partial_0^2 \vec{\psi}_\mathrm{L}
    +\I\lambda \partial_0(\vec{\partial} \times \vec{\psi}_\mathrm{L})
    +\vec{\partial}\times(\I q^{-1} \partial_0 \vec{\psi}_\mathrm{L})
    -\vec{\partial}(\vec{\partial}\cdot \vec{\psi}_\mathrm{L}) \\
  &=\partial_0^2 \vec{\psi}_\mathrm{L}
    -q^{-1}\lambda \partial^2_0 \vec{\psi}_\mathrm{L}
    - q^{-2} \partial^2_0 \vec{\psi}_\mathrm{L}
    -\vec{\partial}(\vec{\partial}\cdot \vec{\psi}_\mathrm{L}) \\
  &= -\vec{\partial}(\vec{\partial}\cdot \vec{\psi}_\mathrm{L}) \,. 
\end{split}
\end{equation}
Contracting the wave equation with $\vec{\partial}$ 
\begin{equation}
  \vec{\partial}\cdot(\vec{\partial}\times \vec{\psi}_\mathrm{L})
  =(\vec{\partial}\times\vec{\partial})\cdot \vec{\psi}_\mathrm{L}
  = -\I\lambda \partial_0 (\vec{\partial}\cdot \vec{\psi}_\mathrm{L})
  = \I q^{-1}\partial_0 (\vec{\partial}\cdot \vec{\psi}_\mathrm{L})\,,
\end{equation}
we see that $\partial_0 (\vec{\partial}\cdot \vec{\psi}_\mathrm{L}) =
0$ if $\vec{\psi}_\mathrm{L}$ is to satisfy the wave equation.
Together with Eq.~\eqref{eq:Butschi} this means that the mass zero
condition is equivalent to $\partial_\mu (\vec{\partial}\cdot
\vec{\psi}_\mathrm{L}) = 0$, that is, $\vec{\partial}\cdot
\vec{\psi}_\mathrm{L}$ must be a constant number. In a momentum
eigenspace we have $\partial_0 (\vec{\partial}\cdot
\vec{\psi}_\mathrm{L}) = k(\vec{\partial}\cdot
\vec{\psi}_\mathrm{L})$, so this constant number must be zero. The
same reasoning applies for the right handed spinor
$\vec{\psi}_\mathrm{R}$.

We conclude that the wave equations~\eqref{eq:Maxwell5} together with
the mass zero condition $\partial_\mu\partial^\mu \psi = 0$ are
equivalent to 
\begin{xalignat}{2}
\label{eq:Maxwell6}
  \vec{\partial}\times \vec{\psi}_\mathrm{L}
  &= \I q^{-1}\partial_0 \vec{\psi}_\mathrm{L}\,, &
  \vec{\partial}\cdot \vec{\psi}_\mathrm{L} &=0 \\
  \vec{\partial}\times \vec{\psi}_\mathrm{R}
  &= -\I q\,\partial_0 \vec{\psi}_\mathrm{R}\,, &
  \vec{\partial}\cdot \vec{\psi}_\mathrm{R} &=0 \,,
\end{xalignat}
which we will call the $q$-Maxwell equations.

\subsection{The $q$-Electromagnetic Field}

Finally, we write the $q$-Maxwell equations in a more familiar form,
that is, in terms of the $q$-deformed electric and magnetic fields.
In the undeformed case the electric and magnetic fields can --- up to
constant factors --- be characterized within the $D^{(1,0)}\oplus
D^{(0,1)}$ representation as eigenstates of the parity operator. The
electric field should transform like a polar vector $\mathcal{P}\vec{E} =
- \vec{E}$, while the magnetic field must be an axial vector
$\mathcal{P}\vec{B} = \vec{B}$. Recall, that the parity operator $\mathcal{P}$
acts on $q$-spinors by exchanging the left and the right handed parts
$\mathcal{P}\psi_\mathrm{L} = \psi_\mathrm{R}$,
$\mathcal{P}\psi_\mathrm{R} = \psi_\mathrm{L}$.  This fixes the fields
\begin{xalignat}{2}
  \vec{E} &= \I(\vec{\psi}_\mathrm{R} - \vec{\psi}_\mathrm{L}) \,,&
  \vec{B} &= \vec{\psi}_\mathrm{R} + \vec{\psi}_\mathrm{L}
\end{xalignat}
up to constant factors which have been chosen to give the right
undeformed limit. Spinor conjugation of the fields is now the same as
ordinary conjugation $\bar{E}^A = (E^A)^*$, $\bar{B}^A = (B^A)^*$.  In
terms of these fields, the $q$-Maxwell equations~\eqref{eq:Maxwell6}
take the form
\begin{xalignat}{2}
\label{eq:Maxwell7}
  \vec{\partial}\times \vec{E} &= \tfrac{1}{2}[2]\,\partial_0 \vec{B}
     -\tfrac{1}{2}\I\lambda\,\partial_0 \vec{E} \,, &
  \vec{\partial}\cdot \vec{E} &=0 \\
  \vec{\partial}\times \vec{B} &= -\tfrac{1}{2}[2]\,\partial_0 \vec{E}
     -\tfrac{1}{2} \I\lambda \,\partial_0 \vec{B} \,, &   
  \vec{\partial}\cdot \vec{B} &=0 \,.
\end{xalignat}
We would also like to express the $q$-Maxwell equations in terms of a
field strength tensor $F^{\mu\nu}$. According to the Clebsch-Gordan
series the left and right chiral $3$-vectors $\psi_\mathrm{L}$ and
$\psi_\mathrm{R}$ can be each identified with a 4-vector matrix. If we
replace in $\psi_\mathrm{L} = e_C \otimes \psi_\mathrm{L}^C$ the
spinor basis $e_C$ with $\Bas{1,0}{C,0}$ from
formula~\eqref{eq:chiralbases},
\begin{equation}
\begin{split}
  \psi_\mathrm{L} &= e_C \otimes  \psi_\mathrm{L}^C 
  = (E_A\otimes E_B \,\varepsilon^{AB}{}_C
  +q E_0\otimes E_C - q^{-1} E_C\otimes E_0) \otimes \psi_\mathrm{L}^C \\
  &= (E_\mu\otimes E_\nu)\otimes F_\mathrm{L}^{\mu\nu} \,,
\end{split}
\end{equation}
this defines the matrix
\begin{equation}
  F_\mathrm{L}^{\mu\nu} := \begin{pmatrix}
    F_\mathrm{L}^{00} & F_\mathrm{L}^{0N} \\
    F_\mathrm{L}^{M0} & F_\mathrm{L}^{MN}
  \end{pmatrix}
  = \begin{pmatrix}
    0 & q \psi_\mathrm{L}^N \\
    -q^{-1} \psi_\mathrm{L}^M &
    \varepsilon^{MN}{}_C \, \psi_\mathrm{L}^C
  \end{pmatrix},
\end{equation}
where $M$, $N$ run through $\{-,+,3\}$. In the same manner we obtain
for the right handed part
\begin{equation}
  F_\mathrm{R}^{\mu\nu} := \begin{pmatrix}
    0 & -q^{-1} \psi_\mathrm{R}^N \\
    q \psi_\mathrm{R}^M & \varepsilon^{MN}{}_C \, \psi_\mathrm{R}^C
  \end{pmatrix}.
\end{equation}
In terms of these matrices the $q$-Maxwell
equations~\eqref{eq:Maxwell6} take the form $\partial_\nu
F_\mathrm{L}^{\mu\nu} = 0$ and $\partial_\nu F_\mathrm{R}^{\mu\nu} =
0$. This suggests to introduce the field strength tensor and its dual
\begin{xalignat}{2}
  F^{\mu\nu}
  &:= \I(F_\mathrm{L}^{\mu\nu} + F_\mathrm{R}^{\mu\nu}) \,, &
  \tilde{F}^{\mu\nu}
  &:= \I(F_\mathrm{L}^{\mu\nu} - F_\mathrm{R}^{\mu\nu}) \,, 
\end{xalignat}
where the factor $\I$ is needed for the right undeformed limit. In
terms of the electric and the magnetic field, we have
\begin{equation}
\begin{aligned}
  F^{\mu\nu} &:= \begin{pmatrix}
    0 & -\tfrac{1}{2}[2]E^N +\tfrac{1}{2}\I\lambda B^N \\
    \tfrac{1}{2}[2]E^M +\tfrac{1}{2}\I\lambda B^M &
    \I\varepsilon^{MN}{}_C \, B^C
  \end{pmatrix} \\
  \tilde{F}^{\mu\nu} &:= \begin{pmatrix}
    0 & \tfrac{1}{2}[2]\I B^N -\tfrac{1}{2}\lambda E^N \\
    -\tfrac{1}{2}[2]\I B^M -\tfrac{1}{2}\lambda E^M &
    -\varepsilon^{MN}{}_C \, E^C
  \end{pmatrix}.
\end{aligned}
\end{equation}
The $q$-Maxwell equations become
\begin{xalignat}{2}
  \partial_\nu F^{\mu\nu} &= 0 \,, &
  \partial_\nu \tilde{F}^{\mu\nu} &= 0 \,,
\end{xalignat}
in complete analogy to the undeformed case.

\appendix

\chapter{Useful Formulas}
\label{sec:Appendix}

\section{Clebsch-Gordan Coefficients}
\label{sec:AppClebsch}

\subsection{Clebsch-Gordan and Racah Coefficients for $\suq$}
\label{sec:AppClebsch1}

We first list some formulas which enable us to calculate some
$q$-Clebsch-Gordan coefficients explicitly
\cite{BiedenharnLohe,Schmuedgen}:
\begin{equation}
\begin{aligned}
  \CGC{0}{j}{j'}{0}{m}{m'} &= \delta_{mm'}\delta_{jj'} \\
  \CGC{j}{\tfrac{1}{2}}{j+\tfrac{1}{2}}%
      {m}{\pm\tfrac{1}{2}}{m\pm\tfrac{1}{2}}
  &= q^{\pm(j\mp m)/2}([j\pm m + \tfrac{1}{2}][2j+1]^{-1})^{\tfrac{1}{2}}\\
  \CGC{j}{\tfrac{1}{2}}{j-\tfrac{1}{2}}%
      {m}{\pm\tfrac{1}{2}}{m\pm\tfrac{1}{2}}
  &= \mp q^{\mp(j\pm m +1)/2}([j \mp m][2j+1]^{-1})^{\tfrac{1}{2}}
\end{aligned}
\end{equation}
For $\CGC{1}{j}{j+\Delta j}{\Delta m}{m}{m+\Delta m}$ we have the
formulas
\\[4ex]
{ \setlength{\extrarowheight}{4pt} \renewcommand{\arraystretch}{1.5}
\begin{tabular}{>{$}r<{$}|>{$}c<{$}>{$}c<{$}>{$}c<{$}}
 \Delta j &\Delta m= -1 & \Delta m= 0  & \Delta m= +1 \\ \hline
 -1 & 
q^{m-j-1}\sqrt{\frac{[j+m][j+m-1]}{[2j+1][2j]}}& 
 -q^{m}\sqrt{\frac{[2][j+m][j-m]}{[2j+1][2j]}}&
q^{m+j+1}\sqrt{\frac{[j-m][j-m-1]}{[2j+1][2j]}}  \\ 
 0  &
-q^{m-1}\sqrt{\frac{[2][j+m][j-m+1]}{[2j+2][2j]}} &
q^{m}\,\frac{q^{(j+1)}[j-m]- q^{-(j+1)}[j+m]}{\sqrt{[2j+2][2j]}} &
q^{m+1}\sqrt{\frac{[2][j+m+1][j-m]}{[2j+2][2j]}} \\ 
 +1 & 
q^{m+j}\sqrt{\frac{[j-m+2][j-m+1]}{[2j+2][2j+1]}} &
q^{m}\sqrt{\frac{[2][j+m+1][j-m+1]}{[2j+2][2j+1]}} &
q^{m-j}\sqrt{\frac{[j+m+2][j+m+1]}{[2j+2][2j+1]}} \\ 
\end{tabular} }
\\[4ex]
The $q$-Clebsch-Gordan coefficients obey the symmety
\begin{equation}
\label{eq:ClebschSymmetrie}
  \CGC{n}{j}{j'}{\nu}{m}{m'} = (-1)^{j'-j}(-q)^{\nu}
  \sqrt{\frac{[2j'+1]}{[2j+1]}}\CGC{n}{j'}{j}{-\nu}{m'}{m} \,.
\end{equation}
For the $q$-Racah coefficients we have
\begin{equation}
  -\sqrt{[3]} \, \RC{1}{1}{j}{0}{j'}{j'} =
  (-1)^{j'+j} \sqrt{\frac{[2j'+1]}{[2j+1]}} \,.
\end{equation}
For $-\sqrt{\frac{[4]}{[2]}} \, \RC{1}{1}{j}{1}{j'}{j''}$ there are
the formulas
\\[4ex]
\begin{center}
\setlength{\extrarowheight}{4pt}
\renewcommand{\arraystretch}{1.5}
\begin{tabular}{>{$}l<{$}|>{$}c<{$}>{$}c<{$}>{$}c<{$}}
  &j''= j-1 & j''= j  & j'' = j+1 \\ \hline
j' = j-1 &  \sqrt{\frac{[2j-2]}{[2j]}}&
\sqrt{\frac{[2j+2][2j-1]}{[2j+1][2j]}}& 0 \\ 
j'= j  & -\sqrt{\frac{[2j+2]}{[2j]}} &
\frac{[2j]- [2j+2]}{\sqrt{[2j+2][2j]}} &
\sqrt{\frac{[2j]}{[2j+2]}} \\ 
j'= j+1 & 0 & -\sqrt{\frac{[2j+3][2j]}{[2j+2][2j+1]}} &
-\sqrt{\frac{[2j+4]}{[2j+2]}} \\ 
\end{tabular}
\end{center}

\subsection{Metric and Epsilon Tensor}
\label{sec:AppClebsch2}

We define the 3-metric and the epsilon tensor as\footnote{Metric and
  epsilon tensor, $g_{AB}$ and $\varepsilon^{AB}{}_{C}$, as defined
here correspond to $g^{AB}$ and $q\varepsilon_{BA}{}^{C}$ in
\cite{Lorek:1997a}.}
\begin{xalignat}{2}
  g^{AB} &:= -\sqrt{[3]} \CGC{1}{1}{0}{A}{B}{0}\,, &
  \varepsilon^{AB}{}_{C} &=
    -\sqrt{\frac{[4]}{[2]}} \CGC{1}{1}{1}{A}{B}{C} \,,
\end{xalignat}
where the capital roman indices run through $\{-1,0,1\} = \{-,3,+\}$.
The positions of the indices are chosen such that the basis vectors are
written with lower indices. From this definition it is clear that the
projectors on the subspaces on the right hand side of the
Clebsch-Gordan series
\begin{equation}
  D^1\otimes D^1 \cong D^0 \oplus D^1 \oplus D^3,
\end{equation}
that we denote by $\Proj_0$, $\Proj_1$, $\Proj_3$, can
be written as
\begin{equation}
\begin{aligned}
  \Proj_0^{AB}{}_{CD} &= [3]^{-1}\,g^{AB}g_{CD} \\
  \Proj_1^{AB}{}_{CD} &= [2][4]^{-1}\, \varepsilon^{ABX}\varepsilon_{DCX} \\
  \Proj_3^{AB}{}_{CD} &= \delta^A_C\delta^B_D
  - \Proj_0^{AB}{}_{CD} - \Proj_1^{AB}{}_{CD} \,,
\end{aligned}
\end{equation}
where the projectors act on lower indices, $\Proj\tr E_C E_D := E_A
E_B \Proj^{AB}{}_{CD}$. The nonzero values of the metric are
\begin{xalignat}{3}
  g^{-+} &= -q^{-1}\,, & g^{+-} &= -q\,, &  g^{33} &= 1 \,. 
\end{xalignat}
By definition $g_{AB}$ is the inverse of $g^{AB}$
\begin{equation}
  g_{AB} g^{BC} = \delta^C_A  = g^{CB} g_{BA}\,,
\end{equation}
implying $g_{AB} = g^{AB}$. This means that we can not raise and lower
the indices of the metric as usual. Instead, we get
\begin{equation}
  g_{AA'} g_{BB'} g^{A'B'} = g_{BA}\,.
\end{equation}
The nonzero values of the epsilon tensor are
\begin{subequations}
\begin{xalignat}{3}
  \varepsilon^{-3}{}_- &= q^{-1} & \varepsilon^{3-}{}_- &= -q &&\\
  \varepsilon^{-+}{}_3 &= 1 & \varepsilon^{+-}{}_3 &= -1 &
  \varepsilon^{33}{}_3 &= -\lambda \\
  \varepsilon^{3+}{}_+ &= q^{-1} & \varepsilon^{+3}{}_+ &= -q \,.&&
\end{xalignat}
\end{subequations}
Lowering the first index as usual by $\varepsilon_A{}^B{}_C :=
g_{AA'}\,\varepsilon^{A'B}{}_{C}$ we get
\begin{subequations}
\begin{xalignat}{3}
  \varepsilon_+{}^3{}_- &= -1 & \varepsilon_3{}^-{}_- &= -q &&\\
  \varepsilon_+{}^+{}_3 &= -q & \varepsilon_-{}^-{}_3 &= q^{-1} &
  \varepsilon_3{}^3{}_3 &= -\lambda \\
  \varepsilon_3{}^+{}_+ &= q^{-1} & \varepsilon_-{}^3{}_+ &= 1 \,. &&
\end{xalignat}
\end{subequations}
Lowering the second index 
\begin{subequations}
\begin{xalignat}{3}
  \varepsilon^{-}{}_{3-} &= q^{-1} & \varepsilon^{3}{}_{+-} &= q^2 &&\\
  \varepsilon^{-}{}_{-3} &= -q^{-1} & \varepsilon^{+}{}_{+3} &= q &
  \varepsilon^{3}{}_{33} &= -\lambda \\
  \varepsilon^{3}{}_{-+} &= -q^{-2} & \varepsilon^{+}{}_{3+} &= -q \,. &&
\end{xalignat}
\end{subequations}
With all indices down $\varepsilon_{ABC}:=
g_{AA'}\varepsilon^{A'}{}_{BC}$
\begin{subequations}
\begin{xalignat}{3}
  \varepsilon_{+3-} &= -1 & \varepsilon_{3+-} &= q^2 && \\
  \varepsilon_{+-3} &= 1 & \varepsilon_{-+3} &= -1 &
  \varepsilon_{333} &= -\lambda \\
  \varepsilon_{3-+} &= -q^{-2} & \varepsilon_{-3+} &= 1 \,. &&
\end{xalignat}
\end{subequations}
Various contractions of $\varepsilon$-tensor and metric yield useful
identities
\begin{equation}
\label{eq:epsidentities}
\begin{gathered}
  \varepsilon^{AB'C}g_{B'\!B} = \varepsilon^{CA}{}_{B}\,, \quad
  \varepsilon^{AB}{}_{C'}g^{C'\!C} = \varepsilon^{BCA} \\
  \varepsilon^{A'B'C} g_{A'\!A} g_{B'\!B} = \varepsilon_B{}^C{}_A \,,
  \quad
  \varepsilon_{A'BC'} g^{A'\!A} g^{C'\!C} = \varepsilon^{CA}{}_B \\
  g_{AB}\varepsilon^{ABC} = 0\,, \quad
  g_{CA}\varepsilon^{ABC} = 0\,, \quad
  \varepsilon_{ABC}\, g^{BA} = 0 \,,\quad
  \varepsilon_{ABC}\, g^{AC} = 0\\
  \varepsilon^{AXB}\varepsilon_{CXD}
  = \varepsilon^{BA}{}_X \varepsilon_C{}^X{}_D
  = \varepsilon^{BAX}\varepsilon_{DCX} \\
  \varepsilon_A{}^B{}_C \,\varepsilon^{AC}{}_D
  = [4][2]^{-1} \delta_D^B\,,\quad
  \varepsilon_{ABC} \,\varepsilon^{BAD}
  =\varepsilon_{BCA} \,\varepsilon^{ADB}
  = [4][2]^{-1} \delta^D_C\\
  \varepsilon^{AB}{}_X\varepsilon^{XC}{}_D + g^{AB}\delta^C_D
  = \varepsilon^{AX}{}_D\varepsilon^{BC}{}_X + \delta^A_D g^{BC}\,.
\end{gathered}
\end{equation}
There are relations between $\varepsilon$-tensors with the same index
in an upper and a lower position
\begin{equation}
\begin{gathered}
  \varepsilon_{ABC} = \varepsilon^{ACB}\,, \quad
  \varepsilon_A{}^B{}_C = \varepsilon^{AC}{}_B \,,\quad
  \varepsilon^A{}_{BC} = \varepsilon_A{}^{CB} \,.
\end{gathered}
\end{equation}
With the metric and the epsilon tensor we can define a scalar and a
vector product. Note, that if we defined real coordinates by $X_1 :=
\I(X_+ - X_+^*)$, $X_2 := X_+ + X_+^*$ we would get, e.g.,
$\varepsilon^{123} = -q\I$.  In the limit $q\rightarrow 1$ our epsilon
tensor will tend to $-\I$ times the usual undeformed epsilon tensor.
We therefore define for $3$-vector operators $X_A$ and $Y_B$
\begin{xalignat}{2}
\label{eq:epsidentities3}
  \vec{X}\cdot \vec{Y} &:= g^{AB} X_A Y_B \,,&
  (\vec{X}\times \vec{Y})_C &:= \I\,X_A Y_B\varepsilon^{AB}{}_C \,,
\end{xalignat}
where we use arrows to indicate the $3$-vector operators. Raising and
lowering the indices we get
\begin{xalignat}{2}
  \vec{X}\cdot \vec{Y} &:= g_{BA}\, X^A Y^B \,,&
  (\vec{X}\times \vec{Y})^C &:= \I\,X^A Y^B\varepsilon_B{}^C{}_A \,.
\end{xalignat}
With this notation some of the identities~\eqref{eq:epsidentities}
take on a very intuitive form
\begin{equation}
\label{eq:epsidentities2}
\begin{gathered}
  \vec{X} \cdot(\vec{Y} \times \vec{Z})
  = (\vec{X} \times \vec{Y}) \cdot \vec{Z} \\
  (\vec{X} \times \vec{Y}) \times \vec{Z}
  -(\vec{X}\cdot \vec{Y}) \vec{Z}
  = \vec{X} \times (\vec{Y} \times \vec{Z})
  -\vec{X} (\vec{Y} \cdot \vec{Z}) \,, 
\end{gathered}
\end{equation}
from which more relations can be deduced very easily.
Finally, we apply Eq.~\eqref{eq:RacahRed} to the scalar and
the vector product
\begin{subequations}
\label{eq:AppRedrel}
\begin{equation}
  \rBraket{j}{\vec{X} \cdot \vec{Y}}{j}
  =\sum_{j'}(-1)^{j'+j} \sqrt{\frac{[2j'+1]}{[2j+1]}}\,
  \rBraket{j}{\vec{X}}{j'} \rBraket{j'}{\vec{Y}}{j}
\end{equation}
\begin{equation}
\begin{split}
  \rBraket{j-1}{\vec{X} \times \vec{Y}}{j}
  = \I\sqrt{\frac{[2j-2]}{[2j]}}\,
    &\rBraket{j-1}{\vec{X}}{j-1} \rBraket{j-1}{\vec{Y}}{j}\\
    -\I\sqrt{\frac{[2j+2]}{[2j]}}\,
    &\rBraket{j-1}{\vec{X}}{j} \rBraket{j}{\vec{Y}}{j}
\end{split}
\end{equation}
\begin{equation}
\begin{split}
  \rBraket{j}{\vec{X} \times \vec{Y}}{j}
  = \I\sqrt{\frac{[2j+2][2j-1]}{[2j+1][2j]}}\,
    &\rBraket{j}{\vec{X}}{j-1} \rBraket{j-1}{\vec{Y}}{j}\\
    +\I\frac{[2j]- [2j+2]}{\sqrt{[2j+2][2j]}}\,
    &\rBraket{j}{\vec{X}}{j} \rBraket{j}{\vec{Y}}{j}\\
    -\I\sqrt{\frac{[2j+3][2j]}{[2j+2][2j+1]}}\,
    &\rBraket{j}{\vec{X}}{j+1} \rBraket{j+1}{\vec{Y}}{j}
\end{split}
\end{equation}
\begin{equation}
\begin{split}
  \rBraket{j+1}{\vec{X} \times \vec{Y}}{j}
  = \I\sqrt{\frac{[2j]}{[2j+2]}}\,
    &\rBraket{j+1}{\vec{X}}{j} \rBraket{j}{\vec{Y}}{j}\\
    -\I\sqrt{\frac{[2j+4]}{[2j+2]}}\,
    &\rBraket{j+1}{\vec{X}}{j+1} \rBraket{j}{\vec{Y}}{j} \,.
\end{split}
\end{equation}
If furthermore there is a $*$-structure $X_A^* = Y^A$, this implies
for the reduced matrix elements of a $*$-representation
\begin{equation}
\label{eq:RedKonjugation}
  \rBraket{j'}{\vec{X}}{j}
  = (-1)^{j'-j} \sqrt{\frac{[2j+1]}{[2j'+1]}}\,
  \overline{\rBraket{j}{\vec{Y}}{j'}} \,.
\end{equation}
\end{subequations}

\subsection{Clebsch-Gordan Coefficients for the $q$-Lorentz Algebra}
\label{sec:AppClebsch3}

The Clebsch-Gordan Coefficients for the $q$-Lorentz algebra can be
read off the formula for the basis vectors of the irreducible
subrepresentations
\begin{multline}
  \Ket{(k_1,k_2),(n_1,n_2)} =
  \sum
  \CGC{j_1}{j'_1}{k_1}{m_1}{b}{n_1} 
  \CGC{j_2}{j'_2}{k_2}{a}{m'_2}{n_2}\\
  \times (R^{-1})^{m_2 m'_1}_{ab}\,  
  \Ket{(j_1,j_2),(m_1,m_2)} \otimes
  \Ket{(j'_1,j'_2),(m'_1,m'_2)} \,,
\end{multline}
As the $R$-matrix is in general not unitary, these basis vectors have
yet to be normalized. We are in particular interested in the
$q$-Clebsch-Gordan coefficients for the decomposition of a tensor
product of two vector representations
\begin{equation}
  D^{(\frac{1}{2},\frac{1}{2})} \otimes D^{(\frac{1}{2},\frac{1}{2})}
  \cong
  D^{(0,0)} \oplus D^{(1,0)} \oplus D^{(0,1)} \oplus D^{(1,1)}\,.
\end{equation}
For a more compact notation we write
\begin{xalignat}{2}
  \Bas{j_1j_2}{m_1m_2} &:=\Ket{(j_1,j_2),(m_1,m_2)} \,, &
  \Base{ab} &:=\Ket{(\tfrac{1}{2},\tfrac{1}{2}),(a,b)} \,,
\end{xalignat}
where $a$, $b$ run through $\{-\tfrac{1}{2},\tfrac{1}{2}\} = \{-,+\}$.
We get for the \emph{unnormalized} basis vectors of the $D^{(1,0)}$
subrepresentation
\begin{equation}
\label{eq:AppMinkbasis1}
\begin{aligned}
  \Bas{1,0}{-1,0} &= q \Base{-+}\otimes \Base{--}
  - q^{-1}\Base{--}\otimes \Base{-+} \\
  \Bas{1,0}{0,0} &= \Base{++}\otimes \Base{--}
  - \Base{--}\otimes \Base{++} + \lambda \Base{-+}\otimes \Base{-+}\\
  & \quad +\Base{-+}\otimes \Base{+-} - q^{-2}\Base{+-}\otimes\Base{-+}\\
  \Bas{1,0}{+1,0} &= \Base{++}\otimes \Base{+-}
  - \Base{+-}\otimes \Base{++} + \lambda \Base{++}\otimes \Base{-+}\,,
\end{aligned}
\end{equation}
for the $D^{(0,1)}$ subrepresentation
\begin{equation}
\label{eq:AppMinkbasis2}
\begin{aligned}
  \Bas{0,1}{0,-1} &= \Base{+-}\otimes \Base{-+}
  - \Base{--}\otimes \Base{+-} + \lambda \Base{-+}\otimes \Base{--} \\
  \Bas{0,1}{0,0} &= \Base{++}\otimes \Base{--}
  - \Base{--}\otimes \Base{++} + \lambda \Base{-+}\otimes \Base{-+} \\
    &\quad + \Base{+-}\otimes \Base{-+} - q^{-2}\Base{-+}\otimes \Base{+-}\\
  \Bas{0,1}{0,+1} &= q \Base{++}\otimes \Base{-+}
  -q^{-1}\Base{-+}\otimes \Base{++} \,,
\end{aligned}
\end{equation}
and for the $D^{(0,0)}$ subrepresentation
\begin{equation}
\label{eq:AppMinkbasis3}
\begin{aligned}
  \Bas{0,0}{0,0} &= q\Base{++}\otimes \Base{--}
  + q^{-1}\Base{--}\otimes \Base{++} -q^{-1}\Base{-+}\otimes \Base{+-} \\
    &\quad - q^{-1}\Base{+-}\otimes \Base{+-}
      -q^{-1}\lambda \Base{-+}\otimes \Base{-+} \,.
\end{aligned}
\end{equation}
Expressed in terms of a $4$-vector basis we find bases for the
$D^{(1,0)}$, $D^{(0,1)}$, and $D^{(0,0)}$ subrepresentations
\begin{equation}
\label{eq:chiralbases}
\begin{aligned}
  \Bas{1,0}{C,0} &= E_A\otimes E_B \,\varepsilon^{AB}{}_C
    +q E_0\otimes E_C - q^{-1} E_C\otimes E_0 \\
  \Bas{0,1}{0,C} &= E_A\otimes E_B \,\varepsilon^{AB}{}_C
    +q E_C\otimes E_0 - q^{-1} E_0\otimes E_C \\
  \Bas{0,0}{0,0} &= E_\mu\otimes E_\nu \,\eta^{\mu\nu} \,,  
\end{aligned}
\end{equation}
which are neither orthogonal nor normalized. The last equation defines
up to a constant factor the 4-vector metric $\eta^{\mu\nu}$ whose
non-zero values are
\begin{xalignat}{4}
  \eta^{00} &= 1\,, & \eta^{-+} &= q^{-1}\,, &
  \eta^{+-} &= q\,, & \eta^{33} &= -1 \,,  
\end{xalignat}
which means in particular that $\eta^{AB} = - g^{AB}$. Let us denote
the projectors on the subrepresentations of the
$D^{(\frac{1}{2},\frac{1}{2})}$ representation in an obvious
notation\footnote{Elsewhere \cite{Lorek:1997a} the same projectors
  have been denoted by $P_T$, $P_+$, $P_-$, $P_S$, in that order.} by
\begin{equation}
\label{eq:Proj1}
  1 = \Proj_{(0,0)} + \Proj_{(1,0)} + \Proj_{(0,1)} + \Proj_{(1,1)}\,.
\end{equation}
The projectors on the symmetric and antisymmetric part are denoted by
\begin{xalignat}{2}
\label{eq:AppClebsch1}
  \Proj_\mathrm{S} &:= \Proj_{(0,0)} + \Proj_{(1,1)} \,, &
  \Proj_\mathrm{A} &:= \Proj_{(1,0)} + \Proj_{(0,1)} \,.
\end{xalignat}
These projectors can be determined from the bases of the corresponding
spaces, which we just computed. Note however that $D^{(1,0)}$ and
$D^{(0,1)}$ are not mutually orthogonal, so we have to project on
$D^{(1,0)}$ along $D^{(0,1)}$ and vice versa. We obtain for the trace part
\begin{equation}
  (\Proj_{(0,0)})^{ab}{}_{cd} = [2]^{-2}\,\eta^{ab}\eta_{cd} \,,
\end{equation}
for the left chiral and right chiral part:
\\[2ex]
{ 
\renewcommand{\arraystretch}{1.5}
\begin{tabular}{>{$\!\!\!\!}r<{$}|>{$}c<{$}>{$}c<{$}>{$}c<{$\!}}
\multicolumn{4}{c}{ $[2]^2 (\Proj_{(1,0)})^{ab}{}_{cd}=$ }\\ 
& C0 & 0D  & CD \\ \hline
A0 & \delta^A_C & -q^{-2}\delta^A_D & -q^{-1} \varepsilon_C{}^A{}_D\\
0B & -q^2\delta^B_C & \delta^B_D & q \varepsilon_C{}^B{}_D \\
AB & -q \varepsilon^{AB}{}_C & q^{-1} \varepsilon^{AB}{}_D &
\varepsilon^{AB}{}_X\varepsilon_C{}^X{}_D
\end{tabular} }~
{ 
\renewcommand{\arraystretch}{1.5}
\begin{tabular}{>{$\!\!\!\!}r<{$}|>{$}c<{$}>{$}c<{$}>{$}c<{\!$}}
\multicolumn{4}{c}{ $[2]^2 (\Proj_{(0,1)})^{ab}{}_{cd}=$ }\\ 
& C0 & 0D  & CD \\ \hline
A0 & \delta^A_C & -q^{2}\delta^A_D & q \varepsilon_C{}^A{}_D\\
0B & -q^{-2}\delta^B_C & \delta^B_D & -q^{-1} \varepsilon_C{}^B{}_D \\
AB & q^{-1}\varepsilon^{AB}{}_C & -q \varepsilon^{AB}{}_D &
\varepsilon^{AB}{}_X\varepsilon_C{}^X{}_D
\end{tabular} }
\\[2ex]
For the anti-symmetrizer this yields
\begin{center}
\renewcommand{\arraystretch}{1.5}
\begin{tabular}{>{$\!\!\!\!}r<{$}|>{$}c<{$}>{$}c<{$}>{$}c<{\!$}}
\multicolumn{4}{c}{ $[2]^2 (\Proj_{\mathrm{A}})^{ab}{}_{cd}=$ }\\ 
& C0 & 0D  & CD \\ \hline
A0 & 2\delta^A_C & -[4][2]^{-1}\delta^A_D & \lambda \varepsilon_C{}^A{}_D\\
0B & -[4][2]^{-1}\delta^B_C & 2\delta^B_D & \lambda \varepsilon_C{}^B{}_D \\
AB & -\lambda\varepsilon^{AB}{}_C & -\lambda \varepsilon^{AB}{}_D &
2\varepsilon^{AB}{}_X\varepsilon_C{}^X{}_D
\end{tabular}
\\[2ex]
\end{center}
The traceless symmetric part is given by Eq.~\eqref{eq:Proj1}.

\section{Representations}
\label{sec:AppRep}

\subsection{Representations of $\suq$}
\label{sec:AppRep1}

The action of the generators within the $D^j$ representation of $\suq$
is given by
\begin{equation}
\begin{aligned}
  E\Ket{j,m} &= q^{(m+1)} \sqrt{[j+m+1][j-m]}\,\Ket{j,m+1}\\
  F\Ket{j,m} &= q^{-m} \sqrt{[j+m][j-m+1]}\,\Ket{j,m-1}\\
  K\Ket{j,m} &= q^{2m}\Ket{j,m} \,.
\end{aligned}
\end{equation}
For for the vectorial generators this means
\begin{equation}
\label{eq:AppJrep}
  J_A \Ket{j,m} = -[2]^{-1}\sqrt{[2j+2][2j]}\,
  \CGC{1}{j}{j}{A}{m}{m+A}\,\Ket{j,m+A} \,.
\end{equation}
The value of the Casimir $W$ within such a representation is given by
\begin{equation}
\label{eq:AppWrep}
  \rho^j(W) = [2]^{-1}\bigl(q^{(2j+1)}+q^{-(2j+1)}\bigr)\,.
\end{equation}
For $j=\frac{1}{2}$ the generators are represented by 
\begin{xalignat}{3}
  E &:=\begin{pmatrix} 0 & 0 \\ q^{\frac{1}{2}} & 0 \end{pmatrix}, &
  F &:=\begin{pmatrix} 0 & q^{-\frac{1}{2}} \\ 0 & 0 \end{pmatrix},&
  K &:=\begin{pmatrix} q^{-1} & 0 \\ 0 & q \end{pmatrix},
\end{xalignat}
with respect to the $\{-,+\}$ basis.  The representation of the vector
generators $J_A$ is proportional to the $q$-Pauli matrices
\begin{xalignat}{2}
  \sigma\!_A &= [2]\,\rho^{\frac{1}{2}}(J_A) \,,&
  \tilde{\sigma}\!_A &= -[2]\,\rho^{\frac{1}{2}}(SJ_A) \,,
\end{xalignat}
where
\begin{xalignat}{3}
\label{eq:qPauli2}
  \sigma_- &= {[2]}^{\frac{1}{2}}
    \begin{pmatrix} 0 & q^{-\frac{1}{2}} \\ 0 & 0 \end{pmatrix}, &
  \sigma_+ &= [2]^{\frac{1}{2}}
    \begin{pmatrix} 0 & 0 \\ -q^{\frac{1}{2}} & 0 \end{pmatrix},&
  \sigma_3 &=
    \begin{pmatrix} -q & 0 \\ 0 & q^{-1} \end{pmatrix} \\ 
  \tilde{\sigma}_- &= {[2]}^{\frac{1}{2}}
    \begin{pmatrix}
      0 & q^{\frac{1}{2}} \\ 0 & 0
    \end{pmatrix}, &
  \tilde{\sigma}_+ &= [2]^{\frac{1}{2}}
    \begin{pmatrix}
      0 & 0 \\ -q^{-\frac{1}{2}} & 0
    \end{pmatrix},&
  \tilde{\sigma}_3 &=
    \begin{pmatrix} -q^{-1} & 0 \\ 0 & q \end{pmatrix},
\end{xalignat}
with respect to the $\{-\frac{1}{2},\frac{1}{2}\} = \{-,+\}$ basis.
The $q$-Pauli matrices satisfy the relations
\begin{xalignat}{2}
  \sigma\!_A\, \sigma_B \,\varepsilon^{AB}{}_C
  &= [4][2]^{-1}\, \sigma_C \,, &
  \sigma\!_A \sigma_B
  &= g_{AB} + \sigma_C\, \varepsilon_A{}^C{}_B\\
  \tilde{\sigma}\!_A\, \tilde{\sigma}_B \,\varepsilon^{BA}{}_C
  &= -[4][2]^{-1}\, \tilde{\sigma}_C \,, &
  \tilde{\sigma}\!_A \tilde{\sigma}_B
  &= g_{BA} - \tilde{\sigma}_C\, \varepsilon_B{}^C{}_A
\label{eq:sigmatilderel},
\end{xalignat}
For the $j=1$ vector representations we get
\begin{xalignat}{3}
  E &:= [2]^{\frac{1}{2}}
  \begin{pmatrix} 0&0&0 \\ 1&0&0 \\ 0&q&0 \end{pmatrix}, &
  F &:= [2]^{\frac{1}{2}}
  \begin{pmatrix} 0&1&0 \\ 0&0&q^{-1} \\ 0&0&0 \end{pmatrix}, &
  K &:= 
  \begin{pmatrix} q^{-2}&0&0 \\ 0&1&0 \\ 0&0&q^2 \end{pmatrix} 
\end{xalignat}
\begin{xalignat}{3}
  J_- &:= 
  \begin{pmatrix} 0&q^{-1}&0 \\ 0&0&1 \\ 0&0&0 \end{pmatrix}, &
  J_+ &:= 
  \begin{pmatrix} 0&0&0 \\ -1&0&0 \\ 0&-q&0 \end{pmatrix}, &
  J_3 &:= 
  \begin{pmatrix} -q&0&0 \\ 0&-\lambda&0 \\ 0&0&q^{-1} \end{pmatrix} 
\end{xalignat}
with respect to the $\{-1,0,1\} = \{-,3,+\}$ basis. The matrix
representations of the vector generator is proportional to the epsilon
tensor,
\begin{equation}
\label{eq:AppRep1}
  \rho^1(J_A)^B{}_C = \varepsilon_A{}^B{}_C \,.
\end{equation}

\subsection{Representations of the $q$-Lorentz Algebra}

The representation maps for the $D^{(j_1,j_2)}$ representations of the
$q$-Lorentz algebra, $\slC$, are composed of the representation maps
of $\suq$ according to $\rho^{(j_1,j_2)} :=\rho^{j_1}\otimes
\rho^{j_2}$.  Particularly simple are the chiral representations
$D^{(j,0)}$ and $D^{(0,j)}$. For any rotations $l\in\suq$ and for the
boosts as defined in Eq.~\eqref{eq:Appboostdef1} we get
\begin{equation}
\begin{gathered}
  \rho^{(j,0)}(l) = \rho^j(l) = \rho^{(0,j)}(l) \\
  \rho^{(j,0)}(B^a{}_b) = \rho^j\bigl((L_-^\frac{1}{2})^a{}_b\bigr) \,,\qquad 
  \rho^{(0,j)}(B^a{}_b) = \rho^j\bigl((L_+^\frac{1}{2})^a{}_b\bigr) \,.
\end{gathered}
\end{equation}
If we denote the basis of a $D^{(\frac{1}{2},\frac{1}{2})}$
representation as in Eq.~\eqref{eq:MinkowskiStar} by
\begin{equation}
  \Base{ab}=
  \begin{pmatrix}\Base{--}&\Base{-+}\\ \Base{+-}&\Base{++}\end{pmatrix}
  =:\begin{pmatrix}A&B\\C&D\end{pmatrix} \,,
\end{equation}
we get for the action
\begin{xalignat}{2}
  E\otimes 1 \tr \begin{pmatrix} A & B \\ C & D \end{pmatrix} &=
  \begin{pmatrix} C & D \\ 0 & 0 \end{pmatrix}, &
  1\otimes E \tr \begin{pmatrix} A & B \\ C & D \end{pmatrix} &=
  \begin{pmatrix} B & 0 \\ D & 0 \end{pmatrix} \notag \\
  F\otimes 1 \tr \begin{pmatrix} A & B \\ C & D \end{pmatrix} &=
  \begin{pmatrix} 0 & 0 \\ A & B \end{pmatrix}, &
  1\otimes F \tr \begin{pmatrix} A & B \\ C & D \end{pmatrix} &=
  \begin{pmatrix} 0 & A \\ 0 & C \end{pmatrix} \\
  K\otimes 1 \tr \begin{pmatrix} A & B \\ C & D \end{pmatrix} &=
  \begin{pmatrix} q^{-1}A & q^{-1}B \\ qC & qD \end{pmatrix}, &
  1\otimes K \tr \begin{pmatrix} A & B \\ C & D \end{pmatrix} &=
  \begin{pmatrix} q^{-1}A & qB \\ q^{-1}C & qD \end{pmatrix} \,.\notag
\end{xalignat}
For the boost generators~\eqref{eq:Appboostdef1} this means in particular
\begin{subequations}
\begin{xalignat}{2}
  \begin{pmatrix} a & b \\ c & d \end{pmatrix} \tr A &=
  \begin{pmatrix} A & q^{-1}\lambda B \\ 0 & A \end{pmatrix}, &
  \begin{pmatrix} a & b \\ c & d \end{pmatrix} \tr B &=
  \begin{pmatrix} q^{-1} B & 0\\ 0 & qB \end{pmatrix} \\
  \begin{pmatrix} a & b \\ c & d \end{pmatrix} \tr C &=
  \begin{pmatrix} qC & \lambda D \\ -q\lambda A &
     q^{-1}C - \lambda^2 B\end{pmatrix}, &
  \begin{pmatrix} a & b \\ c & d \end{pmatrix} \tr D &=
  \begin{pmatrix} D & 0 \\ -\lambda B & D \end{pmatrix} 
\end{xalignat}
\end{subequations}
In terms of the 4-vector basis $E_\mu$ of
$D^{(\frac{1}{2},\frac{1}{2})}$, defined as in
Eq.~\eqref{eq:FourVectorDef} by
\begin{equation}
\begin{aligned}
  E_0 &= q^{-1}[2]^{-1}(q^{\frac{1}{2}}C - q^{-\frac{1}{2}}B) \\
  E_- &= [2]^{-\frac{1}{2}} A\\
  E_+ &= [2]^{-\frac{1}{2}} D\\
  E_3 &= [2]^{-1}(q^{-\frac{1}{2}}C + q^{\frac{1}{2}}B) \,.
\end{aligned}
\end{equation}
the action becomes
\begin{equation}
\label{eq:AppBoostAct}
\begin{aligned}
  \begin{pmatrix} a & b \\ c & d \end{pmatrix} \tr E_0 &=
    \begin{pmatrix}
      [2]^{-1}\left(\frac{[4]}{[2]}E_0 + q^{-1}\lambda E_3 \right)  &
      q^{-\frac{1}{2}}\lambda [2]^{-\frac{1}{2}} E_+ \\
      -q^{\frac{1}{2}}\lambda [2]^{-\frac{1}{2}} E_- &
      [2]^{-1}\left( \frac{[4]}{[2]} E_0 - q \lambda E_3 \right)
    \end{pmatrix} \\
  \begin{pmatrix} a & b \\ c & d \end{pmatrix} \tr E_- &=
    \begin{pmatrix} E_- \quad &
      q^{-\frac{1}{2}}\lambda [2]^{-\frac{1}{2}}(E_3 - E_0) \\
      0 & E_-
    \end{pmatrix} \\
  \begin{pmatrix} a & b \\ c & d \end{pmatrix} \tr E_+ &=
    \begin{pmatrix} E_+ & 0 \\
      -q^{\frac{1}{2}}\lambda [2]^{-\frac{1}{2}}(E_3 - E_0) & \quad E_+
    \end{pmatrix}\\
  \begin{pmatrix} a & b \\ c & d \end{pmatrix} \tr E_3 &=
    \begin{pmatrix} [2]^{-1}(2E_3 + q\lambda E_0)&
      q^{-\frac{1}{2}}\lambda [2]^{-\frac{1}{2}} E_+ \\
      -q^{\frac{1}{2}}\lambda [2]^{-\frac{1}{2}} E_- & \quad
      [2]^{-1}(2E_3 - q^{-1} \lambda E_0)
    \end{pmatrix}  
\end{aligned}
\end{equation}
Now we can calculate the 4-vector matrix representation $\Lambda$
defined by
\begin{equation}
  h\tr E_\mu = E_{\mu'}\,\Lambda(h)^{\mu'}{}_\mu
\end{equation}
for all $q$-Lorentz transformations $h$. For the rotations $l\in\suq$
we get by construction of the 4-vector basis
\begin{equation}
\label{eq:AppFourVector1}
  \Lambda(l) =
  \begin{pmatrix} \rho^0(l) & 0 \\ 0 & \rho^1(l)\end{pmatrix} \,.
\end{equation}
For the boost we calculate
\begin{equation}
\label{eq:AppFourVector2}
\begin{gathered}
  \Lambda(a) =
  \begin{pmatrix}
  [4][2]^{-2} & 0 & 0 & q\lambda [2]^{-1} \\
  0 & 1 & 0 & 0 \\
  0 & 0 & 1 & 0 \\
  q^{-1}\lambda [2]^{-1} & 0 & 0 & 2 [2]^{-1} \\
  \end{pmatrix},\quad
  \Lambda(b) = q^{-\frac{1}{2}}\lambda[2]^{-\frac{1}{2}}
  \begin{pmatrix}
  0 & -1 & 0 & 0 \\
  0 & 0 & 0 & 0 \\
  1 & 0 & 0 & 1 \\
  0 & 1 & 0 & 0 \\
  \end{pmatrix} \\
  \Lambda(c) = -q^{\frac{1}{2}}\lambda[2]^{-\frac{1}{2}}
  \begin{pmatrix}
  0 & 0 & -1 & 0 \\
  1 & 0 & 0 & 1 \\
  0 & 0 & 0 & 0 \\
  0 & 0 & 1 & 0 \\
  \end{pmatrix},\quad
  \Lambda(d) =
  \begin{pmatrix}
  [4][2]^{-2} & 0 & 0 & -q^{-1}\lambda [2]^{-1} \\
  0 & 1 & 0 & 0 \\
  0 & 0 & 1 & 0 \\
  -q \lambda [2]^{-1} & 0 & 0 & 2 [2]^{-1}
  \end{pmatrix} \,,
\end{gathered}
\end{equation}
with respect to the $\{0,-,+,3\}$ basis.

\section{$\R$-matrices}
\label{sec:AppR}

For a Hopf algebra $H$ a universal $\R$-matrix is an invertible
element $\R\in H\otimes H$, which we will also write in a Sweedler
like notation as $\R := \R_{[1]}\otimes \R_{[2]}$, with
\begin{equation}
\label{eq:Rmatrix1}
\begin{gathered}
  (\tau\circ\Delta)(h) = \R\,\Delta(h)\,\R^{-1}\\
  (\Delta\otimes\id)(\R)=\R_{13}\R_{23}\,,\qquad
  (\id\otimes\Delta)(\R)=\R_{13}\R_{12} \,,
\end{gathered}
\end{equation}
where the indices indicate the position of the tensor factors,
$\R_{13} := \R_{[1]}\otimes 1 \otimes \R_{[2]}$ etc. If there is a
$*$-structure on $H$ the $\R$-matrix is said to be real if
$\R^{*\otimes *} = \R_{21}$ and anti-real if $\R^{*\otimes *} = \R^{-1}$.
There are some useful properties of $\R$ that can be deduced from
Eqs.~\eqref{eq:Rmatrix1}:
\begin{equation}
\label{eq:Rmatrix2}
\begin{gathered}
  \R_{12}\R_{13}\R_{23} = \R_{23}\R_{13}\R_{12}\,,\quad
  (\varepsilon\otimes\id)(\R)=1 \,\quad
  (\id\otimes\varepsilon)(\R)=1 \\
  (S\otimes\id)(\R)= \R^{-1}\,,\quad
  (\id\otimes S)(\R^{-1})= \R\,,\quad
  (S\otimes S)(\R)= \R \,. 
\end{gathered}
\end{equation}

\subsection{The $\R$-Matrix of $\suq$}
\label{sec:AppR1}

There is a universal $\R$-matrix for $\suq$,
\begin{equation}
\label{eq:Rmatrix3}
  \R = q^{(H\otimes H)/2} \sum_{n=0}^{\infty} R_n(q) (E^n \otimes F^n) \,
\end{equation}
which is not an element $\suq\otimes\suq$ proper, since it is
described as an infinite power series. For our purposes this does not
raise serious problems. This $\R$-matrix is real. For representations
$\rho^j$, $\rho^{j'}$ of $\suq$ we can define $R$-matrices and a
variant, the $\hat{R}$-matrices, by
\begin{xalignat}{2}
  R^{(j,j')} &:= (\rho^j\otimes \rho^{j'})(\R) \,,&
  (\hat{R}^{(j,j')})^{ab}{}_{cd} &:= (R^{(j,j')})^{ba}{}_{cd} \,.
\end{xalignat}
Traditionally, the $R$-matrices are normalized differently. We will use
\begin{xalignat}{2}
\label{eq:AppR4}
  R_{\mathrm{su}_2} &:= q^{\frac{1}{2}} R^{(\frac{1}{2},\frac{1}{2})}\,, &
  R_{\mathrm{so}_3} &:= q^{-2} R^{(1,1)} \,.
\end{xalignat}
Explicitly, we get
\begin{equation}
   (R_{\mathrm{su}_2})^{ab}{}_{cd} =
  \begin{pmatrix}
    q & 0 & 0 & 0 \\
    0 & 1 & 0 & 0 \\
    0 & \lambda & 1 & 0 \\
    0 & 0 & 0 & q 
  \end{pmatrix},
\end{equation}
with respect to the basis $\{--,-+,+-,++\}$, and
\begin{equation}
\label{eq:AppR6}
\begin{aligned}
 (R_{\mathrm{so}_3})^{AB}{}_{CD}
   &= \delta^B_C \delta^A_D - q^{-3}\lambda\, g^{BA}g_{CD}
   - q^{-2} \varepsilon^{BAX}\varepsilon_{DCX} \\
 (R^{-1}_{\mathrm{so}_3})^{AB}{}_{CD}
   &= \delta^A_D \delta^B_C - q^{3}\lambda\, g^{AB}g_{DC}
   - q^{2} \varepsilon^{ABX}\varepsilon_{CDX}  \,.
\end{aligned}
\end{equation}
This means that we have a projector decomposition
\begin{equation}
  \hat{R}_{\mathrm{so}_3}
  = 1 -q^{-3}\lambda[3]^{-1}\Proj_0 -q^{-2}[4][2]^{-1} \Proj_1
  = -q^{-6}\Proj_0 -q^{-4}\Proj_1 + \Proj_3 \,.
\end{equation}
Applying a representation to one half of the $\R$-matrix only leads to
the definition of the $L$-matrices
\begin{xalignat}{2}
  (L^j_+)^a{}_b &:= \R_{[1]} \rho^j(\R_{[2]})^a{}_b \,,&
  (L^j_-)^a{}_b &:= \rho^j(\R^{-1}_{[1]})^a{}_b
    \,\R^{-1}_{[2]} \,.
\end{xalignat}
We calculate the $L$-matrices for $j=\frac{1}{2}$ and $j=1$,
explicitly.
\begin{xalignat}{2}
  L_+^\frac{1}{2} &= \begin{pmatrix}
    K^{-\frac{1}{2}} &
    q^{-\frac{1}{2}}\lambda K^{-\frac{1}{2}}E \\
              0      & K^{\frac{1}{2}}   \end{pmatrix}, &
  L_-^\frac{1}{2} &= \begin{pmatrix}
   K^{\frac{1}{2}} & 0 \\ -q^{\frac{1}{2}}\lambda FK^{\frac{1}{2}}
   & K^{-\frac{1}{2}} \end{pmatrix} \\
\label{eq:AppR5}
  L_+^1 &= \begin{pmatrix}
    K^{-1} & \lambda[2]^{\frac{1}{2}} K^{-1}E & \lambda^2 K^{-1}E^2\\
     0 & 1 & q^{-1}\lambda[2]^{\frac{1}{2}} E \\
     0 & 0 & K 
   \end{pmatrix}, &
  L_-^1 &= \begin{pmatrix}
     K & 0 & 0\\
     -\lambda[2]^{\frac{1}{2}} FK & 1 & 0 \\
     \lambda^2 F^2 K & -q\lambda[2]^{\frac{1}{2}}F & K^{-1}
     \end{pmatrix}
\end{xalignat}
These results are being used in Eq.~\eqref{eq:boostdef} to calculate
the boost generators defined as
\begin{equation}
\label{eq:Appboostdef1}
  B^a{}_c :=
  \bigl(L_-^\frac{1}{2}\bigr)^a{}_b \otimes
  \bigl(L_+^\frac{1}{2} \bigr)^b{}_c
  =: \begin{pmatrix} a & b \\ c & d \end{pmatrix}\,,
\end{equation}
which yields
\begin{equation}
\label{eq:Appboostdef2}
  \begin{pmatrix} a & b \\ c & d \end{pmatrix} =
  \begin{pmatrix}
    K^{\frac{1}{2}}\otimes K^{-\frac{1}{2}} &
    q^{-\frac{1}{2}}\lambda K^{\frac{1}{2}}
    \otimes K^{-\frac{1}{2}} E \\
    -q^{\frac{1}{2}}\lambda F K^{\frac{1}{2}}
    \otimes K^{-\frac{1}{2}} & \quad
    K^{-\frac{1}{2}}\otimes K^{\frac{1}{2}}
    - \lambda^2 F K^{\frac{1}{2}}\otimes K^{-\frac{1}{2}} E
  \end{pmatrix}.
\end{equation}

\subsection{The $\R$-Matrices of the $q$-Lorentz Algebra}
\label{sec:AppR2}

There are two universal $\R$-matrices of the $q$-Lorentz algebra,
which are composed of the $\R$-matrix of $\slq$ according to
\begin{xalignat}{2}
  \RI &= \R^{-1}_{41}\R^{-1}_{31}\R_{24}\R_{23}\,, &
  \RII &= \R^{-1}_{41}\R_{13}\R_{24}\R_{23} \,.
\end{xalignat}
$\RI$ is anti-real while $\RII$ is real. Their vector representations
are normalized as
\begin{xalignat}{2}
\label{eq:AppR3}
  R_\mathrm{I} &:= (\Lambda \otimes \Lambda)(\RI) \,,&
  R_\mathrm{II} &:= q (\Lambda \otimes \Lambda)(\RII) \,,
\end{xalignat}
where $\Lambda$ is the 4-vector representation map of the $q$-Lorentz
algebra. These matrices can be decomposed into projectors
\begin{equation}
\begin{aligned}
  \hat{R}_\mathrm{I} &= \Proj_{(0,0)} - q^{-2} \Proj_{(1,0)}
  - q^{2} \Proj_{(0,1)} + \Proj_{(1,1)} \\
  \hat{R}_\mathrm{II} &= q^{-2}\Proj_{(0,0)} - \Proj_{(1,0)}
  - \Proj_{(0,1)} + q^{2} \Proj_{(1,1)} \,.
\end{aligned}
\end{equation}
The $L_+$-matrix of $\RI$ has a simple form:
\begin{equation}
\label{eq:AppR2}
  \bigl(L_{\mathrm{I}+}^\Lambda\bigr)^a{}_b :=
  \RI{}_{[1]} \,\Lambda(\RI{}_{[2]})^a{}_b =
  \begin{pmatrix} 1 & 0 \\ 0 & t^A{}_B \end{pmatrix} \,,
\end{equation}
where $t^A{}_B$ is the vector corepresentation matrix of $\SUq^\op$,
\begin{equation}
  t = \begin{pmatrix}
    a^2 & q^{\frac{1}{2}}[2]^{\frac{1}{2}} ab & b^2\\
    q^{\frac{1}{2}}[2]^{\frac{1}{2}} ac & (1+[2]bc) &
    q^{\frac{1}{2}}[2]^{\frac{1}{2}} bd \\
    c^2 & q^{\frac{1}{2}}[2]^{\frac{1}{2}} cd & d^2
  \end{pmatrix}
\end{equation}
with respect to the basis $\{-1,0,1\} = \{-,3,+\}$. For chiral
representations we get
\begin{xalignat}{2}
\label{eq:AppR1}
  \rho^{(j,0)}(t^A{}_B) &= \rho^j\bigl((L_-^1)^A{}_B\bigr)\,, & 
  \rho^{(0,j)}(t^A{}_B) &= \rho^j\bigl((L_+^1)^A{}_B\bigr) \,.
\end{xalignat}

\end{document}